\documentclass[a4paper, 11pt]{amsart}
\usepackage{amsmath}
\usepackage{amssymb}
\usepackage{epsfig}
\usepackage{eepic}
\usepackage{latexsym
}
\usepackage{color}

\usepackage{subcaption}
\usepackage{verbatim}
\usepackage{graphicx}
\usepackage{enumerate}
\usepackage{xspace}
\usepackage[ruled]{algorithm}
\usepackage{booktabs}
\usepackage{tabularx}
\usepackage{natbib}
\usepackage{times}
\usepackage{amsmath,amssymb,amsxtra}

\usepackage{hyperref}
\usepackage{commath}
\usepackage{tikz}

\usepackage{amsmath,amssymb,amsxtra}

\usepackage{amsthm}

\usepackage{dsfont}
\newtheorem{theorem}{Theorem}[section]

\newtheorem{definition}[theorem]{Definition}

\newtheorem{remark}{Remark}[section]

\newcommand{\IR}{\mathbb{R}}

\def \ve{\varepsilon}
\usepackage[latin1]{inputenc}

\usepackage[top=2 cm,bottom=2 cm,left=2. cm, right=2. cm]{geometry}
\def \O{\Omega}

\def\XXint#1#2#3{{\setbox0=\hbox{$#1{#2#3}{\int}$}
\vcenter{\hbox{$#2#3$}}\kern-.87\wd0}}

\def\XXiint#1#2#3{{\setbox0=\hbox{$#1{#2#3}{\int}$}
\vcenter{\hbox{$#2#3$}}\kern-1.05\wd0}}

\def\XXintt#1#2#3{{\setbox0=\hbox{$#1{#2#3}{\int}$}
\vcenter{\hbox{$#2#3$}}\kern-.72\wd0}}

\def\Xinttt#1{\mathchoice
{\XXinttt\displaystyle\textstyle{#1}}%
{\XXinttt\textstyle\scriptstyle{#1}}%
{\XXinttt\scriptstyle\scriptscriptstyle{#1}}%
{\XXinttt\scriptscriptstyle\scriptscriptstyle{#1}}%
\!\int}
\def\XXinttt#1#2#3{{\setbox0=\hbox{$#1{#2#3}{\int}$}
\vcenter{\hbox{$#2#3$}}\kern-.52\wd0}}

\def\XXintttr#1#2#3{{\setbox0=\hbox{$#1{#2#3}{\int}$}
\vcenter{\hbox{$#2#3$}}\kern-.6\wd0}}

\def\XXintttt#1#2#3{{\setbox0=\hbox{$#1{#2#3}{\int}$}
\vcenter{\hbox{$#2#3$}}\kern-.78\wd0}}

\def\sqr#1#2{{\vcenter{\vbox{\hrule height.#2pt\hbox{\vrule width.#2pt height#1pt \kern#1pt\vrule width.#2pt}\hrule height.#2pt}}}}

\def\ddashinttt{\Xinttt-}

\renewcommand{\vec}[1]{\ensuremath{\boldsymbol{#1}}}
\newcommand{\myall}{\ensuremath{\quad \forall}}

\newcommand{\Reals}{\ensuremath{{\mathbb{R}}}}

\newcommand{\diff}{\ensuremath{{\operatorname{d}}}}

\renewcommand{\O}{\ensuremath{{\Omega}}}

\newcommand{\Hil}[1]{\ensuremath{\operatorname{H}^{#1}}}

 \usepackage{transparent}

\newcommand{\A}{\ensuremath{\vec{\mathcal{A}}}}

\usepackage{upgreek}

\usepackage{color}
\definecolor{MyGreen}{rgb} {0.05,0.4,0.05}
\definecolor{RedViolet}{rgb} {0.1,0.1,0.75}

\usepackage{tabularx}

\newcommand{\Th}{\ensuremath{{\mathcal{S}_{h,\O}}}}

\newcommand{\She}{\ensuremath{{\mathcal{S}_{h,e}}}}
\newcommand{\Vho}{\ensuremath{{\mathbb{V}_{h,\Omega}}}}
\newcommand{\Vhe}{\ensuremath{{\mathbb{V}_{h,e}^\#}}}

\newcommand{\Ltwop}[3]{\ensuremath{\left\langle#1,#2\right\rangle}_{#3}}

\newcommand{\edit}[1]{{#1}}

\begin{document}

\title[Multiscale modelling of cancer invasion]{Multiscale modelling, analysis and simulation of  cancer invasion mediated by bound and soluble enzymes} 
\author[M.~Ptashnyk, C.~Venkataraman]{{Mariya Ptashnyk, Chandrasekhar Venkataraman}\medskip\\
  Department of Mathematics,  Heriot-Watt University, The Maxwell Institute for Mathematical Sciences, Edinburgh, Scotland,  UK\\
  Department of Mathematics, University of Sussex, UK  \\
   }
\date{}
\maketitle

\begin{abstract}
We formulate a cell-scale model for the degradation of the extra-cellular matrix by membrane-bound and soluble matrix degrading enzymes produced by cancer cells. Based on the microscopic model and using tools from the theory of homogenisation we propose a macroscopic model for cancer cell invasion into the extra-cellular matrix mediated by bound and soluble matrix degrading enzymes.
For suitable and biologically relevant initial data we prove the macroscopic model is well-posed. We propose a finite element method for the numerical approximation of the macroscopic model and report on simulation results illustrating the role of the bound and soluble enzymes in cancer invasion processes.
\end{abstract}

\section{Introduction} 

Invasion is one of the hallmarks of cancer \citep{HanahanWeinberg2000, HanahanWeinberg2011, Hanahan2022}. It is a complex process involving numerous interactions between cancer cells and the extracellular matrix (cf.~the tumour microenvironment) facilitated by matrix degrading enzymes. Along with active cell migration, both individual and collective, and increased/excessive proliferation, these processes enable the local spread of cancer cells into the surrounding tissue. Any encounter with blood or lymphatic vessels (cf.~tumour-induced angiogenesis, lymph-angiogenesis) in the tumour microenvironment offers the opportunity of intravasation, which together with subsequent extravasation, determines the dissemination of the cancer to secondary locations in the host, i.e., metastasis or metastatic spread. Further details of the invasion-metastasis process, and also extensive biological/clinical references, can be found in, e.g.,~\cite{Fidler2003, friedl2003tumour, HanahanWeinberg2000, HanahanWeinberg2011, klein2008, Talmadge2010}. 

Mathematical modelling of cancer invasion can be traced back to the model formulated in~\cite{Gatenby1995} using an approach from theoretical ecology and a modified predator-prey system, which focussed on the interactions between cancer cells and normal, healthy cells competing for space and other resources within a tissue. Spatial models of cancer invasion were then developed using systems of reaction-diffusion-taxis equations \citep{Orme1996b, Gatenby1996, PerumpananiSherrattNorburyEtAl1996}, focussing on  haptotaxis as the key cell migration mechanism. A hybrid discrete-continuum approach was developed by \cite{anderson2000} which enabled the depth of penetration of the ECM by individual cancer cells to be studied. The reaction-diffusion-taxis equations formulated in \cite{anderson2000} were analysed in \cite{ptashnyk2010}. The effects of cell-cell and cell-matrix adhesion  were incorporated in~\cite{gerisch2008} using a non-local model.  A number of studies considered the derivation of mixture or multi-phase models for tumour growth~\cite{byrne2003two,franks2003interactions}. The phase-field modelling approach to describe tumour growth and invasion was proposed in~\cite{cristini2009nonlinear,fritz2019, fritz2023} where the Cahn-Hilliard equation was used to model the invasion dynamics. An analysis of a diffusive interface model for tumour growth describing evolutions driven by long-range interactions was considered in \cite{scarpa2021}, while the Cahn-Hilliard-Darcy model for tumour growth including nutrient diffusion, chemotaxis, active transport, adhesion, apoptosis and proliferation was derived in \cite{garcke2016}  with extensions to account for viscoelastic effects  considered in \cite{garcke2022viscoelastic}. A recent review of mathematical models of cancer invasion can be found in~\cite{sfakianakischaplain2021}.

A crucial part of the invasive/metastatic process is the ability of the cancer cells to degrade the surrounding tissue or extracellular matrix (ECM). The ECM is a complex mixture of macromolecules, some of which, like the collagens, are believed to play a structural role and others, such as laminin, fibronectin and vitronectin, are important for cell adhesion, spreading and motility. All of these macromolecules are  bound within the tissue, i.e., they are non-diffusible. The extracellular matrix can also sequester growth factors and itself be degraded to release fragments which can have growth-promoting activity. Thus, while extracellular matrix may have to be physically removed in order to allow a tumour to spread or intra- or extravasate, its degradation may in addition have biological effects on tumour cells. Degradation of the ECM is achieved through the proteolytic activity of various matrix degrading enzymes (MDE). The experiments of \cite{sabeh2009protease} demonstrated that there are two important types of MDE - {\em membrane-bound} and {\em diffusible} metalloproteinases - involved in cancer invasion and their study focussed on the membrane-bound membrane-type-1 matrix metalloproteinase (MT1-MMP) and the diffusible (soluble) matrix metalloproteinase-2 (MMP-2), and their interactions with the extracellular matrix.

Following the work in~\cite{deakin2013mathematical}, which was  based on the experimental results in~\cite{sabeh2009protease}, in this paper we consider models for cancer cell invasion  mediated by enzyme based degradation of the ECM. Specifically our focus is on ECM degradation via both soluble and membrane-bound matrix metalloproteinases (MMPs). We propose a microscopic model which  accounts for matrix degradation via these two pathways and formulate an effective macroscopic model valid in the limit where the size of the tissue is large relative to the size of the individual cells.
Assuming a periodic or locally periodic distribution of cancer cells,  the rigorous derivation of the macroscopic model can be obtained using  homogenization methods  namely, two-scale convergence and the unfolding operator, see e.g.~\cite{allaire1992, cioranescu2018, nguetseng1989, ptashnyk2013, ptashnyk2015}, using similar methods as developed in \cite{marciniak2008derivation, ptashnyk2020multiscale}.
It is however extremely challenging, and beyond the scope of this work, to incorporate tumour invasion via movement and proliferation of the cancer cells in the  microscopic model such that effective macroscopic equations can be rigorously derived from the microscopic description of the processes. We remark that a major obstacle in this direction is the lack of a sufficient theory for continuum modelling of the dynamics of cellular processes during proliferation.  We  instead propose a phenomenological approach to derive an effective macroscopic model for invasion in which microscopic features, specifically degradation of the ECM by membrane-bound MMPs, can be accounted for. Furthermore, we account for the influence of the microstructure on the transport of the soluble MMPs through the ECM by deriving an effective diffusion coefficient  using   homogenization methods.  We note that the framework and methodology we propose for computing the effective diffusivity is applicable generally to multiphase models of tumour growth  such as those mentioned previously, as well as other applications involving heterogeneous mixtures. Our numerical results demonstrate  qualitative differences between the computed effective diffusion tensors  in two versus three dimensions. We show existence of nonnegative and bounded solutions of the macroscopic model using a fixed-point argument together with the Galerkin method and the method of positive invariant regions. The model consists of a coupled system of ordinary differential equations and an anisotropic  degenerate parabolic equation. The variable doubling method can be applied to show well-posedness  for anisotropic  degenerate parabolic equations with diffusion coefficients in $W^{1,\infty}(\Omega)$, see e.g.~\cite{chen2003, chen2005}.  However, in our case the dependence of the diffusion coefficient on the ECM density, which satisfies an ordinary differential equation, does not imply the required regularity. To prove the well-posedness results we use the fact that for strictly positive initial ECM density and for any finite time interval the ECM density is uniformly bounded from below by a positive constant.
Numerical simulations of the macroscopic model illustrate the crucial role bound MMPs appear to play in cancer invasion, in agreement with previous results~\cite{deakin2013mathematical, sabeh2009protease}, as well as the fact that the action of bound MMPs can lead to spatial heterogeneity in the invasive front   whilst the action of the soluble MMPs seems to generate more radially symmetric invasive profiles.  
 
 In summary, the key results of the work are
 \begin{itemize}
 \item Development of a framework in which detailed cell-scale models for ECM degradation can be  incorporated in a macroscopic model for cancer cell invasion of the ECM that can be simulated at tissue scales.
 \item A methodology for computing the effective diffusivity of soluble molecules in heterogeneous medium representing cancer cells and ECM with varying volume fractions.  
 \item The analysis of a coupled system of semilinear ODE-PDEs that serves as an effective model for cancer cell invasion of the ECM.
 \item Numerical simulation results of invasive processes illustrating qualitative features of the modelling such as the crucial role bound MMPs appear to play in determining both the speed of invasion and spatial heterogeneity in the invasive front.
 \end{itemize} 

The paper is organised as follows. In section~\ref{sec:micro} we formulate a microscopic model for ECM degradation through the action of bound and soluble MMPs which consists of a system of coupled bulk-surface equations. In section~\ref{sec:macro} we propose effective macroscopic models for cancer cell invasion of the ECM. In section~\ref{sec:analysis} we prove the well posedness of the macroscopic model, i.e., the coupled semilinear ODE-PDE system.  
In section~\ref{sec:FEM} we formulate a finite element scheme for the approximation of the macroscopic model. In section~\ref{sec:numerics} we discuss the parameterisation of the model based, where possible, on experimentally measured parameter values,  present numerical simulation results for the effective diffusivity of soluble MMPs 
for different volume fractions of the ECM in $2{\rm d}$ and $3{\rm d}$, and   report on the results
of numerical simulations of the macroscopic model that illustrate various aspects of cancer invasion. 
Finally, in section~\ref{sec:conc} we discuss our results as well as potential directions for future work.

\section{Microscopic description of cancer cell and ECM interactions}\label{sec:micro}

\begin{figure}[htbp]
  \includegraphics[trim = 0mm 0mm 0mm 0mm,  clip, width=0.73\linewidth]{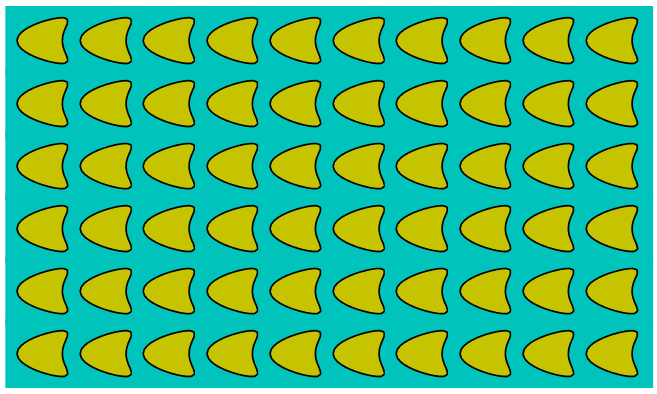}
  \put(-348,68){ $\Omega^\ve \hspace{-0.1 cm}  $}
   \put(-338,87){ $\Omega^\ve_i$}
\put(-350.4,34.){ $\Gamma^\ve$}
\put(-342.4,36.5){ $\nearrow  $}
\caption{A sketch of the geometry consisting of  cancer cells surrounded by the ECM. The domains $\Omega_i^{\ve}, \Gamma^\ve$ and $\Omega^\ve$ corresponding to the cancer cells, cancer cell membranes and ECM are shaded yellow, black and blue respectively. }
  \label{fig:micro_domain}
  \end{figure}
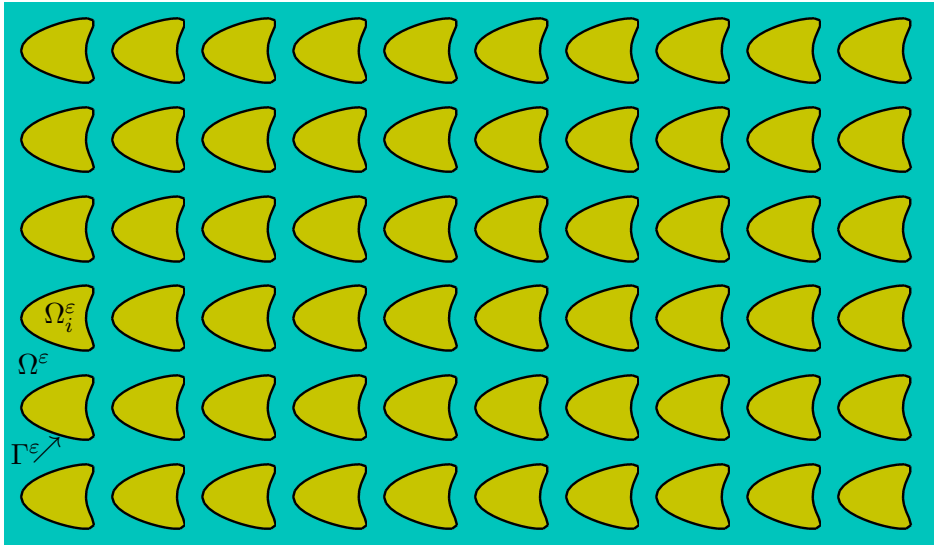

To formulate the microscopic model, we consider processes at the level of single cells where ECM degradation occurs due to the action of soluble and membrane bound MMPs. 
  We consider a convex Lipschitz, or $C^{1,1}$, domain $\Omega\subset \IR^d$, with $d=2,3$, representing a part of a biological tissue and  assume a  time independent geometry and a periodic distribution of cells of the same shape.  To describe the microscopic structure of the tissue, we consider a `reference domain' $Y=(0,1)^d$,  and the subdomains  $\overline Y_i\subset Y$ and $Y_e = Y \setminus \overline Y_i$, together with  the boundary  $\Gamma = \partial Y_i$.   The domain occupied  by the cells  is given by $\Omega_i ^\ve= \bigcup_{\xi \in \Xi^\ve} \ve(Y_i + \xi) $,
where $\Xi^\ve = \{ \xi \in \mathbb Z^d, \; \; \ve(\overline Y_i + \xi) \subset \Omega \}$, the nonnegative parameter $\ve$ represents the ratio between the diameter of a cancer cell and the tissue.  The extracellular space is denoted by  $\Omega^\ve= \Omega\setminus \overline \Omega_i^\ve$.  The surfaces that describe cell membranes are denoted by $\Gamma^\ve= \bigcup_{\xi \in \Xi^\ve} \ve(\Gamma + \xi)$ and  we assume they are smooth surfaces without boundary, see Figure~\ref{fig:micro_domain} for a sketch of the geometry. As we are interested in tissues containing many cells we have $\ve \ll 1$.

We let $c^\ve_s:\Omega^\ve\times[0,T]\to\Reals$ denote the concentration of soluble MMPs that are produced by the cancer cells and diffuse in the ECM, $c_b^\ve:\Gamma^\ve\times[0,T]\to\Reals$ denotes the surface density of bound MMPs on  the cancer cell membranes and we let $m^\ve: \Omega^\ve \times[0,T]\to\Reals$  and  \edit{$m^\ve_\Gamma: \Gamma^\ve \times[0,T]\to\Reals$  denote the ECM density in the extracellular domain and at the cell surface, respectivelly}.  

For $c_s^\ve$ we propose the following dynamics
\begin{equation}\label{micro_cs}
\begin{aligned}
& \partial_t c_s^\ve   - \nabla\cdot(\bar{D}_s(x)\nabla c_s^\ve) = - g_s(c_s^\ve, m^\ve) - \beta_s c_s^\ve && \text{ in } (0,T] \times  \Omega^\ve,\\
&\bar{D}_s(x)\nabla c_s^\ve \cdot \nu = \ve \kappa_s \edit{f_s(c_s^\ve, m^\ve_\Gamma) } && \text{ on } (0,T] \times  \Gamma^\ve, \\
& \bar{D}_s(x)\nabla c_s^\ve \cdot \nu = 0 && \text{ on } (0,T] \times \partial \Omega,  \\
& c_s^\ve (0) = c_{s, \rm in} && \text{ in }  \Omega^\ve,
\end{aligned} 
\end{equation}
where $\nu$ is the normal pointing outwards of $\Omega^\ve$, $\bar D_s$ corresponds to the diffusivity of soluble MMPs in the ECM, $\kappa_s f_s(c_s^\ve, m_\Gamma^\ve)$ describes the production \edit{ and/or degradation} of soluble enzymes \edit{ by cells}, $g_s(c_s^\ve, m^\ve)$ models the use of enzymes in degradation of ECM process and $\beta_s$ is the natural decay rate.
 Although, for simplicity, we consider $\bar D_s$ to depend only on the space variable,  it could depend on the ECM density to account for effects due to changes in the pore size of the ECM.
For $c_b^\ve$ we consider 
\begin{equation}\label{micro_cb}
\begin{aligned}
& \partial_t c_b^\ve - \ve^2 \nabla_\Gamma \cdot( D_b \nabla_\Gamma c_b^\ve)  =  \kappa_b \edit{ f_b(c_b^\ve, m^\ve_\Gamma) }- g_b(c_b^\ve, m^\ve_\Gamma) - \beta_b c_b^\ve && \text{ on } (0,T]\times  \Gamma^\ve,\\
& c_b^\ve (0) = c_{b, \rm in}^\ve && \text{ on } \Gamma^\ve,
\end{aligned} 
\end{equation}
where $D_b$ is the surface diffusivity of $c_b^\ve$, $\kappa_b f_b(c_b^\ve, m^\ve_\Gamma)$ describes production of bound enzymes through activation of the corresponding receptors due to  interactions between cancer cells and  ECM, $\beta_b$ is the natural decay rate, $ g_b(c_b^\ve, m^\ve_\Gamma)$ corresponds to \edit{the use of the cancer cell-bound enzymes in the degradation of ECM}, and \edit{$c_{b, \rm in}^\ve(x) = c_{b, \rm in}(x/\ve)$ for  $c_{b, \rm in}\in L^2(\Gamma)$,  extended $Y$-periodically to $\mathbb R^d$.}
Finally for the ECM densities $m^\ve$ and $m^\ve_\Gamma$ we consider the following equations modelling  the degradation of EMC due to
interactions of bound and free MMPs with the ECM
\begin{equation}\label{micro_m}
\begin{aligned}
&\partial_t m^\ve =  - g_s(c_s^\ve, m^\ve) && \text{ in } (0,T]\times  \Omega^\ve, \\
&\edit{\partial_t m_\Gamma^\ve = -  g_b(c_b^\ve, m^\ve_\Gamma) } && \text{ on } (0,T]\times  \Gamma^\ve, \\
& m^\ve (0) = m_{\rm in} \; \text{ in }   \Omega^\ve, \quad \edit{ m^\ve_\Gamma (0) = m_{\Gamma,\rm in}^\ve} && \text{ on }   \Gamma^\ve,
\end{aligned} 
\end{equation}
\edit{where $m_{\Gamma,\rm in}^\ve(x) = m_{\Gamma,\rm in}(x/\ve)$ for $m_{\Gamma,\rm in} \in L^2(\Gamma)$,  extended $Y$-periodically to $\mathbb R^d$. }
The $\ve$-scalings in equations~\eqref{micro_cs} and \eqref{micro_cb},  ensuring a nontrivial limiting model when $\ve\to 0$, are obtained through nondimensionalisation of the microscopic model. In this regard, we briefly discuss biologically relevant values for the various parameters, which suggest this scaling regime is reasonable, cf., also \cite{ptashnyk2020multiscale}. 
Considering a representative cancer cell diameter of $10 \mu{\rm m}$, we obtain that the parameter $\ve$ representing the ratio between the size of a cell and the size of the tissue is of order $10^{-2}$ or smaller.  
In the microscopic model~\eqref{micro_cs}-\eqref{micro_m}, the dimension for the parameter $\kappa_s$ in addition to the time scale will also depend on the representative size of cancer cells, resulting in the  $\ve$-scaling in the boundary condition on $\Gamma^\ve$ in the non-dimensionalised problem~\eqref{micro_cs}. For the equation on the cell membrane, the diffusion coefficient $D_b$ will be of order $10^{-11}-10^{-12} {\rm cm}^2 {\rm s}^{-1}$, see e.g.~\cite{linderman1986, sako1994, simson1998, valdez2003}, which naturally yields the $\ve^2$-scaling in the non-dimensional equation in~\eqref{micro_cb} and hence dependence of $c_b$ on the microscopic cell-scale variable in the macroscopic model. Fast diffusion on the cell membrane, corresponding to a $\ve$- or $\ve^0$- scaling in the diffusion term of the non-dimensional equation for $c_b^\ve$, would result in the homogeneous distribution of membrane-bound MMPs on the cell membrane in the macroscopic tissue level model. 

We note that dynamics of species in the cell interior that influence the production of MMPs can also be included in the microscopic model and the corresponding modelling and homogenization results are discussed in~\cite{ptashnyk2020multiscale}.

For  locally-Lipschitz continuous functions $f_s$,  $f_b$, $g_b$, and $g_s$, such that
\begin{equation}\label{assump_fg}
\begin{aligned}
&\edit{ |f_b(\eta, \xi)| +  |f_s(\eta, \xi)| \leq C(1+ |\eta|+ |\xi|) \; \text{ for } \eta , \xi\in \mathbb R, \, C>0, \;  f_b(0, \xi), \, f_s(0, \xi) \geq 0 \; \text{ for }  \xi \geq 0, }\\
& 0 \leq g_b(\xi, \eta) \leq  C_b \xi \eta, \; \;  0 \leq g_s(\zeta, \eta) \leq  C_s\zeta \eta  \; \; \text{ for } \; \xi, \eta, \zeta \geq 0 \;  \text{ and  } \; C_s>0, C_b >0, \\
&  g_b(\xi, \eta) =0 , \;    \;
g_s(\xi, \eta) = 0 \;   \text{  if } \;    \xi=0 \;  \text{ or } \;  \eta =0, 
\end{aligned}
\end{equation}
and 
\begin{equation}\label{assump_coef}
\begin{aligned}
& 0 < d_s \leq \bar D_s(x) \leq \hat d_s < \infty,  \; \text{for } x \in \Omega, \quad D_b, \, \beta_s, \, \beta_b, \, \kappa_b, \, \kappa_s >0, \\
& \edit{ c_{s, \rm in}, m_{\rm in} \in L^\infty(\Omega), \, c_{b, \rm in}, m_{\Gamma, {\rm in}} \in L^\infty(\Gamma),} \;  \text{ with } \; c_{s, \rm in},   m_{\rm in} \geq 0 \; \text{ in } \Omega,  \; c_{b, \rm in}, m_{\Gamma, \rm{in}} \geq 0 \; \text{ on } \Gamma,
\end{aligned}
\end{equation}
using homogenization methods, in the same way as in~\cite{marciniak2008derivation,ptashnyk2020multiscale}, one can derive  the corresponding \edit{tissue-level (macroscopic) two-scale model}
\begin{equation}\label{macro_model}
\begin{aligned}
 &\partial_t c_s - \nabla \cdot ( D_{s}^{\rm hom}(x) \nabla c_s) = - g_s (c_s, m) 
 + \kappa_s \frac {1} { |Y_e|}  \int_\Gamma f_s(c_s, m_\Gamma) d\gamma - \beta_s c_s && \text{ in } (0,T] \times \Omega, \\
 & \partial_t c_b - \nabla_{\Gamma, y}\cdot  (D_b \nabla_{\Gamma, y} c_b) = \kappa_b f_b(c_b, m_\Gamma) - g_b(c_b, m_\Gamma) - \beta_b c_b   && \text{ in } (0,T] \times \Omega \times \Gamma, \\
 & \partial_t m \; = \edit{ - g_s(c_s, m)}  && \text{ in } (0,T] \times  \Omega , \\
  &\edit{\partial_t m_\Gamma  = -   g_b(c_b,m_\Gamma) } && \edit{\text{ in } (0,T] \times  \Omega\times \Gamma,} \\
 & D_{s}^{\rm hom} (x)\nabla c_s \cdot \nu = 0  && \text{ on } (0,T] \times \partial\Omega,\\
 & c_s(0) = c_{s, \rm in},  \;  \;  m(0) = m_{\rm in} \; \text{ in } \; \Omega, \quad \;  c_b(0) =c_{b, \rm in},  \; \; m_\Gamma(0) = m_{\Gamma, {\rm in}} \; \text { in } \;  \Omega \times \Gamma.
 \end{aligned}
\end{equation}
The effective diffusion matrix is given by
\begin{equation}\label{D_hom}
D^{\rm hom}_{s,ij}(x)= \frac 1{|Y_e|} \int_{Y_e} \big[ \bar D_{s}(x) \delta_{ij}  + \big(\bar D_{s}(x) \nabla_y w^j(y)\big)_i \big] dy, \quad \text{ for }\;\; i,j=1,\dots,d,
\end{equation}
with $w^j$ being $Y$-periodic solutions of the unit cell problems
\begin{equation}\label{eqn:cell_problems}
\begin{aligned}
&\nabla_y\cdot  \big( \bar D_s(x)(\nabla_y w^j + e_j)\big) = 0 && \text{ in } Y_e \times\O, \quad \int_{Y_e} w^j(x,y) dy = 0, \\
&  \bar D_s(x)(\nabla_y w^j + e_j) \cdot \nu = 0 && \text{ on } \Gamma \times\O, \quad w^j(x, \cdot) \; \; Y-\text{periodic},
\end{aligned}
\end{equation}
for $x\in \Omega$ and $j=1, \ldots, d$, where $\{ e_j\}_{j=1, \ldots, d}$ is the standard basis in $\mathbb R^d$. The above assumptions on the nonlinear functions $f_s, f_b, g_s, g_b$ also ensure  the well-posedness of the microscopic  \eqref{micro_cs}-\eqref{micro_m} and macroscopic \eqref{macro_model} models, see~\cite{marciniak2008derivation,ptashnyk2020multiscale}.

The model \eqref{micro_cs}-\eqref{micro_m} and the corresponding macroscopic model~\eqref{macro_model} describe cancer cell-ECM dynamics but not invasion processes, as in the latter case  we must account for movement and proliferation  of cancer cells. For proliferation, one would have to include topological changes for  which there is as yet no rigorous theory and hence we do not account for this even at the microscopic level.
Movement of the cancer cells could be incorporated by deriving the corresponding equations on time dependent domains, this has been the subject of a number of recent works, e.g., \cite{macdonald2016computational,alphonse2016coupled,alphonse2024free}. However our goal is to understand tissue level models in the limit $\ve\to0$ and to do this we assume separation of scales for processes happening at the cell and tissue level respectively. Thus we do not include movement in the microscopic models described above and instead in section~\ref{sec:macro} we propose a framework that phenomenologically accounts for movement and proliferation of cancer cells at the tissue level inspired by the above microscopic description of processes occurring at the cell level. A rigorous homogenisation of a microscopic model that includes cell movement is an interesting topic which we intend to consider in a future study.

\section{Tissue level modelling of MMP mediated cancer cell invasion}\label{sec:macro}

In this section we propose a macroscopic model for cancer cell invasion into the ECM. 
We combine aspects of mixture and multiscale modelling approaches to model the effect of bound and soluble MMPs on ECM degradation and cancer invasion.
We denote the volume fraction of the ECM by $\phi(t,x)$  and  make the modelling assumption that the corresponding volume fraction of  tumour cells is given by $\psi(t,x) = 1- \phi(t,x)$ for  $ x\in \Omega$ and  $t\in[0,T]$. In effect, we are assuming that as the cancer cells degrade the ECM they instantaneously occupy the resultant space, via proliferation and movement, both of which are not explicitly modelled. 
\edit{ 
 The modelling formalism we employ essentially consists of a mixture model or two-phase model with the tissue consisting of a mixture of cells and ECM such that their combined volume fraction is always equal to one. Such a modelling formalism is  widely used in models for tumour or tissue growth, see e.g., \cite{byrne2003modelling} and in our case serves as a phenomenological framework for modelling invasion in the present study.  Accounting for movement and proliferation of cells is an extremely interesting direction, especially in the latter case as topological changes must be considered in the microscopic model.  We intend to focus on this in future work.
}
We present two modelling frameworks for investigating MMP mediated invasion that take into account the cell-scale features discussed in section \ref{sec:micro}. The macroscopic model described in section~\ref{sec:macro_model} does not include heterogeneity of the dynamics at the cell level but takes into account changes in the tissue microstructure due to cancer invasion, whilst the two-scale model, outlined in section~\ref{sec:two_scale_model} is similar in structure to the models rigorously derived in \cite{ptashnyk2020multiscale} and allows spatial heterogeneity in the intracellular  processes inherent in the model of section~\ref{sec:micro} to be incorporated. 

\subsection{Macroscopic model}\label{sec:macro_model}
In formulation of a  macroscopic model  we assume that the MMPs can be described by concentrations at the tissue scale, i.e., we assume spatial homogeneity at the cell level.  
As in section~\ref{sec:micro}, the concentration of soluble matrix degrading enzymes is denoted by $c_s(t,x)$ and the whilst the surface density of membrane bound matrix degrading enzymes is now given by  $c_b(t,x)(1-\phi(t,x))$ for  $ x\in \Omega$ and  $t\in[0,T]$. We note the  dependence on the volume fraction of cancer cells in the latter expression as the membrane bound species are bound to the membranes of the cancer cells.
The evolution equations for the ECM and MMPs we propose are as follows
\begin{equation}\label{eqn:model_macro}
\left.
\begin{aligned}
\partial_t \phi &= - \mu_b \phi \, c_b(1-\phi) - \mu_s \phi c_s\\
\partial_tc_s - \nabla \cdot (D_s(\phi) \nabla c_s) &= \kappa_s  f_s(c_s)(1- \phi) - \mu_s \phi c_s  - \beta_s c_s \\
\partial_t (c_b(1- \phi))  &= \kappa_b f_b(\phi)(1- \phi) - \mu_b \phi c_b(1-\phi)  - \beta_b c_b  (1- \phi)
\end{aligned}
\right\}
\text{ in }\O\times(0,T],
\end{equation}
where we have accounted for  degradation of the ECM by bound and soluble MMPs, ECM volume fraction dependent diffusion of the soluble MMPs in the tissue, see section \ref{sec:effective_diffusion} for details on how this is modelled, production of soluble and bound MMPs by the  
cancer cells, degradation of the ECM by the soluble and insoluble MMPs \edit{ and the corresponding use of MMPs}, at rates $\mu_s$ and $\mu_b$ respectively, and linear decay of the MMPs.
For the degradation of the ECM by the MMPs we have made the simplest relevant assumption that this is linear with respect to the concentration of  the MMPs.  For the reaction terms governing the production of MMPs by the cancer cells we assume a Hill function form~\cite{corral2024}, i.e.,
\begin{equation}\label{MMP_production}
f_b(\phi) = \frac{\phi}{1+\phi}, \quad
f_s(c_s) = \frac 1 { 1+ c_s}.
\end{equation}
 \edit{More general models where  $f_b$ depends also on $c_b(1-\phi)$ and $f_s$  on $\phi$, as we allow in the microscopic model \eqref{micro_cs} and \eqref{micro_cb}, can  be included without major complications, given suitable growth conditions.} 
To close the model, we consider zero flux boundary condition for $c_s$ and the initial conditions   
\begin{equation}\label{IC_BC}
\begin{aligned} 
&D_s(\phi) \nabla c_s\cdot \nu = 0 \; \; && \text{ on } \partial\Omega\times(0,T], \\
&\phi(0,x) = \phi_0(x), \; \; \; c_s(0,x) = c_{s,0}(x), \; \; c_b(0,x) = c_{b,0}(x) \; \;&& \text{ for } x\in \Omega.
\end{aligned} 
\end{equation}
For biological relevance, we assume $\phi_0, c_{s,0}$ and $c_{b,0}$ to be bounded  nonnegative functions, 
\begin{equation}\label{assum_initial}
\begin{aligned} 
&   0 < \delta_0\leq \phi_0(x)\leq 1, \quad 0 \leq c_{s,0}(x)\leq M_{s,0},  \quad 0 \leq c_{b,0}(x) \leq M_{b,0}  \;  \;\text{ for } \; x \in \Omega, \\
& c_{s,0} \in  H^1(\Omega), \quad  \phi_0, \, c_{b,0}(1-\phi_0) \in W^{1,p}(\Omega)   \text{ for } p\geq 4, \\
&\edit{ D_s \in {\rm Lip}_{\rm loc}(\mathbb R)^{d\times d}, \; D_s(\xi) \eta \eta \geq d(\xi) |\eta|^2 \text{ for } \eta \in \mathbb R^d, \text{ with } d(\xi) \geq d_s \xi  \text{ for  }  d_s >0, \text{ and }  \xi \in \mathbb R_+.   }
  \end{aligned} 
\end{equation}
The assumption that the initial condition for $\phi$ is strictly positive is needed in our proof of the well-posedness of model~\eqref{eqn:model_macro}, \eqref{IC_BC}. \edit{The strict positivity of $\phi_0$ will imply the strict positivity of solutions $\phi$ of model~\eqref{eqn:model_macro}, which in turn  ensures that  $D_s(\phi)$ is uniformly positive definite and a priori estimates for $\nabla c_s$. The uniform estimates for $\nabla c_s$ are essential for the proof of the well-posedness results, in particular for the proof of strong convergence used in the fixed-point argument. Removing the restriction of strict positivity of  initial conditions for $\phi$ leads to a degenerate nonlinear parabolic equation for $c_s$. There are some  results for degenerate nonlinear parabolic equations, which rely on the specific structure of the singularity, that allow one to  prove well-posedness results for  kinetic or entropy solutions in $L^1$, see e.g.~\cite{chen2003, chen2005}.  However similar ideas can not be applied to problem~\eqref{eqn:model_macro}, \eqref{IC_BC} due to strong coupling between solutions of the degenerate  parabolic problem and solutions of the ordinary differential equations and the dependence of the degeneracy in the elliptic operator on the solutions of the ordinary differential equations. Therefore development of a new approach for systems of such type is needed, this is an interesting topic for future work.}  
The $H^1$- and $W^{1,p}$-regularity of initial conditions is required only for the proof of uniqueness of solutions. 

\subsection{Two-scale model}\label{sec:two_scale_model}
If cell-scale heterogeneity in the spatial distribution of membrane bound MMPs is relevant, this is in effect governed by the details of the microscopic model, then a two-scale model of the following form may be considered. In the absence of invasion  this model can be derived from the microscopic model in section \ref{sec:micro} cf., \cite{ptashnyk2020multiscale}.  For the two-scale model we let $c_b(t,x,y)(1-\phi(t,x))$ for  $y\in\Gamma, x\in \Omega$ and  $t\in[0,T]$ denote the density of membrane bound MMPs at a point $x$ in the tissue, $y$ on the membrane and 
time~$t$.  The corresponding system of equations  reads 
\begin{equation}\label{eqn:model_micro}
\begin{aligned}
\partial_t \phi &= -\mu_b \phi \, \ddashinttt_\Gamma c_b \diff\gamma \, (1-\phi) - \mu_s \phi c_s &&\text{ in } \; (0,T]\times \Omega,\\
\partial_tc_s - \nabla \cdot (D_s(\phi) \nabla c_s) &= \kappa_s  f_s(c_s)(1- \phi) - \mu_s \phi c_s  - \beta_s c_s &&\text{ in } \; (0,T]\times\Omega,\\
\partial_t (c_b(1- \phi))  - (1- \phi) \nabla_{\Gamma,y} \cdot (D_b \nabla_{\Gamma, y} c_b) 
&= \kappa_b f_b(\phi)(1- \phi) \\
& \quad - \mu_b \phi c_b(1-\phi)  - \beta_b  c_b  (1- \phi) &&\text{ in } \;  (0,T]\times\Omega \times \Gamma.
\end{aligned}
\end{equation}
Comparing with~\eqref{eqn:model_macro}, system \eqref{eqn:model_micro} includes the dynamics of a spatially heterogenous distribution  of $c_b$ on the cell membrane. 
For simplicity, in the  description above we have assumed that all  cells are of the same shape and size. Using techniques from   locally-periodic homogenization, see e.g.~\cite{ptashnyk2015}, we can also derive  the above system if cells varied in shape with respect to $x$ and $t$, e.g.~$\Gamma_{t,x}= \A(t,x) \Gamma$, with some given diffeomorphism $\A$,  where $\Gamma$ is a reference membrane.  Using this approach we could incorporate  the deformation and/or growth of cells, c.f., \cite{lakkis2024parallel} for an example of such a two-scale model and its numerical resolution.

\subsection{Effective diffusion tensor}\label{sec:effective_diffusion}
Both the single-scale  and two-scale macroscopic models contain an effective diffusion coefficient  that describes the diffusivity of the soluble MMPs in the tissue. This diffusivity is a nonlinear function of $\phi$, the volume fraction of the  extracellular matrix.  For the reference domain  $Y$, the corresponding value of $\phi$ is obtained by computing $\vert Y_e\vert$ and the effective diffusion matrix is computed as in \eqref{D_hom} and \eqref{eqn:cell_problems}. To obtain the function $D_s(\phi)$ we compute $D^{\rm hom}_{s}$ for a number of different values of $\phi$ i.e., on different reference domains,
\edit{ and fit an effective diffusivity on these values, see  section~\ref{sec:homg_diff_sims} and Appendix \ref{app:eff_diff} as well as Figures~\ref{fig:hom_domain_2d}, \ref{fig:hom_domain_3d} and~\ref{fig:d_eff}. We note that the fitted effective diffusivity should satisfy that $D_s(\phi)$ is positive definite for all $\phi\in(0,1]$.  In Appendix~\ref{app:full_eff_diff} we remark on one approach towards ensuring this for the general case of diffusion tensors which are not diagonal.}


\subsection{ \edit{Macroscopic model with} heterogeneous  extracellular matrix}
Following \cite{deakin2013mathematical}, we introduce a matrix suitability factor $s:\Omega\times[0,T]\to [0,1]$, which models heterogeneity in the ECM due, for example, to variations in pore size or fibre density. We assume the ability of cancer cells to invade the ECM is a decreasing function of $s$, i.e., $s=1$ and $s=0$ correspond  \edit{to ECM configurations that are the least and most suitable for cancer cell invasion, respectively}.    We further assume that only the membrane bound MMPs can change the suitability of the ECM for invasion \cite{sabeh2009protease}. Assuming a simple linear dependence of the  invasion related parameters on $s$ and that the only effect of matrix suitability is to inhibit the degradation of matrix by MMPs yields the model:
\begin{equation}\label{het_ecm}
\left.
\begin{aligned}
\partial_t \phi &= - (1-s)\mu_b \phi \, c_b(1-\phi) - (1-s)\mu_s \phi c_s\\
\partial_tc_s -\nabla\cdot (D_s(\phi) \nabla c_s) &= \kappa_s  f_s(c_s)(1- \phi) 
\\&\quad - (1-s)\mu_s \phi c_s - \beta_s c_s \\
\partial_t (c_b(1- \phi))  &= \kappa_b f_b(\phi)(1- \phi) 
\\&\quad- (1-s)\mu_b \phi \, c_b(1-\phi)  - \beta_b c_b  (1- \phi)\\
\partial_t s  & = - \delta_s c_b(1-\phi)s
\end{aligned}
\right\}
\text{ in }\O\times(0,T],
\end{equation}
with boundary and initial conditions as in \eqref{IC_BC} and additional  initial condition  $s(0, x) = s_0(x)$, where  $0 \leq s_0(x) \leq 1$ for $x \in \Omega$. One could also consider that the diffusion coefficient $D_s$ depends on $s$ but we do not include this in the current work to allow for comparison with the results in~\cite{deakin2013mathematical}.

\section{Analysis of the macroscopic models }\label{sec:analysis}
\edit{
Following the same arguments as  in~\cite{ptashnyk2020multiscale} we can prove the following well-posedness  result for the macroscopic two-scale model~\eqref{macro_model}. 
\begin{theorem}\label{existence_two-scale}
   Under assumptions~\eqref{assump_fg}, \eqref{assump_coef}  there exists a unique weak solution $c_s \in L^2(0,T; H^1(\Omega))$, $c_b \in L^2((0,T)\times\Omega; H^1(\Gamma))$, $m \in H^1(0,T; L^2(\Omega ))$, $m_\Gamma \in H^1(0,T; L^2(\Omega \times \Gamma))$ of the two-scale problem~\eqref{macro_model}, with $c_s(t,x), c_b(t,x,y), m(t,x), m_\Gamma (t,x,y) \geq 0$ for 
   $(t,x) \in (0,T)\times \Omega$ and  $y \in \Gamma$,  $c_s, m \in L^\infty((0,T)\times \Omega)$, $c_b, m_\Gamma \in L^\infty((0,T)\times \Omega\times \Gamma)$,  and $\partial_t c_s \in L^2(0,T; H^1(\Omega)^\prime)$, $\partial_t c_b \in L^2((0,T)\times \Omega; H^1(\Gamma)^\prime)$. 
\end{theorem}
}

The challenges in proving the existence and uniqueness of solutions to the macroscopic model~\eqref{eqn:model_macro},~\eqref{IC_BC} arise due to the fact that the system consists of coupled ODEs and PDEs with a  possibly degenerate diffusion. For the model with heterogeneous suitability \eqref{het_ecm} the proof of the well-posedness is similar to that for~\eqref{eqn:model_macro},~\eqref{IC_BC} and hence omitted.

A weak solution of system \eqref{eqn:model_macro} with initial and boundary conditions~\eqref{IC_BC} is defined as follows. 

\begin{definition}
 Functions $\phi, c_b(1-\phi) \in H^1(0,T; L^2(\Omega)) \cap L^\infty((0,T)\times \Omega)$ and $c_s \in L^\infty(0,T; L^2(\Omega)) \cap L^2(0,T; H^1(\Omega))$  are weak solutions of \eqref{eqn:model_macro},~\eqref{IC_BC} if
 \begin{equation}
  \int_0^T\!\!\Big[ -\langle  c_s, \partial_t w \rangle + \langle D_s(\phi) \nabla c_s, \nabla w \rangle  \Big] dt =\int_0^T \!\!\langle \kappa_s f_s(c_s) (1-\phi) - \mu_s \phi c_s - \beta_s c_s, w \rangle dt  + \langle c_{s, 0} , w(0) \rangle
 \end{equation}
 for every $w \in L^2(0,T; H^1(\Omega))$, with $\partial_t w \in L^2((0,T)\times \Omega)$ and $w(T) = 0$, and the \edit{equations}  and initial conditions for $\phi$ and $c_b(1-\phi)$ are satisfied a.e.\ in $(0,T)\times \Omega$.
\end{definition}
Here we introduced the notation $\langle u, v \rangle = \int_\Omega u v dx$, for $u\in L^p(\Omega)$ and $v \in L^q(\Omega)$ with $1/p+1/q =1$.

\begin{theorem}\label{th:wellpos}
 Under the assumptions~\eqref{assum_initial}  there exists a unique weak solution of problem~\eqref{eqn:model_macro},\eqref{IC_BC}, satisfying
 $$
 0\leq \phi(t,x)  \leq 1, \quad 0 \leq c_b(t,x)(1- \phi(t,x)) \leq M_b, \quad 0 \leq c_s(t,x) \leq M_s ,
 $$
 for $x \in \Omega$ and $t\geq 0$, where $M_b= \max\{M_{b,0}, {\kappa_b}/{\beta_b}\}$ and $M_s= \max\{M_{s,0}, {\kappa_s}/{\beta_s}\}$.

 Additionally over any finite time interval $[0,T]$ we have
 \begin{equation}\label{lower_bound_phi}
 \phi(t,x) \geq \tilde \delta =  \delta_0 \exp(- (\mu_b \kappa_b/\beta_b + \mu_s \kappa_s/\beta_s) T) >0.
 \end{equation}
\end{theorem}

\edit{ \begin{remark} The proof of Theorem~\ref{th:wellpos} is structured as follows. 
\begin{itemize}
\item First we use a fixed point argument considering a given $0< \tilde \delta\leq \bar \phi \leq 1$ in the equation for $c_s$, see~\eqref{fixed_point_model}, which decouples it from the  equations for $\phi$ and $c_b$ and ensures uniform ellipticity.  
\item Applying the method of positive invariant regions for ODEs,  e.g.~\cite{amann1990}, and for reaction-diffusion equations,e.g.~\cite{smoller1994}, we show nonegativity and boundedness of solutions of~\eqref{fixed_point_model},~\eqref{IC_BC}, assuming their existence. The same arguments ensure  the nonegativity and boundedness for the Galekin approximations, and hence also for the limits of the approximation sequences.  Note that convex sets are weakly sequentially closed. 
\item Then using the Galerkin method we prove the  existence of a weak solution of the  problem with a given $\bar \phi$ in the equation for $c_s$, see~\eqref{fixed_point_model}.
\item To pass to the limit in the nonlinear terms we prove the strong convergence of the Galerkin approximation by showing the Cauchy properties of the corresponding sequences. 
\item To complete the existence proof we show that the  map $\bar \phi \mapsto \phi$, defined by solutions of problem~\eqref{fixed_point_model},~\eqref{IC_BC}, has a fixed point by proving the strong convergence of the corresponding sequences and applying the Schauder fixed-point theorem.  To show strong convergence in $L^2((0,T)\times \Omega)$ we combine  the Cauchy properties for solutions of the ODEs with the regularity properties and Aubin-Lions compactness lemma for solutions of the reaction-diffusion equation.  
\item Local Lipschitz continuity of the nonlinear functions and boundedness of solutions are applied to prove uniqueness.
\item To estimate the terms arising due to the diffusion matrix $D_s$ being a function of $\phi$ we use higher regularity properties of solutions of the parabolic equations and the structure of the ODEs,  together with the fact that the regularity of $\phi$ and $c_b(1-\phi)$ with respect to the space variables is determined by the regularity of $c_s$.  
\end{itemize} 
\end{remark}
}

\noindent{\bf Proof. }
To prove the existence result we combine the fixed point  and   Galerkin methods.
For given  
$$\bar \phi \in V, \; \; \text{ where } \;  V =\big\{ v \in L^2((0,T)\times\Omega) \,   : \; 0 < \tilde \delta \leq v(t,x) \leq 1  \; \text{ for } \; (t,x) \in (0,T)\times \Omega\big \}, $$
 consider the map $K: V \to V$ with $\phi= K (\bar \phi) $ being the solution of system 
\begin{equation} \label{fixed_point_model}
 \begin{aligned}
 \partial_t \phi & = - \mu_b \phi c_b (1- \phi) - \mu_s \phi c_s, \\
 \partial_t (c_b(1- \phi))  &= \kappa_b f_b(\phi)(1- \phi) - \mu_b  \phi c_b(1-\phi)  - \beta_b c_b  (1- \phi),\\
   \partial_tc_s - \nabla \cdot (D_s(\bar \phi) \nabla c_s) &= \kappa_s  f_s(c_s)(1- \bar\phi) - \mu_s \bar\phi c_s  - \beta_s c_s, 
 \end{aligned}
\end{equation}
with initial and boundary conditions as in \eqref{IC_BC}.

For $\tilde \delta \leq \bar \phi \leq 1$,  applying  the method of positive invariant regions, see e.g.~\cite{smoller1994} and using the  nonnegativity of  $c_{s,0}$  we  obtain  $c_s (t,x)\geq 0$ for $x \in \Omega$ and $t\geq 0$.  
Using the theory of invariant regions for ordinary differential equations (ODEs), see e.g.~\cite{amann1990} and  $\phi_0(x)\geq 0$ for $x \in \Omega$,  from  the equation for $\phi$ it directly follows that $\phi(t,x) \geq 0 $ for $x \in \Omega$ and $t\geq 0$. Considering
first the positive part  $(1-\phi)_{+}$ instead of $(1-\phi)$ in the first term on the right hand side of the  equation for $c_b$ and using the assumptions on $c_{b,0}$ and $\phi_0$, we obtain $c_b(t,x)(1- \phi(t,x)) \geq 0$ for $x \in \Omega$ and $t\geq 0$. 
The non-negativity of $c_s$ and $c_b(1-\phi)$    and the fact that $\phi_0(x) \leq 1$ imply that $\phi(t,x) \leq 1$  for  $x \in \Omega$ and $t\geq 0$. Then the assumption $\phi_0(x) \geq \delta_0 $ ensures the lower bound  \eqref{lower_bound_phi} for any finite time interval $[0,T]$.  Hence, the equations for $c_b$ with $(1-\phi)_+$ and
with $(1-\phi)$ are equivalent.  Considering the properties of the reaction terms and applying  the method of positive invariant regions  we also obtain $c_s\leq M_s $ and $c_b(1- \phi) \leq M_b(1- \phi)$.

To show existence of solutions of \eqref{fixed_point_model} for  $\bar \phi \in V$, we apply the Galerkin method and consider  the sequence of  approximations $\{\phi^k\}$, $\{c_s^k\}$ and $\{c_b^k\}$ for solutions of system \eqref{fixed_point_model} given by
$$
\phi^k(t,x) = \sum_{j=1}^k a_j^k(t) \psi_j(x), \quad c^k_b (t,x)= \sum_{j=1}^k b_j^k(t) \psi_j(x), \quad c^k_s (t,x)= \sum_{j=1}^k d_j^k(t) w_j(x),
$$
and satisfying
\begin{equation}\label{eq:Galerkin}
\begin{aligned}
 & \langle \partial_t \phi^k, \psi_j\rangle  = -   \langle  \mu_b \phi^k \, c_b^k(1-\phi^k) + \mu_s \phi^k  c_s^k, \psi_j \rangle, \\
& \langle \partial_t (c_b^k(1- \phi^k)), \psi_j \rangle  = \langle \kappa_b f_b(\phi^k)(1-  \phi^k) - \mu_b \phi^k c_b^k(1- \phi^k)  - \beta_b c_b^k  (1-  \phi^k), \psi_j \rangle,\\
& \langle \partial_t c_s^k, w_j \rangle + \langle D_s( \bar \phi) \nabla c^k_s, \nabla w_j \rangle  = \langle \kappa_s  f_s(c_s^k)(1- \bar\phi) - \mu_s \bar\phi c^k_s  - \beta_s c_s^k, w_j \rangle, 
\end{aligned}
\end{equation}
for  $t \in (0,T]$ and initial conditions 
\begin{equation}\label{IC_k}
\phi^k(0,x) = \phi^k_0(x), \quad c^k_b(0,x) = c^k_{b,0}(x), \quad 
c_s^k(0,x) = c^k_{s,0}(x),
\end{equation}
where $\{\psi_j\}_{j \in \mathbb N}$ are the basis functions of $L^2(\Omega)$ and $\{w_j\}_{j\in \mathbb N}$ are  the basis functions of $H^1(\Omega)$,  and $\phi^k_0$, $c_{b,0}^k$ and $c_{s,0}^k$ are projections of $\phi$, $c_{b.0}$ and $c_{s,0}$ on $\text{span}\{\psi_1, \ldots, \psi_k\}$ and $\text{span}\{w_1, \ldots, w_k\}$  respectively. 
 Notice that the same arguments as above ensure nonnegativity and boundedness for $\phi^k, c_s^k$ and $c_s^k(1-\phi^k)$.  
 For given $\bar \phi$ the equation for $c_s^k$ in \eqref{eq:Galerkin} is decoupled from the equations for $\phi^k$ and $c_b^k$ and comprises a system of ODEs for $d_j^k$, with $j=1, \ldots, k$. Since $f_s$ is Lipschitz-continuous, ODE-theory, see e.g.~\cite{amann1990}, ensures existence of a unique solution~$(d_1^k, \ldots, d_k^k)$ which is absolutely continuous. The local Lipschitz continuity of the nonlinear functions in the equations for $\phi^k$ and $c_b^k$ and the boundedness of $\phi^k$ and $c_b^k$ ensure the existence of a unique absolutely continuous function satisfying  the systems of ODEs for $(a_1^k,\ldots, a_k^k, b_1^k, \ldots, b_k^k)$. 

Taking $\phi^k|\phi^k|^{p-2}$,  $c_b^k(1-\phi^k)|c_b^k(1- \phi^k)|^{p-2}$,  for $p\geq 2$, and  $c_s^k$ as test functions  in equations for $\phi^k$, $c_b^k$, and  $c_s^k$ respectively, and using the nonnegativity and boundedness of $\phi^k$, $c_b^k(1-\phi^k)$ and $c_s^k$, we obtain the following a priori  estimates
\begin{equation} \label{aprioriestim}
 \begin{aligned}
 &\sup_{(0,T)} \|\phi^k(t) \|^p_{L^p(\Omega)} \leq  \|\phi(0) \|^p_{L^p(\Omega)},   \quad \text{ for } 1<p<\infty, \\
 & \sup_{(0,T)} \|c_b^k(t)(1-\phi^k(t)) \|^p_{L^p(\Omega)} \leq C\big(1+ \|\phi(0) \|^p_{L^p(\Omega)} \big),\\
 &\sup_{ (0,T)} \|c_s^k(t) \|^2_{L^2(\Omega)} + \| (d(\bar\phi))^{\frac 12} \nabla c_s^k \|^2_{L^2((0,T)\times\Omega)} \leq C,\\
 \end{aligned}
\end{equation}
where the  constant $C>0$  depends on $T$, $M_b$, and $M_s$ \edit{and is independent of $k$}. 
Using \eqref{aprioriestim}, from the equations for $\phi^k$, $c_b^k$, and $c_b^k$ we also obtain
\begin{equation}\label{estim:time_diff}
\|\partial_t \phi_k \|_{L^2((0,T)\times \Omega)} + \|\partial_t (c_b^k(1-\phi_k)) \|_{L^2((0,T)\times \Omega)} + \| \partial_t c_s^k\|_{L^2(0,T; (H^1(\Omega)^\prime)}\leq C,
\end{equation}
\edit{see e.g.~\cite{evans2010} for the derivation of the estimate for $\partial_t c_s^k$}. To pass to the limit in the nonlinear terms we require strong or almost everywhere convergence of $\{\phi^k\}$, $\{c_s^k\}$ and $\{c_b^k\}$. To this end we  show that $\{\phi^k\}$, $\{c_s^k\}$ and $\{c_b^k(1- \phi^k)\}$ are Cauchy sequences. The uniform boundedness of $\phi^k$, $c_s^k$ and $c_b^k$ and  the equations for $\phi^k - \phi^l$, $c_s^k - c^l_s$ and $c_b^k(1- \phi^k) - c_b^l(1- \phi^l)$  yield 
\begin{equation} \label{Cauchy}
 \begin{aligned}
\|c_s^k (t)- c_s^l(t) \|_{L^2(\Omega)}^2  + \|(d(\bar \phi))^{\frac 1 2} \nabla (c_s^k - c_s^l) \|_{L^2((0,t)\times\Omega)}^2    +  \| \phi^k (t)- \phi^l(t)\|_{L^2(\Omega)}^2 \qquad \qquad \\
+ \|c_b^k(t)(1- \phi^k(t)) - c_b^l(t)(1-\phi^l(t))  \|_{L^2(\Omega)}^2\leq  \| \phi^k_0- \phi^l_0\|_{L^2(\Omega)}^2 + \|c_{s,0}^k- c_{s,0}^l\|_{L^2(\Omega)}^2\\
+ \|c_{b,0}^k(1-\phi_0^k) - c_{b,0}^l(1-\phi^l_0)  \|_{L^2(\Omega)}^2 
+  C \!\int_0^t \!\Big[ \| \phi^k(\tau)- \phi^l(\tau)\|_{L^2(\Omega)}^2 \qquad 
\\ + \|c_s^k (\tau)- c_s^l(\tau) \|_{L^2(\Omega)}^2     + \|c_b^k(\tau)(1- \phi^k(\tau))- c_b^l(\tau)(1-\phi^l(\tau))\|_{L^2(\Omega)}^2  \Big] d\tau.
 \end{aligned}
\end{equation}
Using the strong convergence in $L^2(\Omega)$ of the Galerkin approximations of the initial conditions, e.g.~$\phi^k_0 \to \phi_0$, $c_{s,0}^k \to c_{s,0}$, and
$c_{b,0}^k \to c_{b,0}$ strongly in $L^2(\Omega)$ as $k \to \infty$, and applying the Gronwall inequality we obtain
\begin{equation}\label{cauchy_property}
\begin{aligned}
 &\sup_{(0,T)}\| \phi^k- \phi^l\|_{L^2(\Omega)}^2  + \sup_{(0,T)} \|c_s^k - c_s^l \|_{L^2(\Omega)}^2 + \sup_{(0,T)} \|c_b^k (1 - \phi^k) - c_b^l(1-\phi^l)  \|_{L^2(\Omega)}^2 \leq \sigma(k,l),
 \end{aligned}
\end{equation}
 where $\sigma(k,l) \to 0$ as $k, l \to \infty$.  The  Cauchy property~\eqref{cauchy_property} implies the strong convergence of  $\{\phi^k\}$, $\{c_s^k\}$ and $\{c_b^k (1- \phi^k)\}$  in $L^2((0,T)\times \Omega)$. The strong convergence of $\phi_k$ and $c_b^k (1- \phi^k)$ ensures the convergence
 $c_b^k (1- \phi^k) \to c_b (1- \phi)$ a.e.\ in $(0,T)\times \Omega$.  Then the a priori estimates~\eqref{aprioriestim} and \eqref{estim:time_diff}, uniform in $k$, together with convergence a.e.~of $c_b^k (1- \phi^k)$, imply
 \begin{equation} \label{convergence}
  \begin{aligned}
 &   c^k_s \rightharpoonup c_s \qquad \qquad \text{weakly$-\ast$ in } L^\infty(0,T; L^2(\Omega)),  \; && \text{ weakly in }  L^2(0,T; H^1(\Omega)), \\
 &   c_s^k \to c_s, \quad \phi^k \to \phi, \quad c_b^k(1-\phi^k) \to c_b(1- \phi) \qquad && \text{ strongly in  } L^2((0,T)\times \Omega), \\
  &   \partial_t c_s^k \rightharpoonup \partial_t  c_s  && \text{ weakly in } L^2(0,T; H^1(\Omega)^\prime), \\
 &   \partial_t \phi^k \rightharpoonup \partial_t \phi, \quad \partial_t (c_b^k(1-\phi^k)) \rightharpoonup \partial_t (c_b (1 -\phi)) \qquad  && \text{ weakly in } L^2((0,T)\times \Omega),
  \end{aligned}
 \end{equation}
as $k \to \infty$. Using the convergence results in \eqref{convergence} we can pass to the limit as $k\to \infty$ in equations~\eqref{eq:Galerkin}  and obtain that $\phi$, $c_s$, and $c_b$ are solutions of \eqref{fixed_point_model}.

To prove the existence of a fixed point of the map $K$, we apply the Schauder fixed point theorem and show that $K: V \to V$ is compact, see e.g.~\cite{evans2010}. For this,  by showing equicontinuity and, applying Riesz-Fr\'eshet-Kolmogorov convergence result, see e.g.~\cite{brezis2011}, we prove strong convergence of a sequence  of solutions $(\phi^n, c_b^n, c_s^n)$ of problem~\eqref{fixed_point_model},~\eqref{IC_BC} with  $\bar \phi^n \in V$.

\edit{Using the a priori estimates for solutions of~\eqref{fixed_point_model},~\eqref{IC_BC}, obtained from~\eqref{aprioriestim} and \eqref{estim:time_diff}  by taking the limit as $k\to \infty$   and employing the lower semi-continuity of a norm,  together with the fact that $\tilde \delta \leq \bar \phi^n \leq 1$ uniformly in $n$,}  and applying the Aubin-Lions compactness lemma~\cite{lions1969},  we obtain the strong convergence of $\{c_s^n\} $ in $L^2((0,T)\times \Omega)$.

Considering equations for $\phi^n(t, x+ h)-\phi^n(t,x)$ and  $c_b^n(t, x+ h) - c_b^n(t,x)$ and taking $\phi^n(t,x+h) - \phi^n(t,x)$ and $c_b^n(t,x+h)(1- \phi^n(t, x+ h)) - c_b^n(t,x) (1- \phi^n(t,x))$ as test functions respectively,  we obtain
$$
\begin{aligned}
 \big\|\phi^n(t,\cdot+h) - \phi^n(t)\big\|_{L^2(\Omega_h)}^2  + \big\|c_b^n(t,\cdot+h)(1-\phi^n(t, \cdot+h)) - c_b^n(t)(1- \phi^n(t))\big\|^2_{L^2(\Omega_h)}\qquad \\
 \leq \big\| \phi_0(\cdot+h) - \phi_0\big\|^2_{L^2(\Omega_h)} + \big\|c_{b,0}(\cdot+h)(1-\phi_0(\cdot+h)) - c_{b,0}(1- \phi_0) \big\|^2_{L^2(\Omega_h)} \qquad \\
  + C \int_0^t\Big[\big\|\phi^n(\tau)(c_s^n(\tau,\cdot+h) - c_s^n(\tau))\big\|_{L^2(\Omega_h)}^2
 +  \big\|\phi^n(\tau,\cdot+h) - \phi^n(\tau)\big\|_{L^2(\Omega_h)}^2 \quad  \\
 + \big \|c_b^n(\tau,\cdot+h)(1-\phi^n(\tau, \cdot+h)) - c_b^n(\tau)(1- \phi^n(\tau))\big\|^2_{L^2(\Omega_h)} \Big] d\tau,
\end{aligned}
$$
for some $h \in \mathbb R^d$, where $\Omega_h = \{ x \in \Omega: {\rm dist}(x, \partial \Omega) >|h| \}$. 
Then, using the estimate
$$
\int_0^t\|\phi^n(\tau)(c_s^n(\tau,\cdot+h) - c_s^n(\tau))\|_{L^2(\Omega_h)}^2 d\tau\leq C_1 h^2 \|\phi^n \nabla c_s^n\|^2_{L^2((0,T)\times \Omega)} \leq C h^2,
$$
where the constant $C>0$ depends on $T$, $\tilde \delta$, $M_b$, and $M_s$, \edit{and is independent of $n$, due to uniform in~$n$ boundedness of $\phi^n$ and $\bar \phi^n\geq \tilde \delta$},  and  applying the Gronwall lemma  yields
\begin{equation}\label{equicontin}
\begin{aligned}
 \sup\limits_{t \in (0,T)}\|c_b^n(t,\cdot+h)(1-\phi^n(t, \cdot+h)) - c_b^n(t)(1- \phi^n(t))\|^2_{L^2(\Omega_h)} \\
 +\sup\limits_{t \in (0,T)}  \|\phi^n(t,\cdot+h) - \phi^n(t)\|_{L^2(\Omega_h)}^2   & \leq \sigma(h),
 \end{aligned} 
 \end{equation}
where $\sigma(h) \to 0 $ as $h\to 0$, uniformly in $n$. This,  together with estimate~\eqref{estim:time_diff} for the time derivatives and boundedness of $\phi^n$ and $c_b^n(1-\phi^n)$,  ensures the equicontinuity  and hence the strong convergence of the sequences $\{\phi^n\}$ and $\{ c_b^n(1-\phi^n)\}$ in $L^2((0,T)\times \Omega)$.
The property  $\tilde \delta \leq \phi(t,x) \leq 1$ for $(t,x) \in (0,T)\times \Omega$ follows directly from the equation for $\phi$ in \eqref{fixed_point_model} and the assumptions on the initial conditions. 

Thus we have that $K$ maps $V$ into itself and compactness is ensured by the strong convergence results. Applying the Schauder fixed point theorem we  obtain existence of solution of model~\eqref{eqn:model_macro},~\eqref{IC_BC}.

To show uniqueness, we consider the equations satisfied by  the difference of two solutions and using energy arguments,  similar as in the derivation of estimates~\eqref{aprioriestim}, obtain
\begin{equation}\label{estim_uniq_1}
\begin{aligned}
  \|\phi_1(t)- \phi_2(t)\|^p_{L^p(\Omega)} + 
  \|c_{b,1}(t)(1-\phi_1(t))- c_{b,2}(1-\phi_2(t))\|^p_{L^p(\Omega)}\quad  \\
  \leq C_p \int_0^t \| c_{s,1} (\tau)- c_{s,2}(\tau)\|^p_{L^p(\Omega)} d\tau,
\end{aligned}
\end{equation}
for $p\geq 2$. Considering $\partial_t c_s $ as test function in equations for $c_s$ in \eqref{het_ecm}, together with a regularisation by a mollifier, see e.g.~\cite{showalter1996}, \edit{and assumptions on $c_{s, 0}$ and $D_s(\phi)$}, yields 
$$
\begin{aligned}
    \| \partial_t c_s \|^2_{L^2((0,T)\times \Omega)} + \sup\limits_{(0,T)} \|\sqrt{d(\phi) }\nabla c_s \|^2_{L^2(\Omega)} \leq  C_1 +C_2 \| \partial_t  \phi\|_{L^\infty((0,T)\times \Omega)} \|\nabla c_s\|^2_{L^2((0,T)\times \Omega)} \leq C.
\end{aligned}
$$
This, together with the $L^p$-regularity for elliptic problems in $C^{1,1}$  or convex domains, yields 
$$
\|c_s\|_{L^2(0,T;W^{1,p}(\Omega))}  \leq C,  \quad 3<p\leq 7,
$$
 see e.g.~\cite{dauge1992}.
For two solutions $c_{s,1}$ and $c_{s,2}$ we have 
\begin{equation}\label{estim_uniq_2}
\begin{aligned} 
\partial_t \|c_{s,1}(t) - c_{s,2}(t)\|_{L^2(\Omega)}^2 
+ \| \sqrt{d(\phi_1(t))} \nabla (c_{s,1} (t)- c_{s,2}(t))\|_{L^2(\Omega)}^2 
\leq  C \Big[ \|c_{s,1}(t) - c_{s,2}(t)\|_{L^2(\Omega)}^2
\\ +\|\nabla c_{s,2}(t)\|^2_{L^4(\Omega)} \|\phi_1(t) - \phi_2(t)\|_{L^4(\Omega)}^2 
+ \|\phi_1(t) - \phi_2(t)\|_{L^2(\Omega)}^2 \Big]. 
\end{aligned}
\end{equation}
Using \eqref{estim_uniq_1} and the estimate 
$$
\|v\|_{L^4(\Omega)} \leq C \left[\|\nabla v \|^\theta_{L^2(\Omega)} \| v \|^{1-\theta}_{L^2(\Omega)} + \| v\|_{L^2(\Omega)}\right], 
$$
with $\theta = d/4$ and $d=2,3$, which follows from the Gagliardo-Nirenberg interpolation inequality, 
and integrating \eqref{estim_uniq_2} with respect to the time variable, we obtain 
$$
\begin{aligned} 
 \|c_{s,1} (t)- c_{s,2}(t)\|_{L^2(\Omega)}^2 
+ 2\| \sqrt{d(\phi_1)} \nabla (c_{s,1} - c_{s,2})\|_{L^2((0,t)\times \Omega)}^2 \leq  \varsigma\|\nabla c_{s,1} - \nabla c_{s,2}\|_{L^2((0,t)\times \Omega)}^2
\\
+ C \|c_{s,1} - c_{s,2}\|_{L^2((0,t)\times \Omega)}^2,  
\end{aligned}
$$
for $t \in (0,T]$. 
Choosing a sufficiently small $0<\varsigma\leq 2 \tilde \delta $ and applying the Gronwall inequality yields
$c_{s,1} = c_{s,2}$ a.e.~in $(0,T)\times \Omega$.  Then estimate \eqref{estim_uniq_1} implies $\phi_1 = \phi_2$ and $c_{b,1}=c_{b,2}$ a.e.~in $(0,T)\times \Omega$, and hence uniqueness of solutions of model~\eqref{eqn:model_macro}, \eqref{IC_BC}. 

The regularity $\phi \in L^\infty(0,T; C(\overline\Omega))$ required for the results in~\cite{dauge1992} is ensured by combining a regularisation argument, the fact  that $W^{1,p}(\Omega) \subset C^{0, \alpha}(\overline \Omega)$ for $p>d$ and $\alpha = 1 - d/p$, and the estimate
$$
\|\nabla \phi(t)\|_{L^p(\Omega)} \leq C \Big(\|\nabla \phi_0 \|_{L^p(\Omega)} +
\|\nabla (c_{b,0} (1-\phi_0) )\|_{L^p(\Omega)}
+ \Big\|\int_0^t |\nabla c_s(\tau)| d\tau \Big\|_{L^p(\Omega)} \Big) \;\text{ for } t \in (0,T]. 
$$
 The  last estimate follows from applying the Gronwall inequality and taking the $L^p(\Omega)$-norm in the      inequality
$$
|\nabla \phi(t,x)| \leq |\nabla \phi_0(x)| + C_1 |\nabla (c_{b,0}(x)(1- \phi_0(x)))| + C_2\int_0^t\big(|\nabla c_s(\tau,x)| +
|\nabla \phi(\tau,x) | \big) d \tau, 
$$ 
for $(t,x) \in (0,T]\times \Omega$, which in turn is obtained from the equations for $\phi$ and $c_b(1- \phi)$ in~\eqref{eqn:model_macro}.  

\edit{\begin{remark}Using similar arguments as in the proof of Theorem~\ref{th:wellpos} we can also show the existence and uniqueness of solutions $\phi \in H^1(0,T; L^2(\Omega))$,  $c_s \in L^\infty(0,T; L^2(\Omega)) \cap L^2(0,T; H^1(\Omega))$, $ c_b(1-\phi) \in L^\infty(0,T; L^2(\Omega\times \Gamma)) \cap L^2((0,T)\times\Omega; H^1(\Gamma))$, with $\partial_t c_s \in L^2(0,T; H^1(\Omega)^\prime)$, $\partial_t(c_b(1-\phi)) \in L^2((0,T)\times \Omega; H^1(\Gamma)^\prime)$,   to the two-scale model~\eqref{eqn:model_micro}, satisfying
 $$
 0\leq \phi(t,x)  \leq 1, \quad 0 \leq c_b(t,x,y)(1- \phi(t,x)) \leq M_b, \quad 0 \leq c_s(t,x) \leq M_s ,
 $$
 for $x \in \Omega$, $y \in \Gamma$,  and $t\geq 0$, where $M_b= \max\{M_{b,0}, {\kappa_b}/{\beta_b}\}$ and $M_s= \max\{M_{s,0}, {\kappa_s}/{\beta_s}\}$, and the lower bound~\eqref{lower_bound_phi} for $t \in [0,T]$. 
 Notice that  $\phi$ is independent of $y \in \Gamma$ and we can use $(1-\phi)c_b$ as a test function in the equation for $c_b$ in~\eqref{eqn:model_micro}. Additionally we have that $\ddashinttt_\Gamma c_b d\gamma (1- \phi)$ satisfies the same ODE equation as $ c_b (1- \phi)$ in~\eqref{eqn:model_macro}. More specifically, integrating equation for $c_b $ in~\eqref{eqn:model_micro} over $\Gamma$ yields  
 $$
 \begin{aligned}
 \partial_t \Big(\ddashinttt_\Gamma c_b d\gamma (1- \phi) \Big) = 
  \kappa_b f_b(\phi)(1- \phi)  - \mu_b \phi \,\ddashinttt_\Gamma c_b d \gamma (1-\phi)  - \beta_b \, \ddashinttt_\Gamma c_b  d\gamma (1- \phi), 
\end{aligned} 
 $$
 where we used that the integral over $\Gamma$ of  $\nabla_{\Gamma, y} \cdot (D_b(1-\phi) \nabla_{\Gamma, y} c_b)$ is zero. 
\end{remark}
}

\section{Finite element approximation}\label{sec:FEM}
In this section we describe a numerical method for the approximation of the macroscopic model formulated in section~\ref{sec:macro_model}. We do not present results for the approximation of the two-scale model as under the simplified intracellular processes considered in this study cell-scale concentrations are effectively spatially uniform at the timescales of invasion. \edit{A number of recent cell-scale models for cellular signalling processes demonstrate inhomogeneous membrane concentrations typically by modelling detailed biochemistry on the membrane or within the cell e.g., \cite{levine2005membrane,neilson2011modeling,jilkine2011comparison,garcke2016coupled}.
}
Our framework readily allows for such complex intracellular dynamics to be incorporated and simulated which we intend to consider in future work. \edit{In particular, developing a model for the intracellular signalling processes which lead to the production of MMPs, both soluble and membrane bound, by cancer cells is an interesting direction for future studies.} The numerical approximation of the two-scale model is described in detail in~\cite{ptashnyk2020multiscale,lakkis2024parallel}.

We define computational domains $\O_h$ and $Y_{h,e}$ by requiring that $\O_h$ and $Y_{h,e}$ are polyhedral approximations to $\O$ and $Y_e$ respectively. We assume that $\O_h$ and $Y_{h,e}$  consist of the union of $d$ dimensional simplices (triangles for $d=2$ and tetrahedra for $d=3$).  

 We define $\Th$ and $\She$ to be triangulations of $\O_h$ and $Y_{h,e}$ respectively and assume that each consists of  closed non-degenerate simplices and that the triangulations are conforming.  
 We define piecewise linear finite element spaces
 as follows
 \begin{align*}
 \Vho&=\left\{\Phi\in C(\O_h) : \; \Phi\vert_k\in\mathbb{P}^1,\myall k\in\Th\right\},\\
 \Vhe&=\left\{\Phi\in \Hil{1}_{\rm{ per}}(Y_{h,e})\cap C(Y_{h,e})\edit{\text{ with }\int_{Y_{h,e}}\Phi \, dx=0} : \;\Phi\vert_k\in\mathbb{P}^1,\myall k\in\She\right\},
\end{align*}
where $\Hil{1}_{\rm{ per}}(Y_{h,e})$ denotes the subspace of $Y$-periodic functions in $\Hil{1}(Y_{h,e})$.
The numerical scheme for the solution of the cell problems to compute the diffusion tensor $D_s^{\rm hom}$ is, for  $j=1,\dots,d$,  determine  $W^j\in\Vhe$ such that  
\begin{align}\label{eqn:scheme_diff_tensor} 
&\Ltwop{\bar{D}_s(\nabla_yW^j+ e_j)}{\nabla_y \Phi}{Y_{h,e}} =0, 
\end{align} 
for all $\Phi\in\Vhe$.

For the approximation of solutions of the macroscopic  model~\eqref{het_ecm}, \eqref{IC_BC} we employ an IMEX scheme  \cite{ruuth1996implicit}. We consider a uniform (purely for simplicity) partition of the time interval $[0,T]$ into $N$ subintervals with step size $\tau=T/N$. 
The scheme is  defined as follows: \\ For $n=1,\dots,N$ and given $(\Phi^{n-1}, C_s^{n-1}, C_b^{n-1}, S^{n-1})\in (\Vho)^4 $ find $(\Phi^{n}, C_s^{n}, C_b^{n}, S^{n})\in (\Vho)^4 $ such that 
\begin{equation}\label{FE_scheme}
\begin{aligned}
&\Ltwop{\Phi^n}{\Psi_1}{\Omega_{h}} =\Ltwop{\Phi^{n-1}}{\Psi_1}{\Omega_{h}}+\tau\Ltwop{-(1-S^{n-1})\mu_s\Phi^{n-1}C_s^{n-1}}{\Psi_1}{\Omega_{h}}\\
&\qquad+\tau\Ltwop{\edit{-}(1-S^{n-1})\mu_b\Phi^{n-1}C_b^{n-1}(1-\Phi^{n-1})}{\Psi_1}{\Omega_{h}},\\
&\Ltwop{C_s^n}{\Psi_2}{\Omega_{h}}+\tau\Ltwop{D_s({\Phi^{n})\nabla C_s^{n}}}{\nabla \Psi_2}{\Omega_{h}} =\Ltwop{C_s^{n-1}}{\Psi_2}{\Omega_{h}}+\tau\Ltwop{\kappa_sf_s(C_s^{n-1})(1-\Phi^{n-1})}{\Psi_2}{\Omega_{h}}\\
&\qquad+\tau\Ltwop{-((1-S^{n-1})\mu_s\Phi^{n-1}+\beta_s)C_s^{n-1}}{\Psi_2}{\Omega_{h}},\\
&\Ltwop{C_b^n(1-{\Phi^{n}})}{\Psi_3}{\Omega_{h}} =\Ltwop{{C_b^{n-1}}(1-\Phi^{n-1})}{\Psi_3}{\Omega_{h}}+\tau\Ltwop{\kappa_bf_b(\Phi^{n-1})(1-\Phi^{n-1})}{\Psi_3}{\Omega_{h}}\\
&\qquad+\tau\Ltwop{-\left(\left(1-S^{n-1}\right)\mu_b\Phi^{n-1}+\beta_b\right)C_b^{n-1}(1-\Phi^{n-1})}{\Psi_3}{\Omega_{h}}, \\
&\Ltwop{S^n}{\Psi_4}{\Omega_{h}} =\Ltwop{S^{n-1}}{\Psi_4}{\Omega_{h}}+\tau\Ltwop{-\delta_s(1-\Phi^{n-1})C_b^{n-1}S^{n-1}}{\Psi_4}{\Omega_{h}},
\end{aligned}
\end{equation}
for all $\Psi_i\in\Vho, i=1,2,3,4$.
The initial approximations at $t=0$ are taken to be the linear Lagrange interpolants of the initial data.  The scheme for model~\eqref{eqn:model_macro}, \eqref{IC_BC}  without matrix suitability corresponds to~\eqref{FE_scheme} with $S^n=0$ for all $n$, hence only three equations are solved per time step.  If greater accuracy is desired in terms of linearisation a fixed-point iteration or even a Newton-type linearisation could be employed as has been considered for semilinear reaction-diffusion systems in previous works, eg. \cite{madzvamuse2014fully}. Our previous studies of such systems suggest the IMEX method we employ above is robust as the diffusive terms, even when treated implicitly as above,  necessitate the use of a small time-step \cite{lakkis2013implicit}.

\section{Model parametrisation and simulation results}\label{sec:numerics}
\subsection{Parametrisation of the model equations and choice of discretisation parameters}\label{sec:param}
To facilitate comparison  with the results in~\cite{deakin2013mathematical},  
for the parameters in the macroscopic models~\eqref{eqn:model_macro}, \eqref{IC_BC} and \eqref{het_ecm}, \eqref{IC_BC} in section~\ref{sec:macro} we use the parameter values of~\cite{deakin2013mathematical}, see Table~\ref{table_1}. It is a worthwhile direction for future work to investigate the sensitivity of the numerical solutions to changes in the parameter values and to understand whether the broad conclusions are robust to variations in the parameters. We checked that  the simulation results reported on below remain qualitatively indistinguishable from those obtained using  numerical parameters (mesh-size and where relevant time step) approximately double those used to generate the simulations reported on in the manuscript.
\begin{table}
\begin{center}
\begin{tabular}{ |c|c|c| c|}
 \hline
 & non-dimensional values & original values & references \\
 \hline
 $\Omega$ & $(0,2)^d$, \, $d=2,3$ & $0.1 {\rm cm}$  & \cite{collier2011} \\
 \hline
 $\kappa_s$ & $4$ &  $4\cdot 10^{-4} {\rm s}^{-1}$ & estimated, \cite{deakin2013mathematical} \\
 \hline
 $\kappa_b  $ & $5$  & $5\cdot 10^{-4} {\rm s}^{-1}$  & estimated, \cite{deakin2013mathematical}\\
 \hline
 $D_s$  &  $1.29\cdot 10^{-2}$  & $1.29 \cdot 10^{-8} \,  {\rm cm}^2 {\rm s}^{-1} $ & \cite{collier2011} \\
 \hline
 $\mu_b$, $\mu_s$ & $ 1$ & $-$ &\cite{anderson2000}\\
 \hline
 $\mu_b$, $\mu_s$ & $8.15$ & $-$ & \cite{lolas2005}\\
 \hline
 $\beta_s$, $\beta_b$   & $0.1$ & $1\cdot 10^{-5} {\rm s}^{-1} $ &  estimated, \cite{anderson2000, deakin2013mathematical}
 \\
  \hline
 $\delta_s$   & $1$ & $ - $ &  \cite{deakin2013mathematical}
 \\
 \hline
\end{tabular}
\end{center}
\caption{Parameter values used for numerical simulations for models \eqref{eqn:model_macro}, \eqref{IC_BC} and \eqref{het_ecm}, \eqref{IC_BC}.} \label{table_1}
\end{table} 

\subsection{Approximation of the effective diffusion tensor for soluble MMPs}\label{sec:homg_diff_sims}
As discussed in section \ref{sec:macro} we obtain an effective diffusion tensor for different volume fractions of the ECM by  solving   problems~\eqref{eqn:cell_problems} on different geometries corresponding to different volume fractions of ECM and fitting, using least squares, a polynomial  based on these values. 

\begin{figure}[htbp]
  \includegraphics[trim = 10mm 20mm 10mm 20mm,  clip, width=0.75\linewidth]{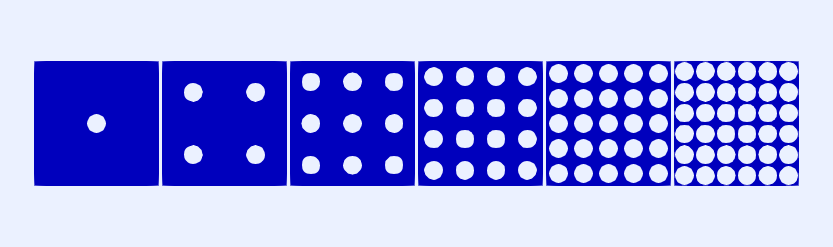}
\caption{Domains used for computing the effective diffusion tensor in $2{\rm d}$ corresponding to volume fractions of ECM of $0.982, 0.929, 0.841, 0.717, 0.558$ and $0.364$ reading from left to right. }
  \label{fig:hom_domain_2d}
  \end{figure}

\begin{figure}[htbp]
  \includegraphics[trim = 50mm 0mm 50mm 0mm,  clip, width=0.24\linewidth]{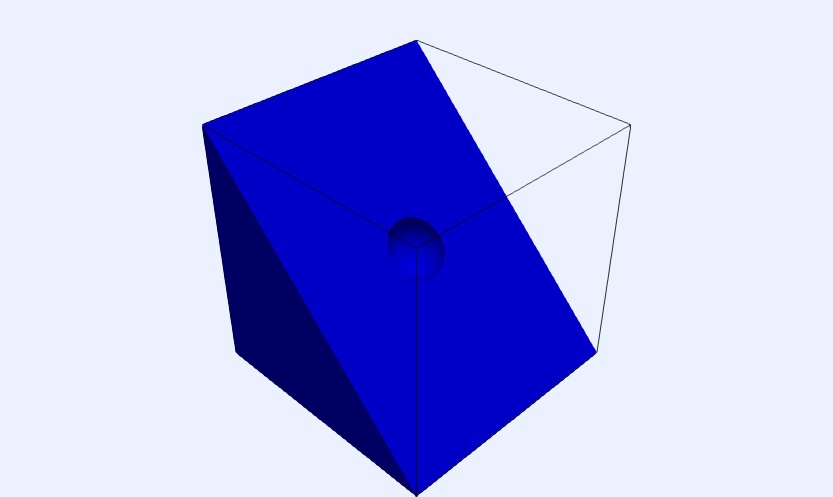}
  \includegraphics[trim = 50mm 0mm 50mm 0mm,  clip, width=0.24\linewidth]{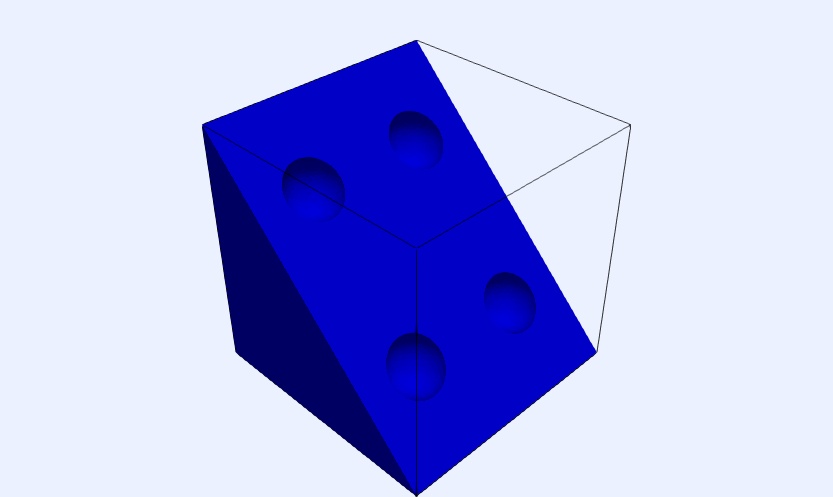}
  \includegraphics[trim = 50mm 0mm 50mm 0mm,  clip, width=0.24\linewidth]{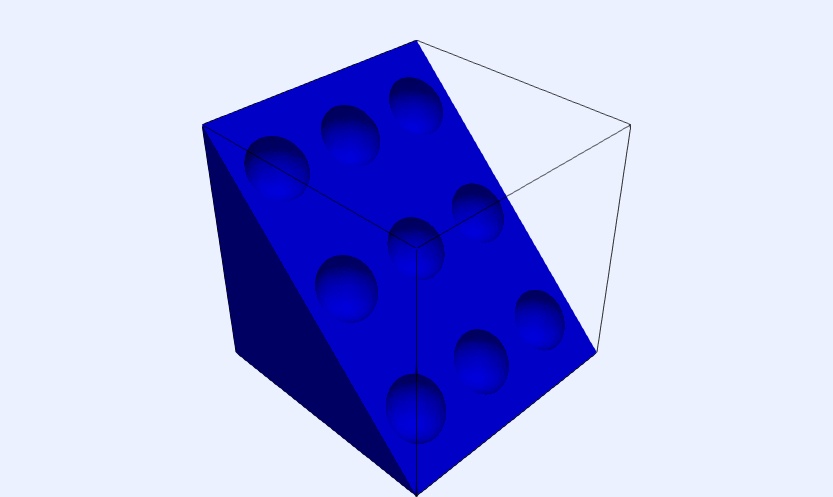}
  \includegraphics[trim = 50mm 0mm 50mm 0mm,  clip, width=0.24\linewidth]{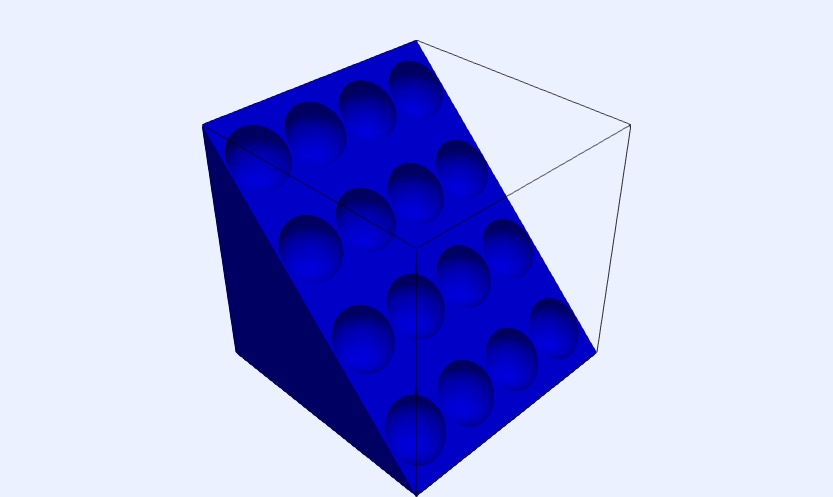}
\caption{Domains used for computing the effective diffusion tensor in $3{\rm d}$ corresponding to volume fractions of ECM of $0.993, 0.942, 0.805$ and $0.537$ reading from left to right.  Half of the domain has been made transparent in each case to visualise the interior boundaries.}
  \label{fig:hom_domain_3d}
  \end{figure}

In the $2{\rm d}$ case the domains we consider correspond to squares of the form $[0,L]^2$ with $n^2$ circle(s) of radius $r=\frac{3}{40}L$, with $n=1,2,3,4,5,6$, removed, i.e., the cells are circular. The centres of the circles are taken to form a square lattice, with the resulting geometry as shown in Figure \ref{fig:hom_domain_2d}. The corresponding volume fractions of ECM are $0.982, 0.929, 0.841, 0.717, 0.558$ and $0.364$ for $n=1,2,3,4,5$ and $6$ respectively.  
  For the $3{\rm d}$ case we adopt a similar geometric setup to the $2{\rm d}$ case using cubes of the form $[0,L]^3$ with $n^3$ spheres of radius  $r=\frac{3}{25}L$, with $n=1,2,3,4$, removed, i.e., the cells are spherical. The centres of the spheres are taken to form a cubic lattice with the resulting geometry shown in Figure~\ref{fig:hom_domain_3d}. The corresponding volume fractions of ECM are $0.993, 0.942, 0.805$ and $0.537$ for $n=1,2,3$ and $4$ respectively. 
  
\begin{figure}[htbp]
  \includegraphics[trim = 0mm 0mm 0mm 0mm,  clip, width=0.45\linewidth]{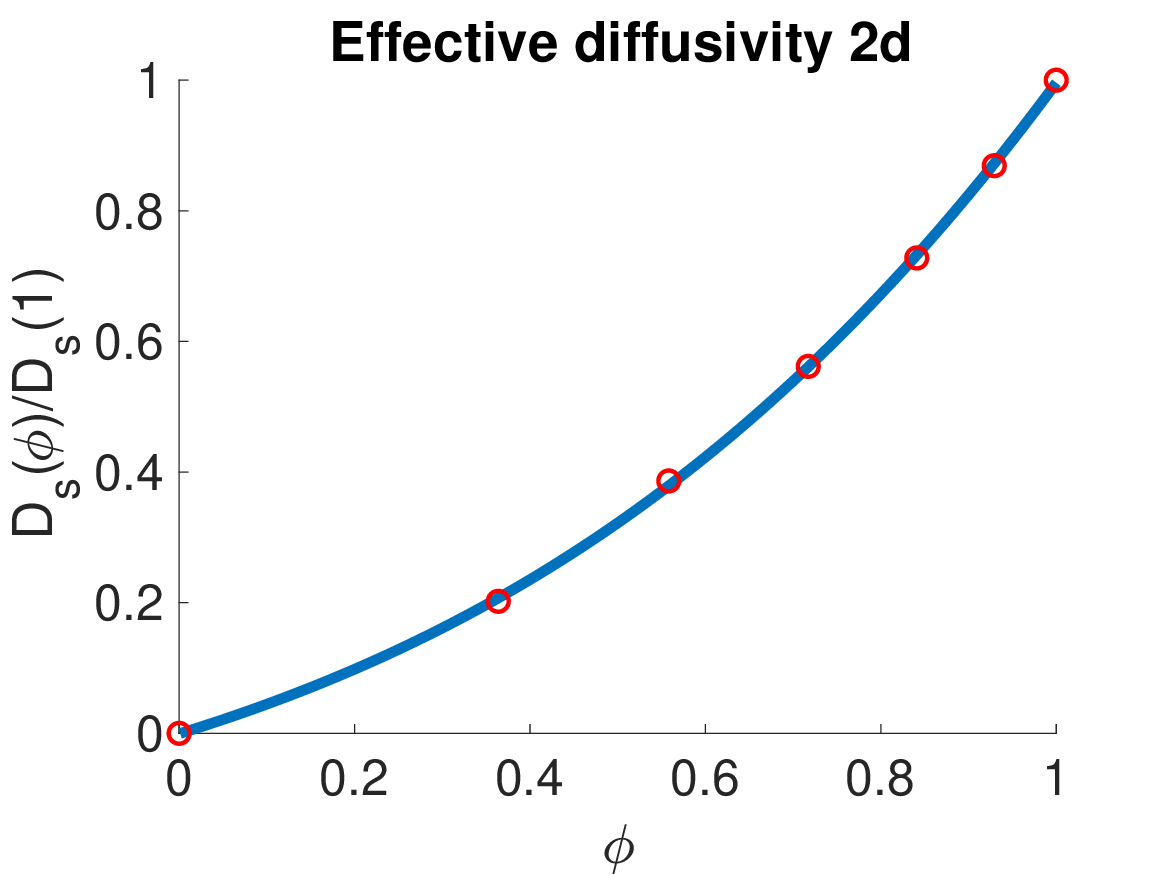}
  \includegraphics[trim = 0mm 0mm 0mm 0mm,  clip, width=0.45\linewidth]{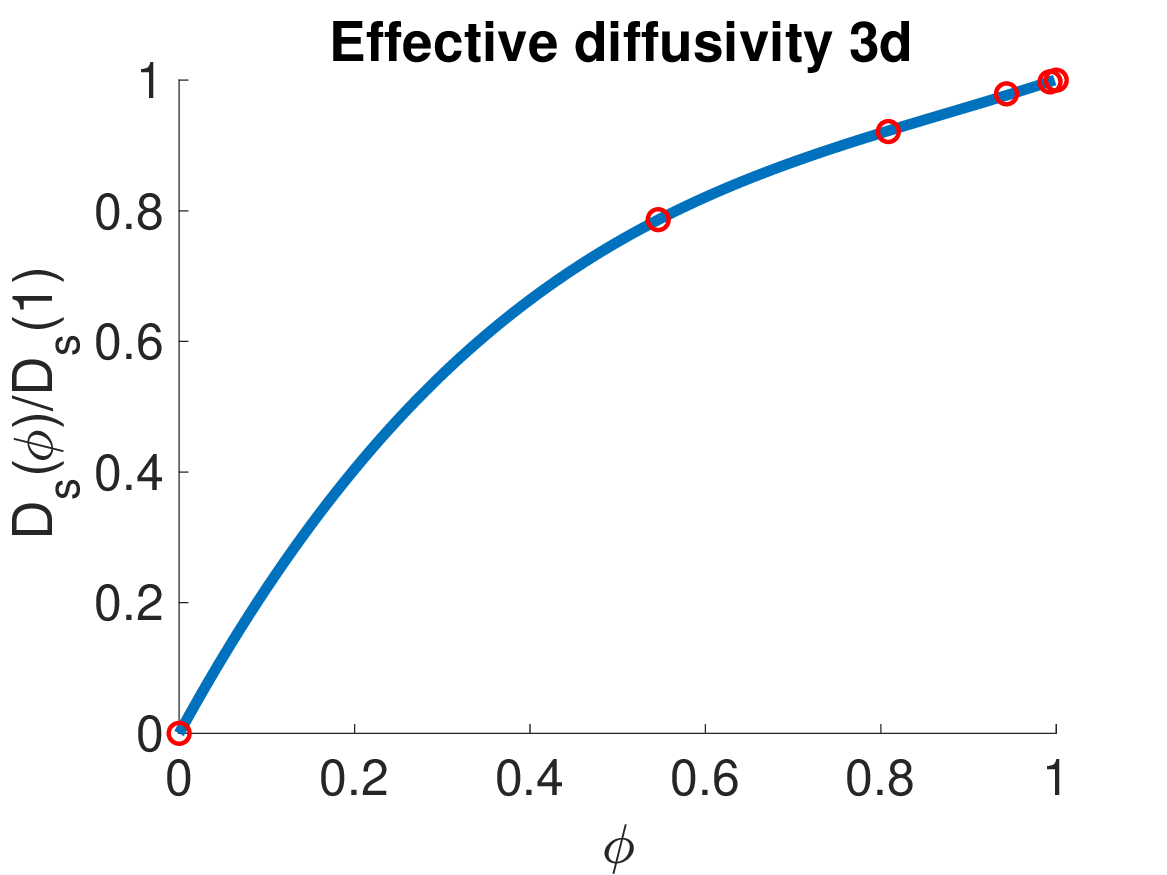}
\caption{Fitted effective diffusivity in $2{\rm d}$ and $3{\rm d}$. The red circles correspond to computed effective diffusivities obtained by solving \eqref{eqn:cell_problems} on domains described in section \ref{sec:homg_diff_sims}. The blue curves are  least squares cubic polynomial fits using these values.}
  \label{fig:d_eff}
  \end{figure}

 For the choices of microscopic geometry described above, we note that, due to symmetry, the effective diffusion tensor will be  isotropic and diagonal.  This is reflected in the numerical results, for which an unstructured triangulation is used, in that the off diagonal components of the computed diffusion tensors for all geometries are  two orders of magnitude smaller than the diagonal components and the diagonal components of the computed diffusion tensors on each specific geometry exhibit only small differences ($<1\%$). We therefore obtain values for fitting an effective diffusion coefficient by  assuming we have a $\phi$ dependent scaling of the diffusion tensor and take the average of the diagonal components of the computed diffusion tensors to approximate this scaling on each microscopic geometry.   We take $D_s(0)=0$ as a further fitting point as we assume the soluble MMPs do not diffuse through the cancer cells. \edit{We use a finite element method implemented in the FEniCS software~\cite{logg2012automated} for the approximation of \eqref{eqn:scheme_diff_tensor}.  In the simplest setting, where the effective diffusion tensor is diagonal and isotropic, we compute a least squares polynomial fit  to obtain the effective diffusivity in $2{\rm d}$ and $3{\rm d}$ (in the examples below a cubic polynomial appears adequate to fit the computed data).
 The fitted polynomials obtained are given by}   
 \begin{align}
 {D_s(\phi)}=\big(0.25\phi^3+0.33\phi^2+0.42\phi\big){D_s(1)},
 \end{align}
 in $2{\rm d}$ and
  \begin{align}
 {D_s(\phi)}=\big(0.89\phi^3-2.35\phi^2+2.46\phi\big){D_s(1)},
 \end{align}
 in $3{\rm d}$. Figure \ref{fig:d_eff} shows plots of the effective diffusivity as  the volume fraction of ECM changes. 
We note that the effective diffusivity appears to be a convex function of the volume fraction in $2{\rm d}$ but not in $3{\rm d}$. \edit{We note that heuristic formulae for the effective diffusivity of a porous medium as a function of volume fraction of the medium versus void are given in many previous works e.g., \cite{weissberg1963effective,du2009modelling} where various (convex) functions are proposed. We note that more complex dependence on volume fraction including  concave functions  are computed in  \cite{doi:10.1021/acs.jpcc.9b03250}[c.f., Figures~11 and~14] albeit in a significantly more complicated setup. 
The result above  may also be linked to interesting results in percolation theory such that the critical exponents that relate diffusivity of the medium with volume fraction are dimension dependent \cite{bunde2005diffusion}. This is an interesting direction for future work.
To investigate further how the effective diffusivity varies depending on the choice of microstructure, in Appendix~\ref{app:square_eff_diff}, we consider a different geometric setup for the microstructure and  the effective diffusivities computed are quantitatively comparable with those presented above.  For  the general case of a  non-diagonal effective diffusion tensor, simple  polynomial fits or affine interpolation may not be sufficient  as positive definiteness can in general not be ensured. This is discussed in Appendix~\ref{app:full_eff_diff} wherein geometries that yield a full effective diffusion tensor are considered. An interesting area for future work would be to consider if those results are robust to irregular shape and distribution of cells within the microstructure.}


 \subsection{Simulation of the macroscopic invasion model \eqref{eqn:model_macro}, \eqref{IC_BC}}\label{subsec:macro_model_numerics} 
We report on simulations of the macroscopic invasion model \eqref{eqn:model_macro}, \eqref{IC_BC} using the parameters in Table~\ref{table_1} in section~\ref{sec:param}. To understand the relative contributions of membrane bound and soluble MMPs we run three simulation \edit{scenarios}  in $2{\rm d}$ and $3{\rm d}$, i.e., six in total.  One where both soluble and membrane bound MMPs degrade the ECM, one where only soluble MMPs degrade the ECM ($\mu_b=0$) and one where only membrane bound MMPs degrade the ECM ($\mu_s=0$). 

For the simulations in $2{\rm d}$ we take $\Omega=[-1,1]^2$ and the initial condition 
\begin{equation}\label{phi_ic_ellipse}
\phi_0(x)=\Big(1-e^{-\left((4x_1)^2+(8x_2)^2\right)}\Big),
\end{equation}
and the initial data for the remaining variables are set to $0$. We take the timestep $\tau=10^{-2}$ and use a mesh with 12909 degrees of freedom (DOFs).
\begin{figure}[htbp]
  \begin{subfigure}[t]{0.49\textwidth}
  \includegraphics[trim = 0mm 200mm 0mm 30mm,  clip, width=\linewidth]{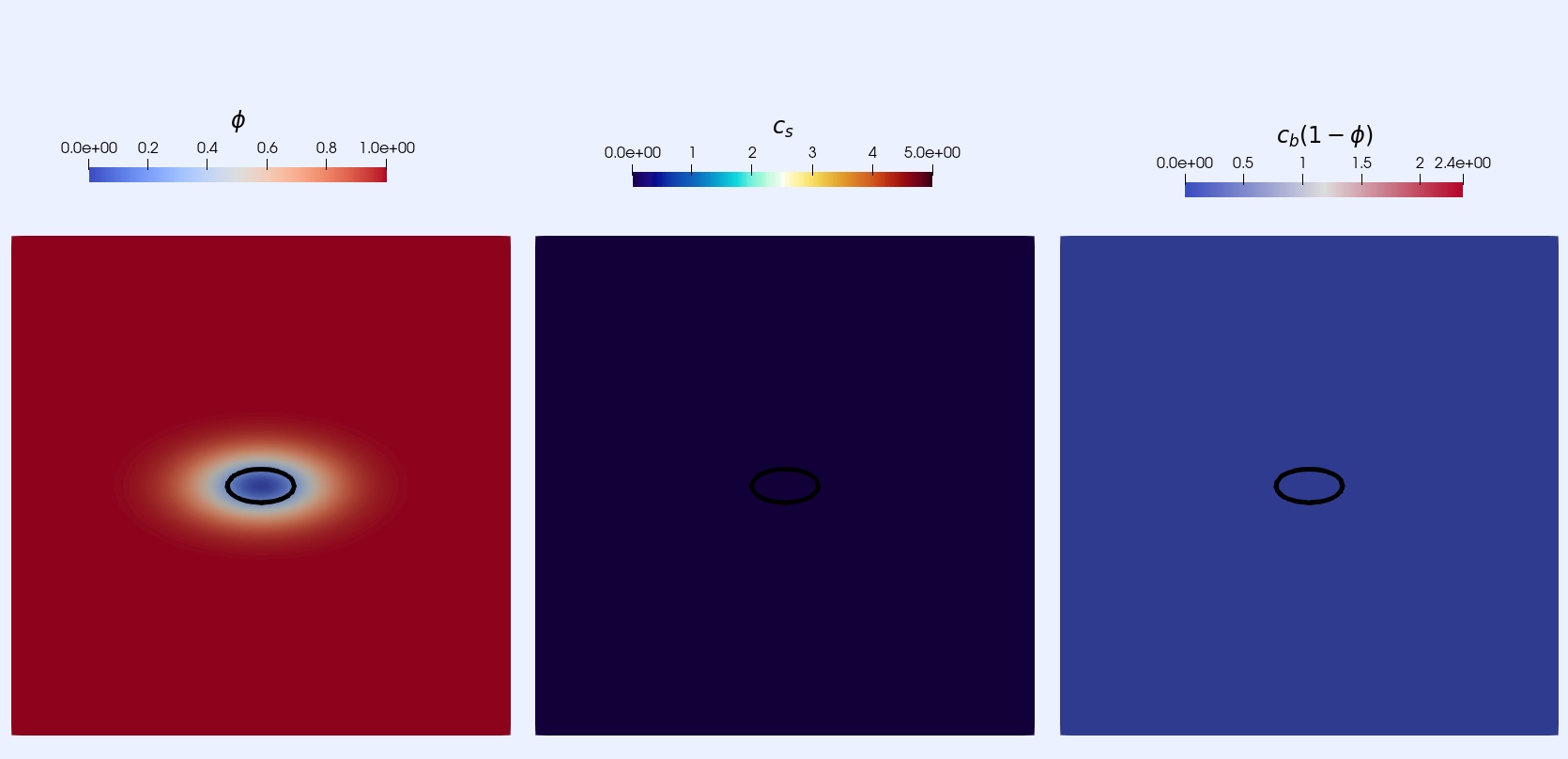}
  \includegraphics[trim = 0mm 0mm 0mm 80mm,  clip, width=\linewidth]{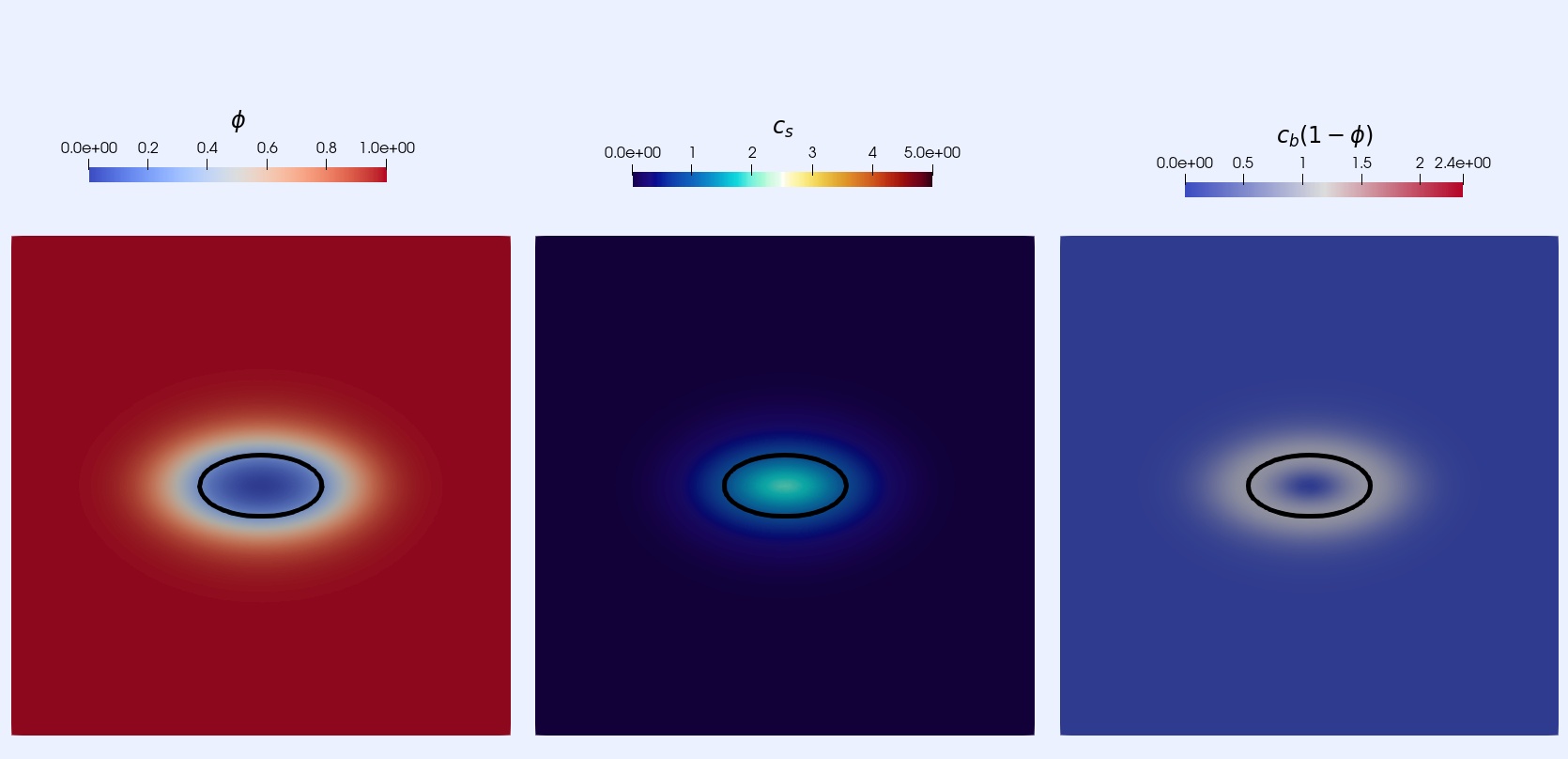}
  \includegraphics[trim = 0mm 0mm 0mm 80mm,  clip, width=\linewidth]{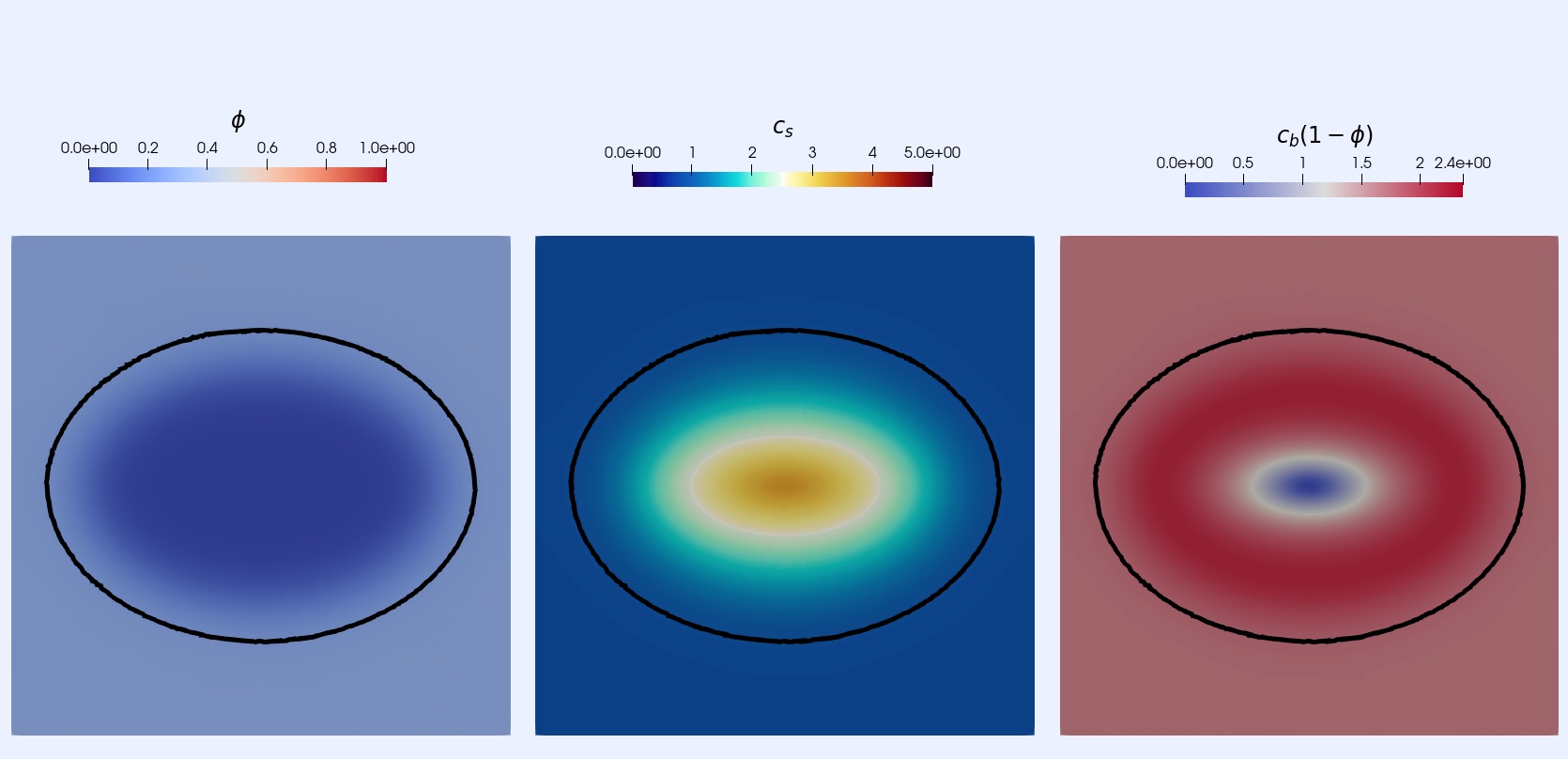}
  \includegraphics[trim = 0mm 0mm 0mm 80mm,  clip, width=\linewidth]{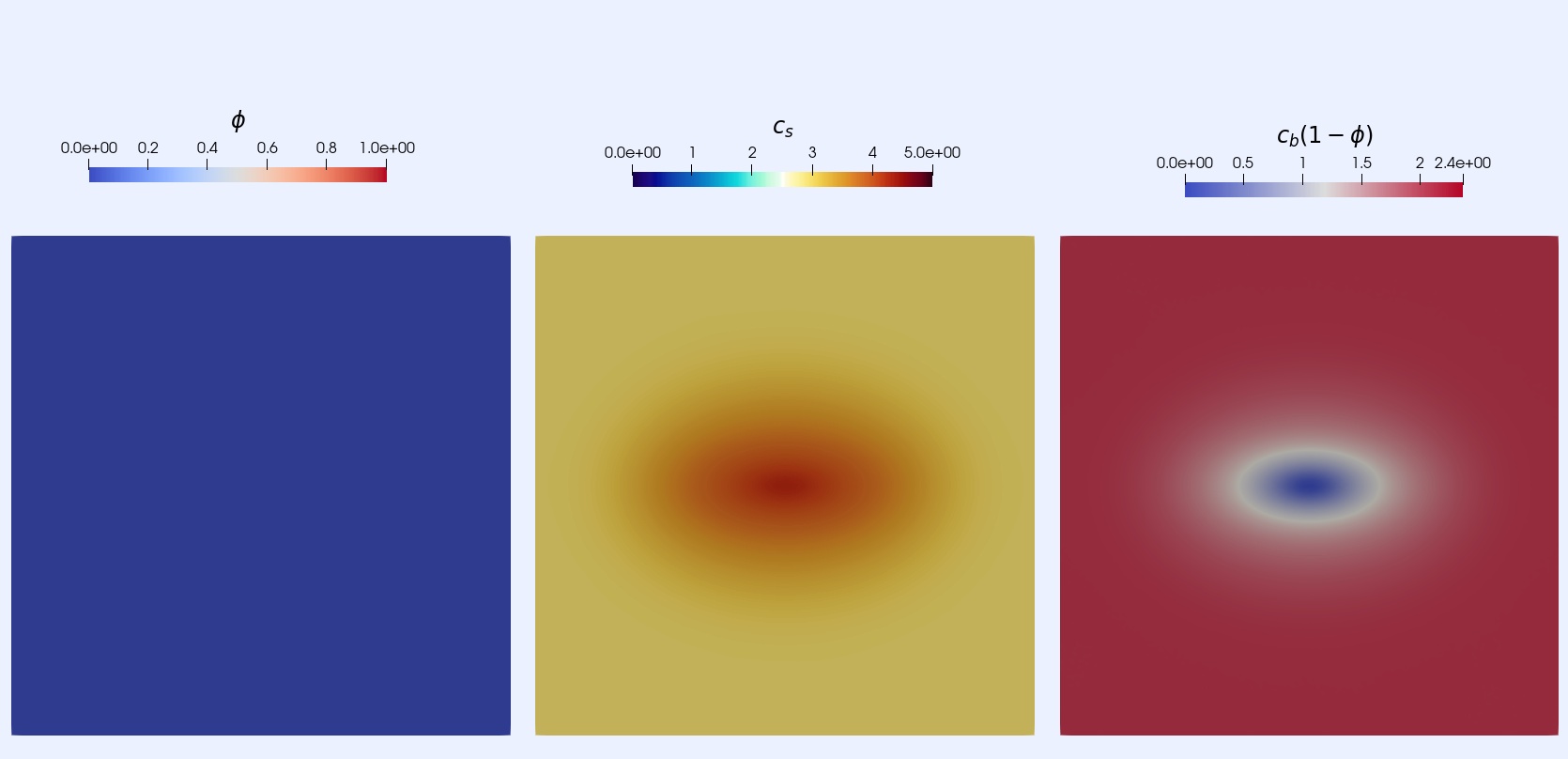}
\caption{$\mu_b=\mu_s=1$. Note for the results in this subfigure both soluble and bound MMPs degrade the matrix.}
  \label{fig:2d_both}
    \end{subfigure}
\hfill
  \begin{subfigure}[t]{0.49\textwidth}
  \includegraphics[trim = 0mm 200mm 0mm 30mm,  clip, width=\linewidth]{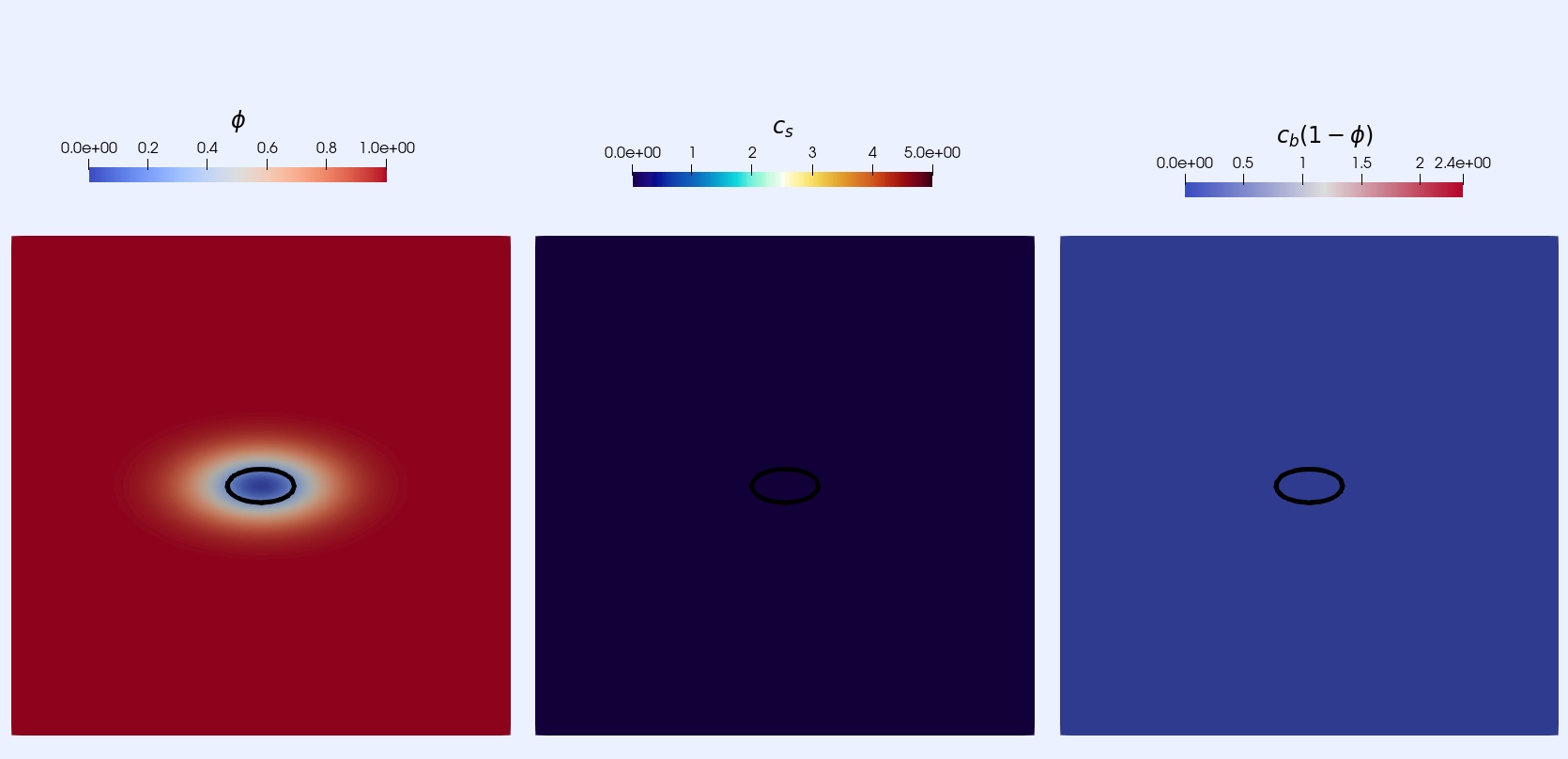}
  \includegraphics[trim = 0mm 0mm 0mm 80mm,  clip, width=\linewidth]{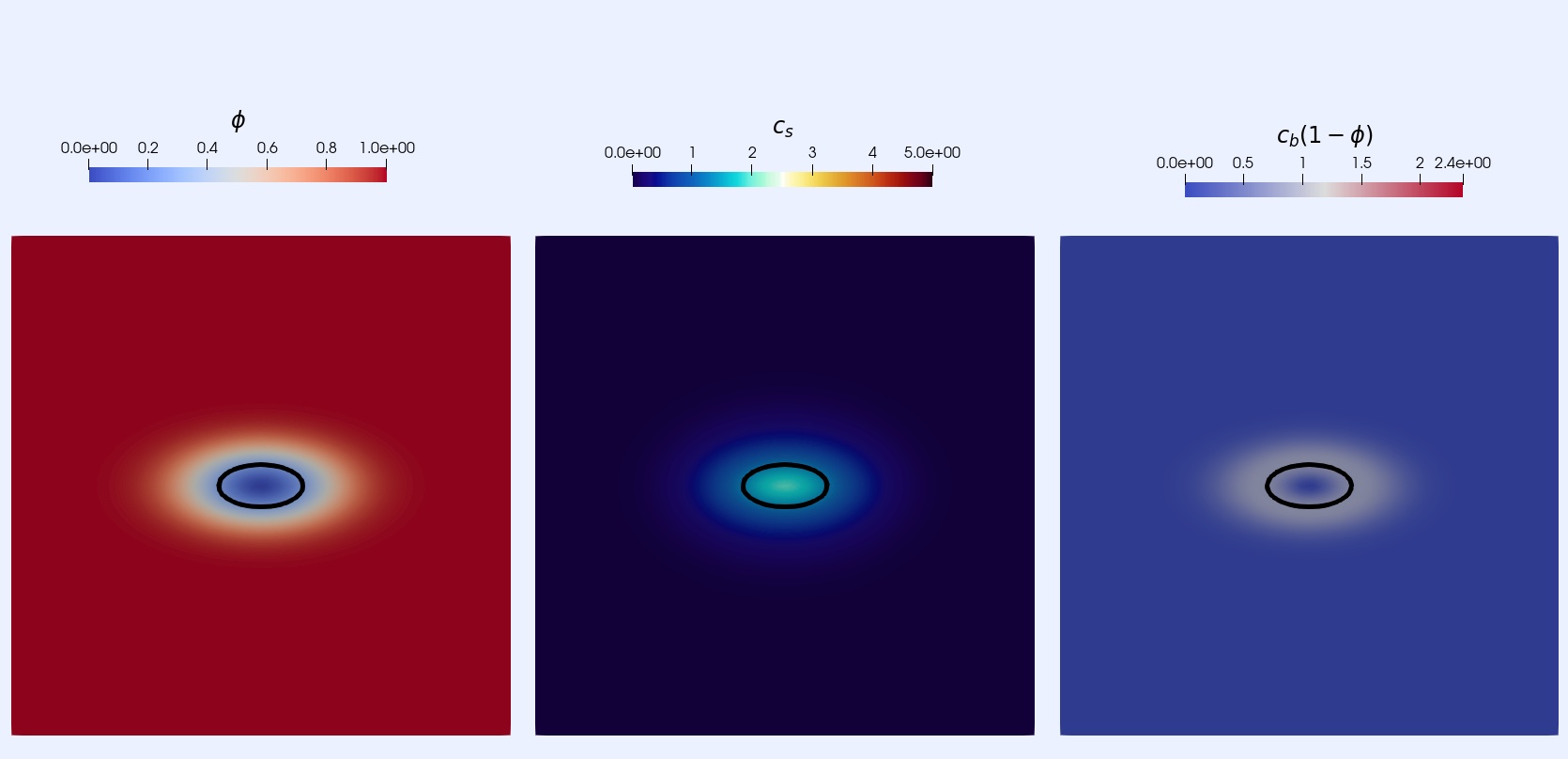}
  \includegraphics[trim = 0mm 0mm 0mm 80mm,  clip, width=\linewidth]{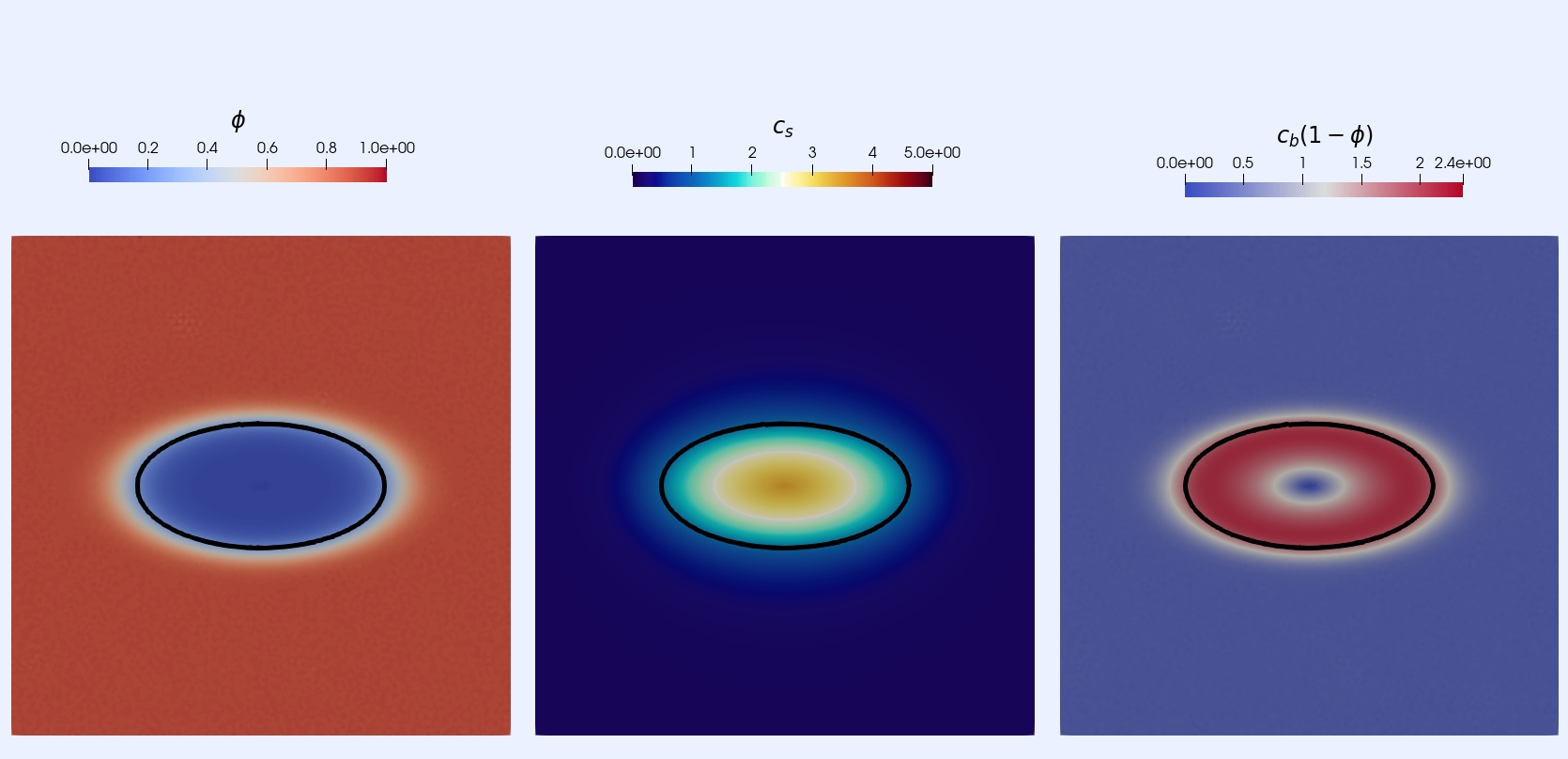}
  \includegraphics[trim = 0mm 0mm 0mm 80mm,  clip, width=\linewidth]{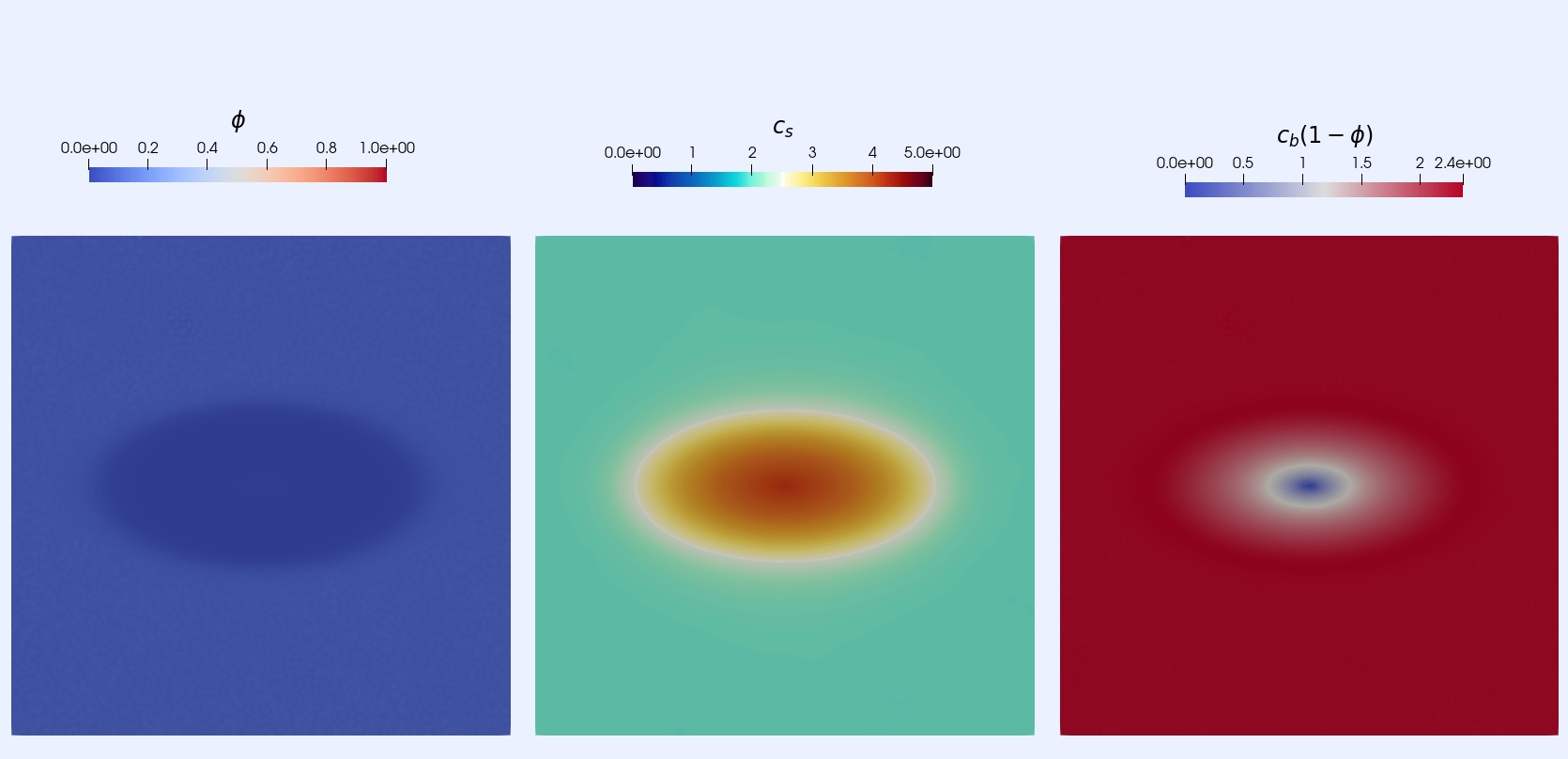}
\caption{$\mu_b=1,\mu_s=0$. Note for the results in this subfigure the soluble MMPs do not influence $\phi$.}
  \label{fig:2d_mub}
    \end{subfigure}
\hfill
  \begin{subfigure}[t]{0.49\textwidth}
  \includegraphics[trim = 0mm 200mm 0mm 30mm,  clip, width=\linewidth]{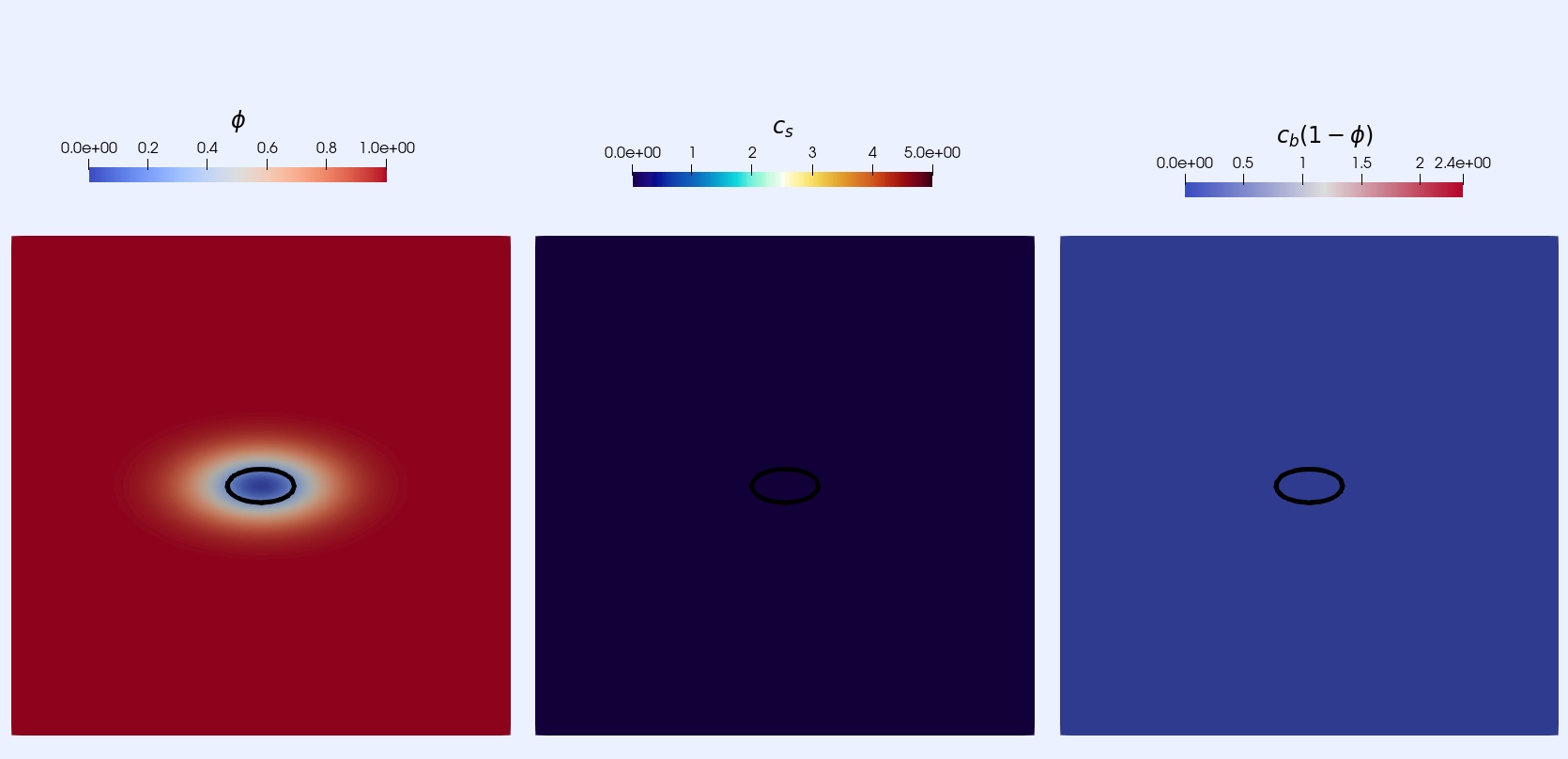}
  \includegraphics[trim = 0mm 0mm 0mm 80mm,  clip, width=\linewidth]{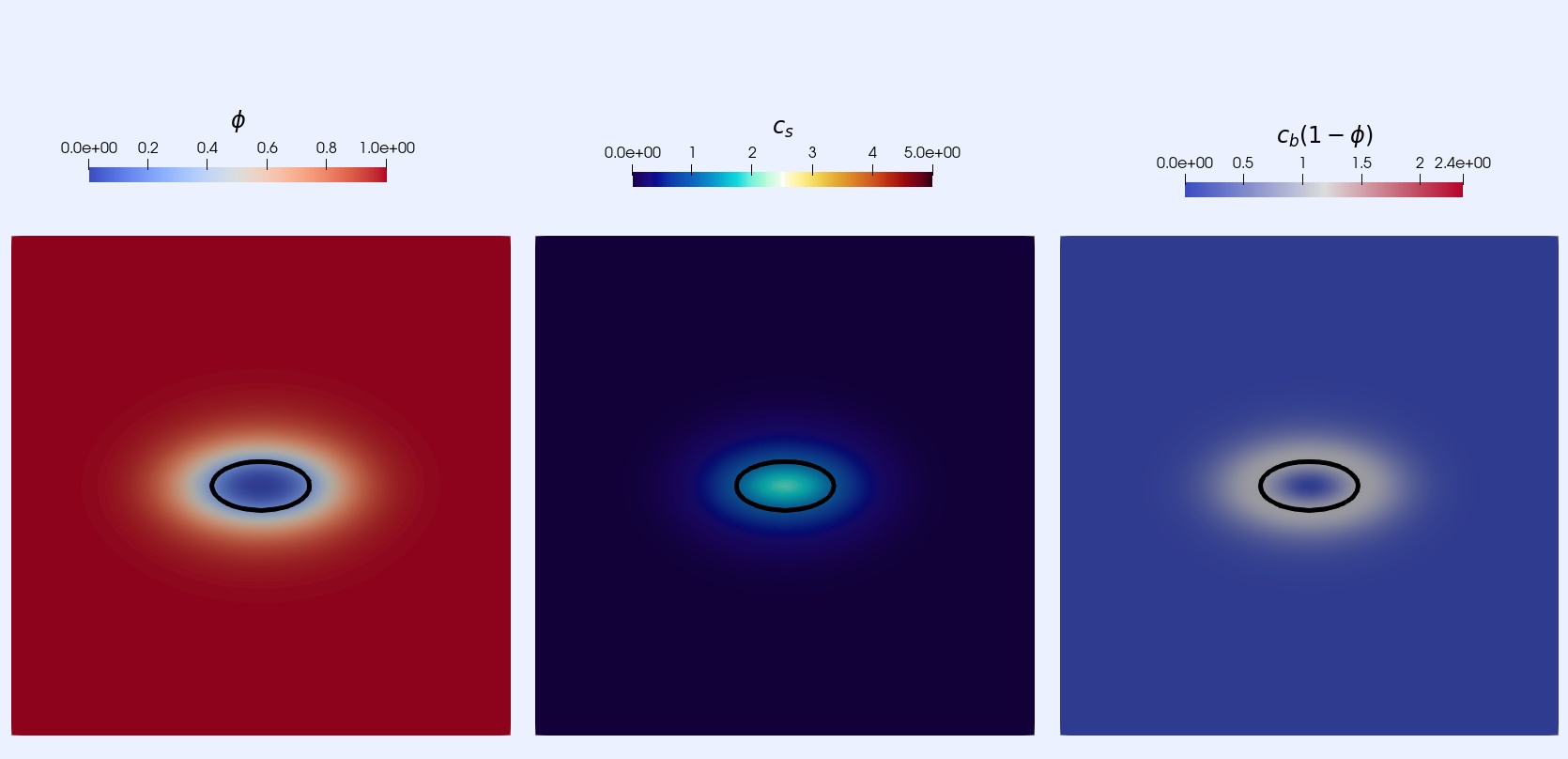}
  \includegraphics[trim = 0mm 0mm 0mm 80mm,  clip, width=\linewidth]{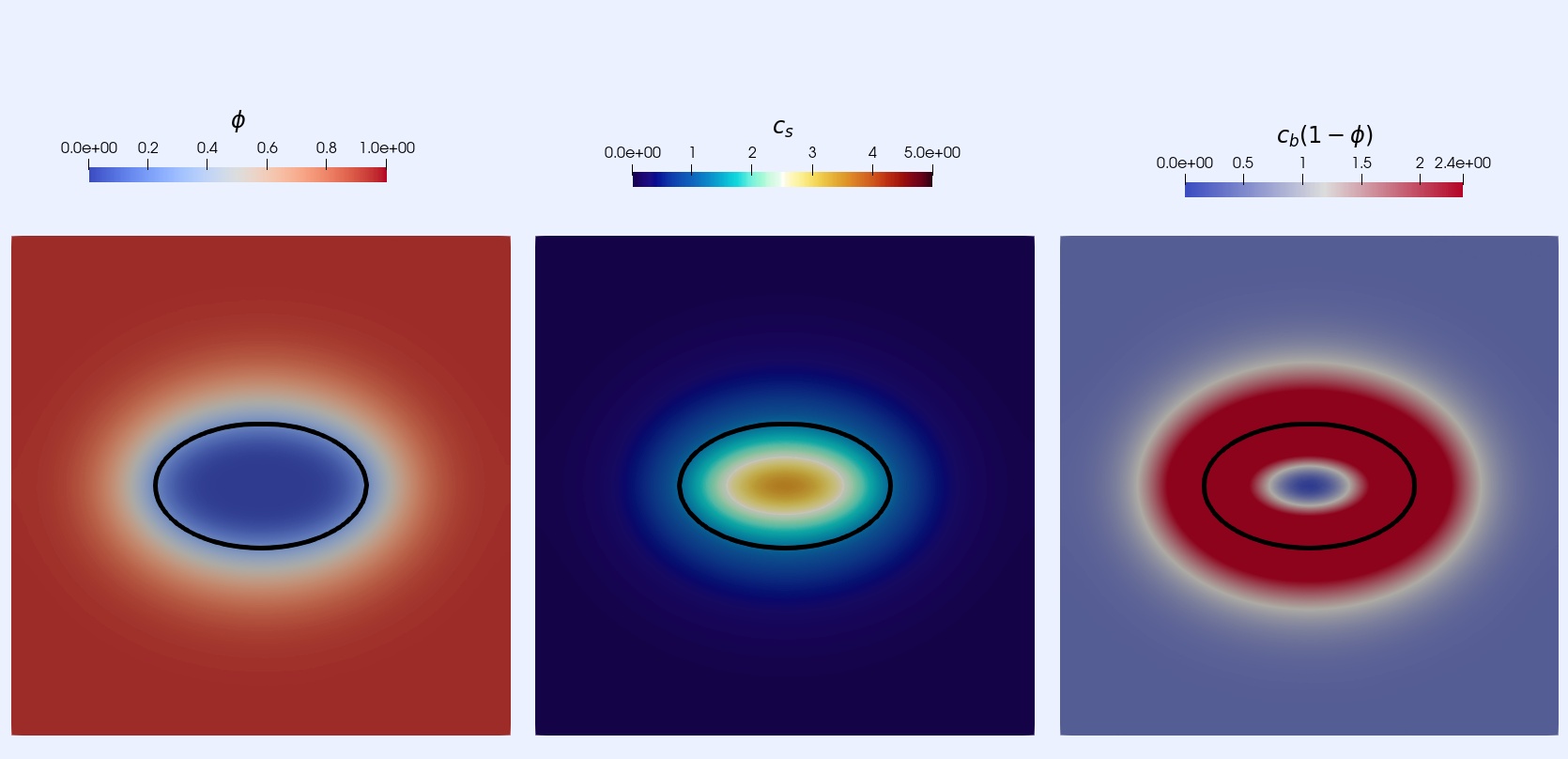}
  \includegraphics[trim = 0mm 0mm 0mm 80mm,  clip, width=\linewidth]{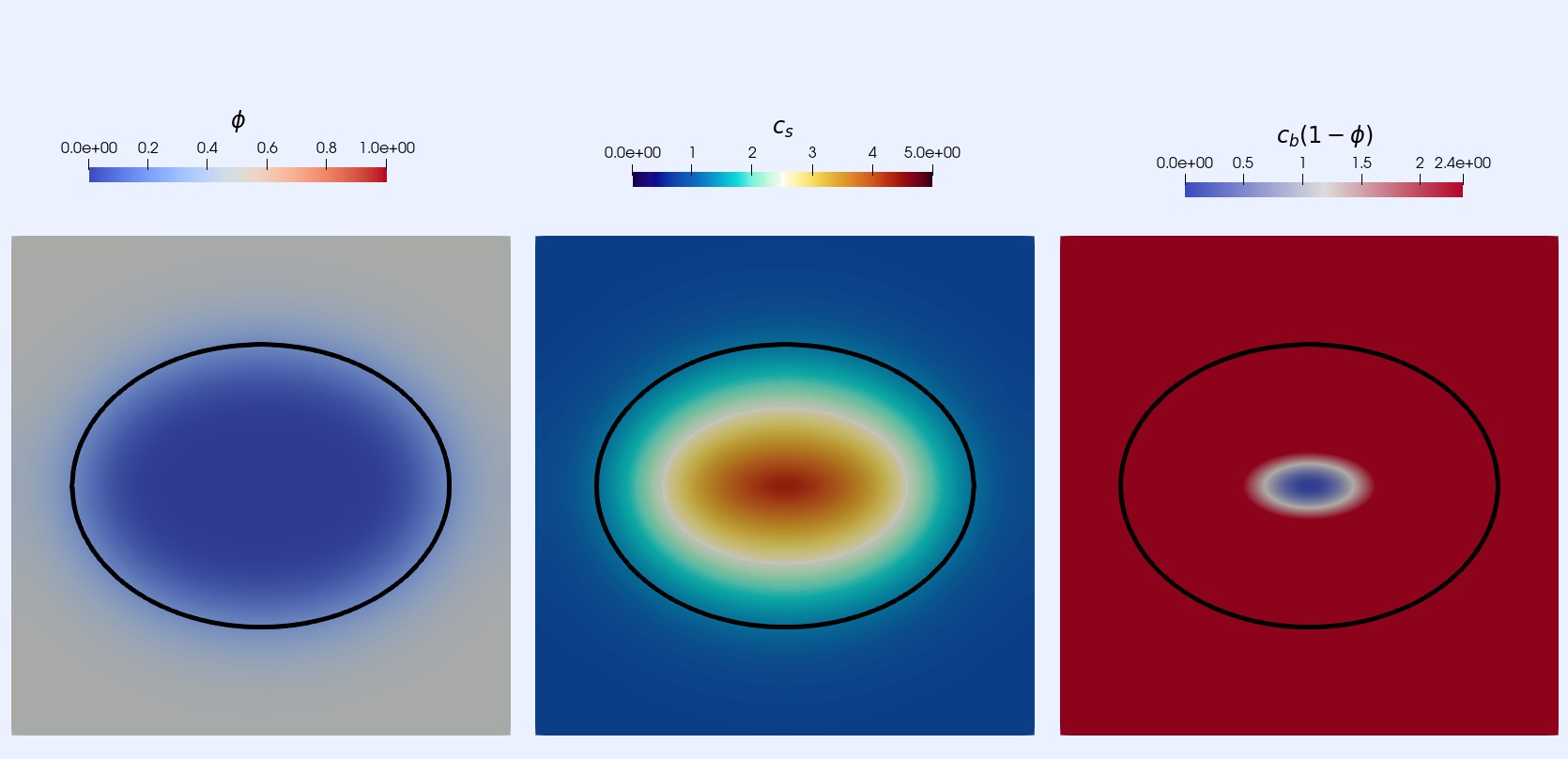}
\caption{$\mu_b=0,\mu_s=1$. Note for the results in this subfigure the bound MMPs do not influence $\phi$.}
  \label{fig:2d_mus}
    \end{subfigure}
    \caption{Simulations of the invasion model \eqref{eqn:model_macro}, \eqref{IC_BC} in $2{\rm d}$. In each subfigure, each row corresponds to times $t=1, 3$ and $5$. The black contour indicates the level set $\phi=0.25$ corresponding to a volume fraction of $75\%$ cancer cells. The effect of membrane bound MMPs  increasing the speed of invasion is apparent whilst the degradation of the ECM by the soluble MMPs appears to generate more radially symmetric invasive profiles. Parameter values as in Table~\ref{table_1}.}
    \label{fig:2d}
  \end{figure}
Figure \ref{fig:2d} shows the results of simulations in $2{\rm d}$. We observe MMP mediated invasion \edit{by cancer cells}. In the case where both membrane bound and soluble MMPs degrade the ECM, see Figure~\ref{fig:2d_both}, we  observe rapid invasion of the ECM and, at least at the initial stages of invasion, the membrane bound MMPs remain localised to a region close to level set $\phi=0.25$.  The soluble MMPs  are primarily located in the region where the volume fraction of cancer cells is high and the influence of slowed transport into the ECM due to the $\phi$ dependence of the diffusivity is evident. In the simulations where the soluble MMPs are unable to degrade the ECM, see Figure~\ref{fig:2d_mub},  we still observe relatively rapid invasion, albeit with reduced speed compared to the case where both MMP species degrade the ECM. The profile of the invasive front appears to match the initial shape of the invasive front.  
In the simulations where the bound MMPs are unable to degrade the ECM, see Figure~\ref{fig:2d_mus},  we observe the slowest invasive speed of all the simulations reported  in Figure~\ref{fig:2d}. At the final time $t=5$ this is the only case in which the entirety of the domain does not consist of at least $75\%$ cancer cells. The invasive  front appears to take on a more radially symmetric shape similar to the case when both species of MMP degrade the matrix. As a whole the results depicted in Figure~\ref{fig:2d} illustrate that the model reflects behaviour observed in experiments such as \cite[Figure 1A]{sabeh2009protease}, i.e.,  the bound MMPs are more significant for \edit{cancer} invasion into the ECM via matrix degradation than the soluble MMPs.

For the simulations in $3{\rm d}$ we take $\Omega=[-1,1]^3$ and the initial condition 
\begin{equation}\label{phi_ic_ellipse}
\phi_0(x)=\Big(1-e^{-\left((4x_1)^2+(4x_2)^2+(8x_3)^2\right)}\Big),
\end{equation}
and the initial data for the remaining variables are set to $0$. We take the time step $\tau=10^{-2}$ and use a mesh with $61905$  DOFs.
\begin{figure}[htbp]
\centering
  \begin{subfigure}[t]{0.49\textwidth}
  \includegraphics[angle=270,trim = 300mm 130mm 0mm 0mm,  clip, width=0.32\linewidth]{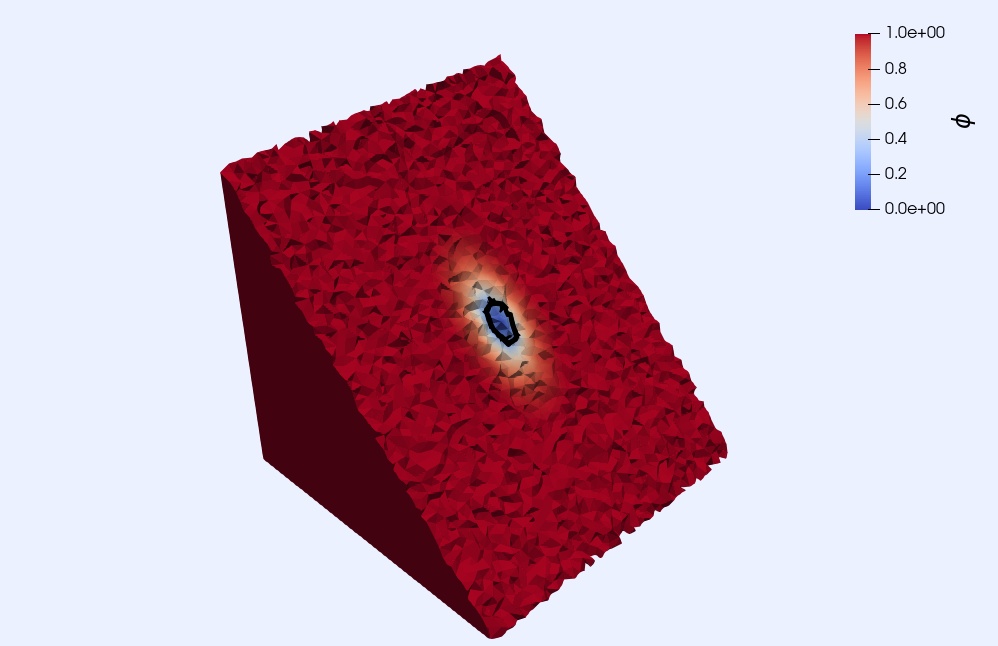}
  \includegraphics[angle=270,trim = 300mm 130mm 0mm 0mm,  clip, width=0.32\linewidth]{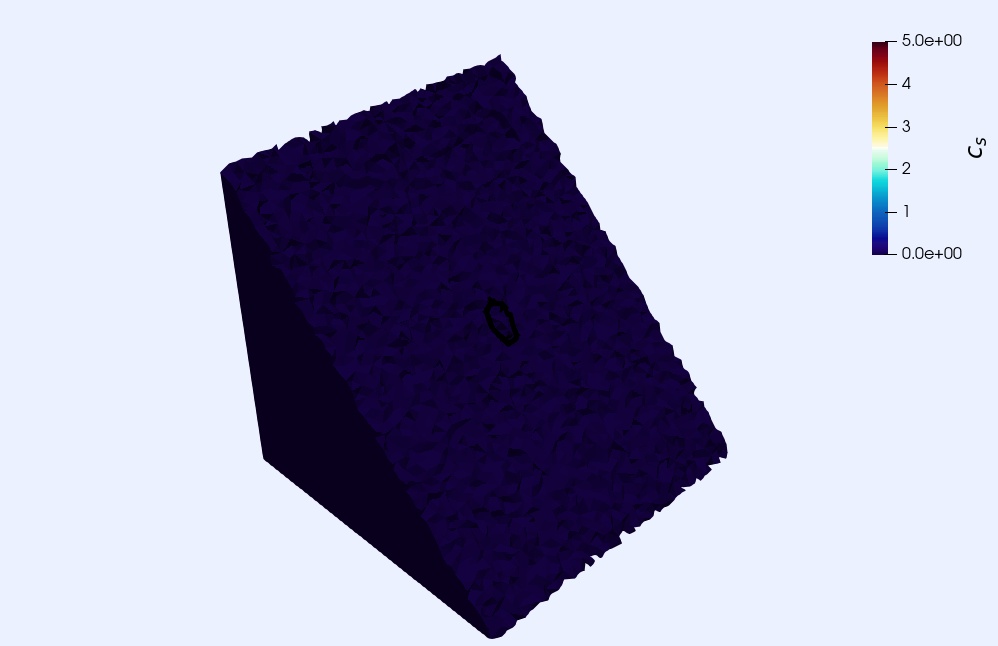}
  \includegraphics[angle=270,trim = 300mm 130mm 0mm 0mm,  clip, width=0.32\linewidth]{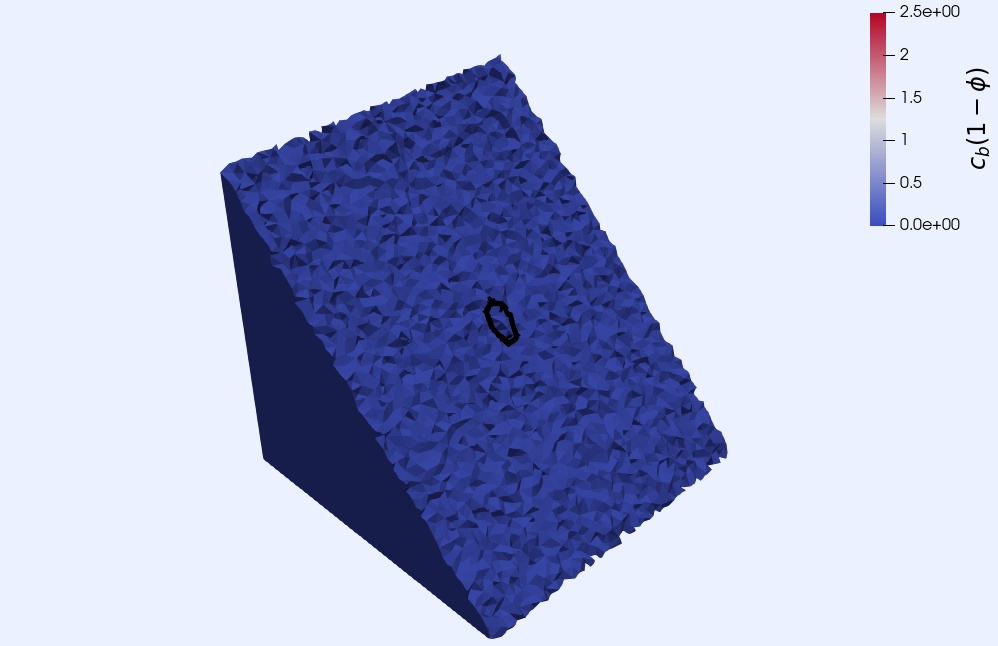}
  \includegraphics[trim = 70mm 0mm 90mm 18mm,  clip, width=0.32\linewidth]{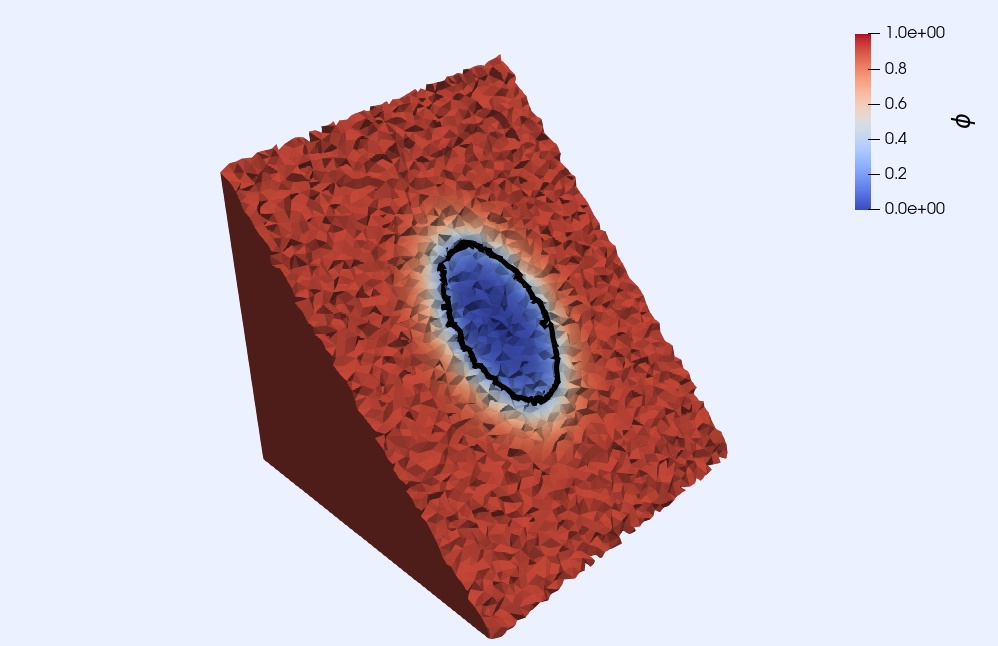}
  \includegraphics[trim = 70mm 0mm 90mm 18mm,  clip, width=0.32\linewidth]{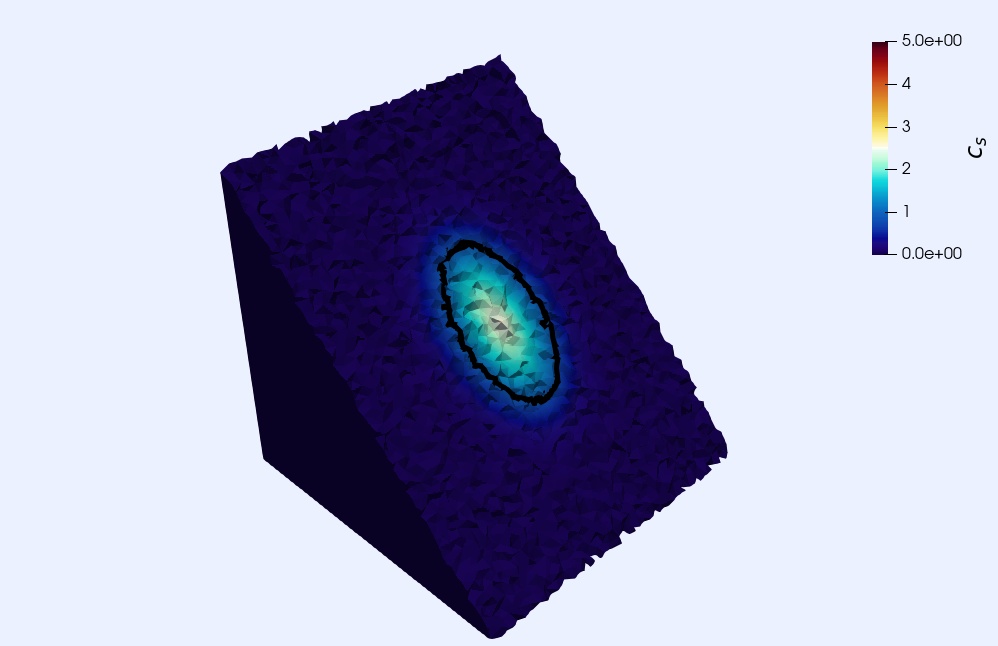}
  \includegraphics[trim = 70mm 0mm 90mm 18mm,  clip, width=0.32\linewidth]{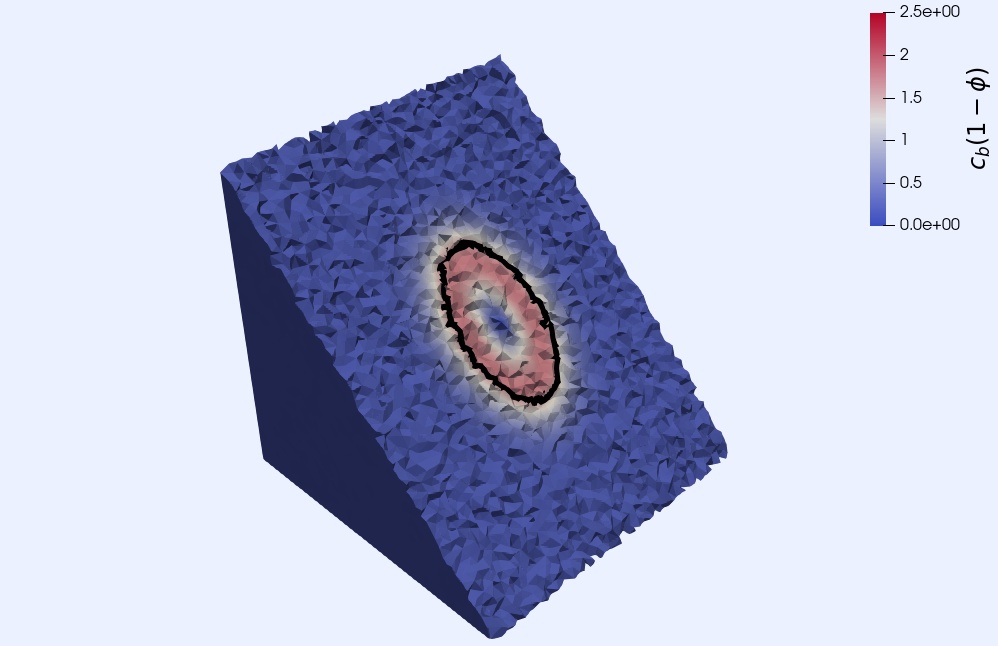}
    \includegraphics[trim = 70mm 0mm 90mm 18mm,  clip, width=0.32\linewidth]{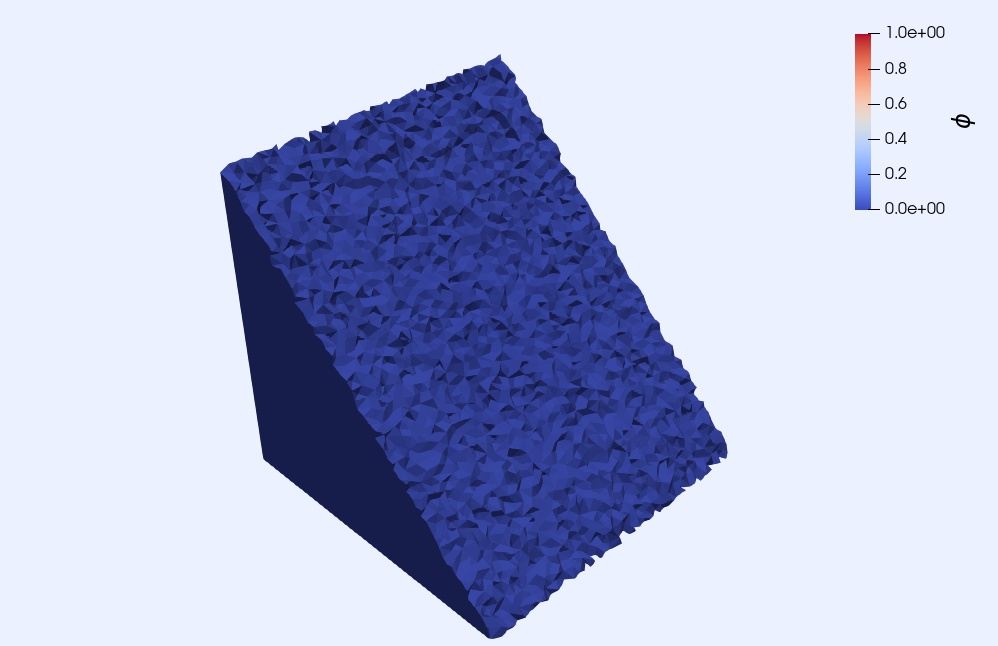}
  \includegraphics[trim = 70mm 0mm 90mm 18mm,  clip, width=0.32\linewidth]{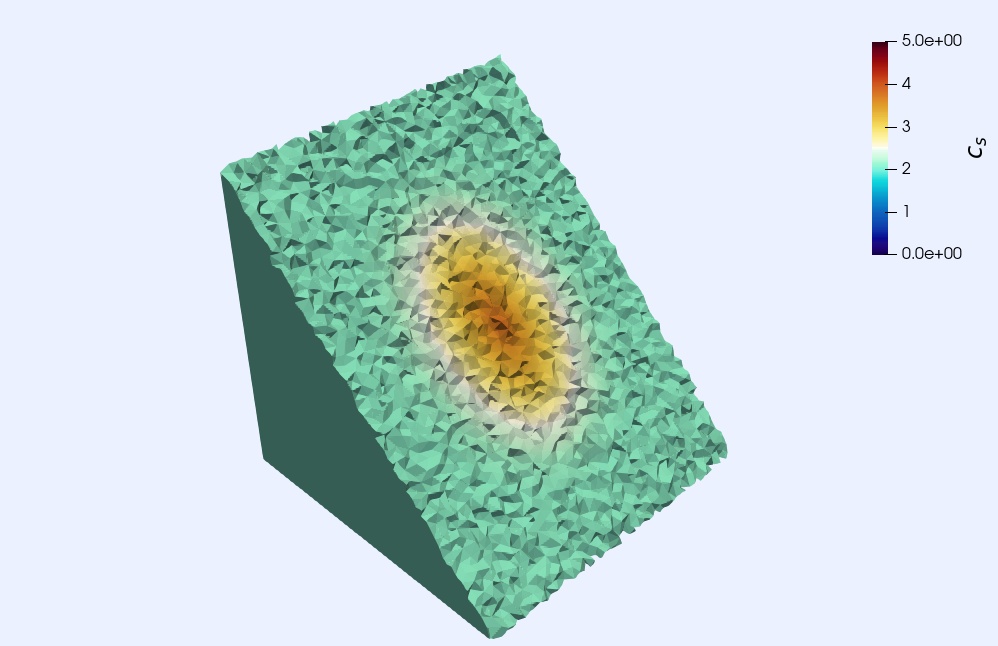}
  \includegraphics[trim = 70mm 0mm 90mm 18mm,  clip, width=0.32\linewidth]{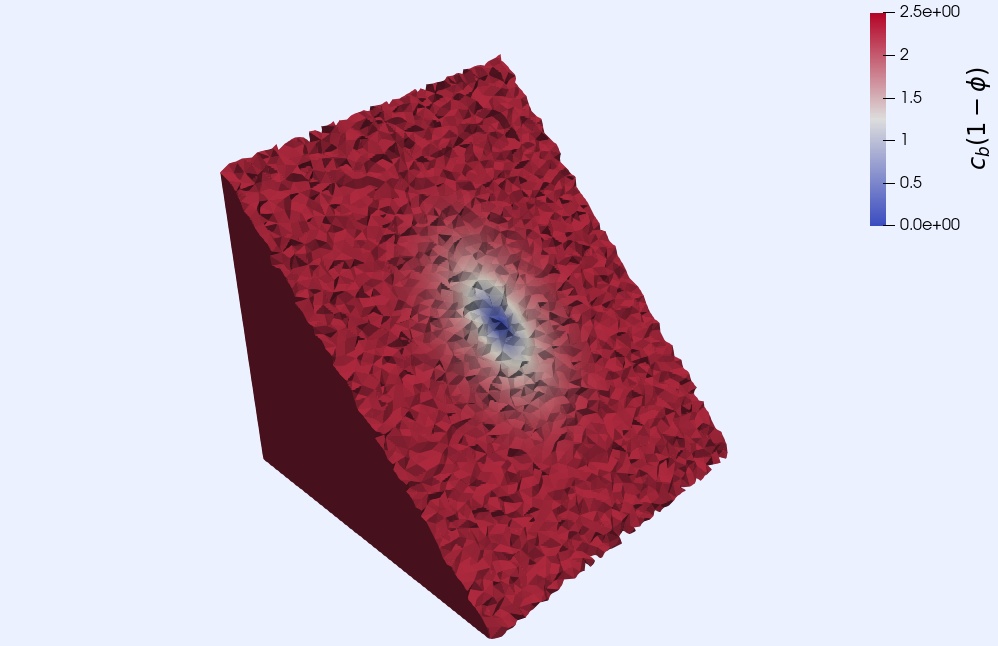}
    \includegraphics[trim = 70mm 0mm 90mm 18mm,  clip, width=0.32\linewidth]{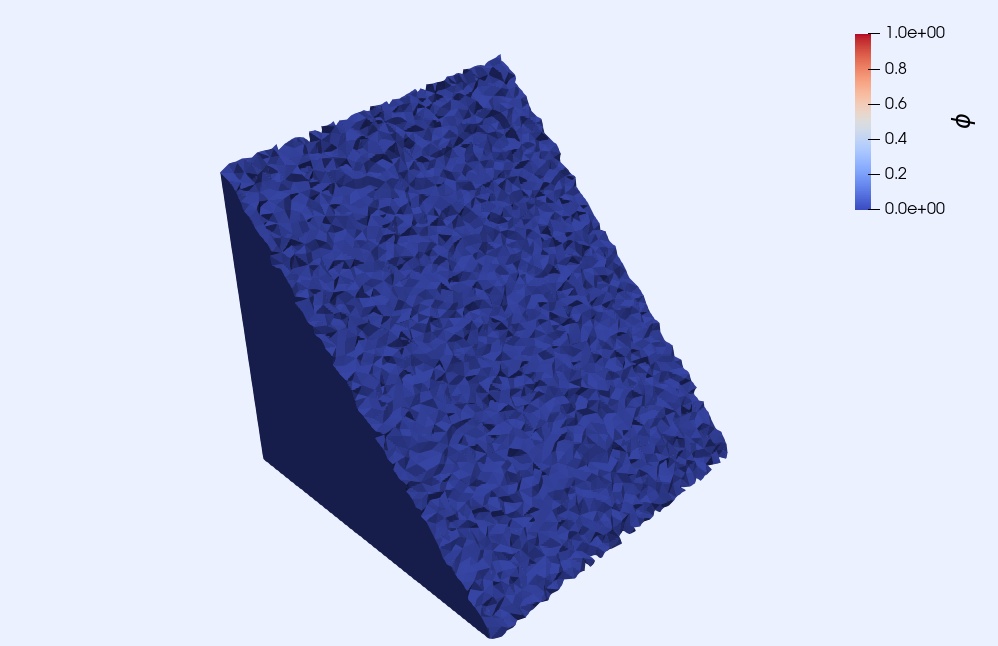}
  \includegraphics[trim = 70mm 0mm 90mm 18mm,  clip, width=0.32\linewidth]{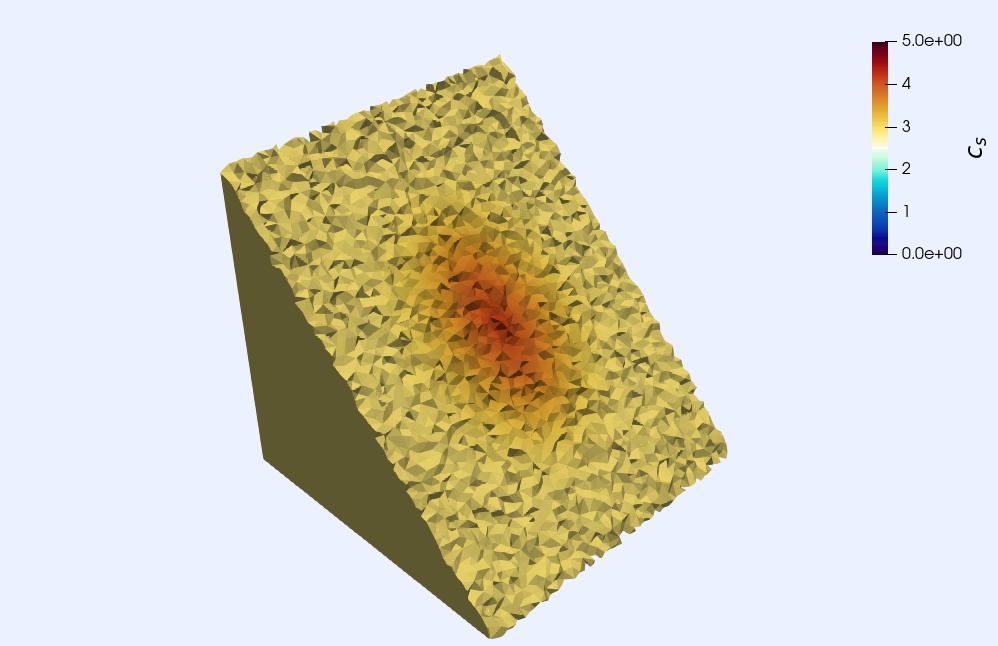}
  \includegraphics[trim = 70mm 0mm 90mm 18mm,  clip, width=0.32\linewidth]{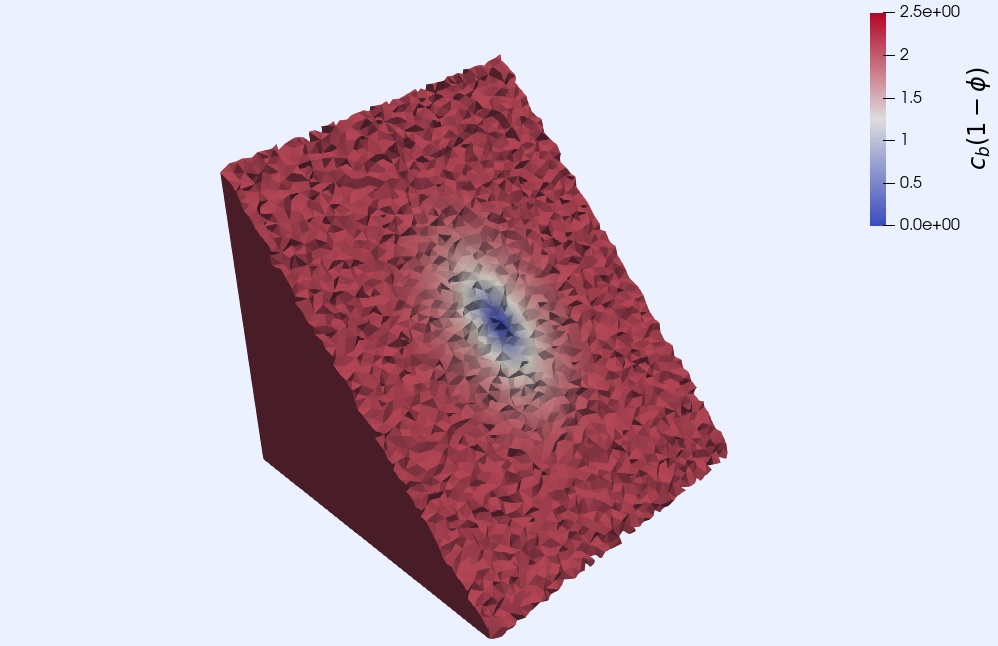}
\caption{$\mu_b=\mu_s=1$. Note for the results in this subfigure both soluble and bound MMPs degrade the matrix.}
  \label{fig:3d_both}
    \end{subfigure}
\hfill
  \begin{subfigure}[t]{0.49\textwidth}
    \includegraphics[angle=270,trim = 300mm 130mm 0mm 0mm,  clip, width=0.32\linewidth]{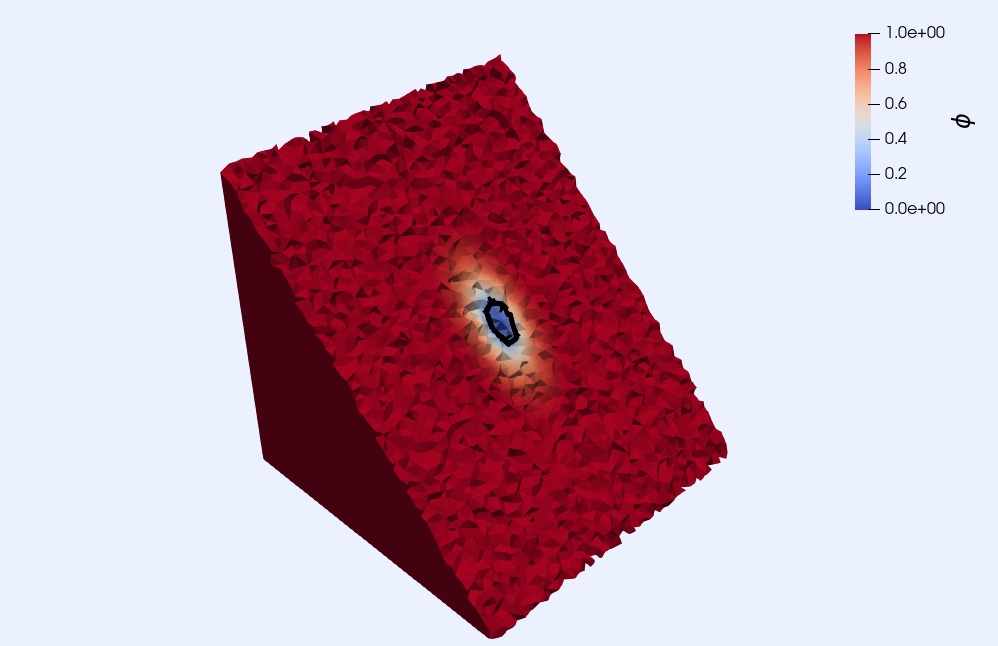}
  \includegraphics[angle=270,trim = 300mm 130mm 0mm 0mm,  clip, width=0.32\linewidth]{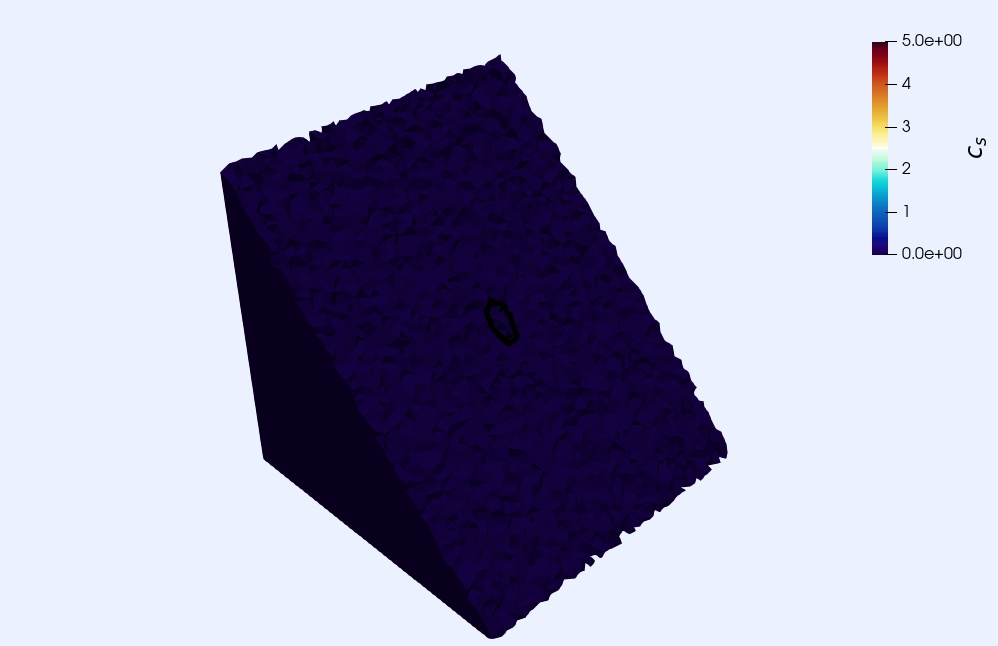}
  \includegraphics[angle=270,trim = 300mm 130mm 0mm 0mm,  clip, width=0.32\linewidth]{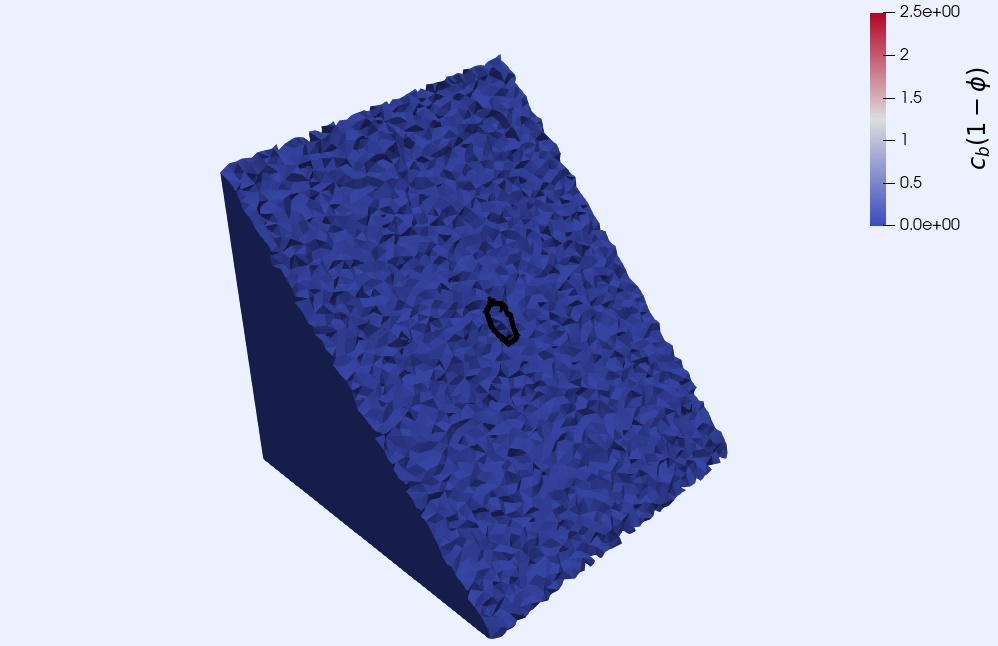}
  \includegraphics[trim = 70mm 0mm 90mm 18mm,  clip, width=0.32\linewidth]{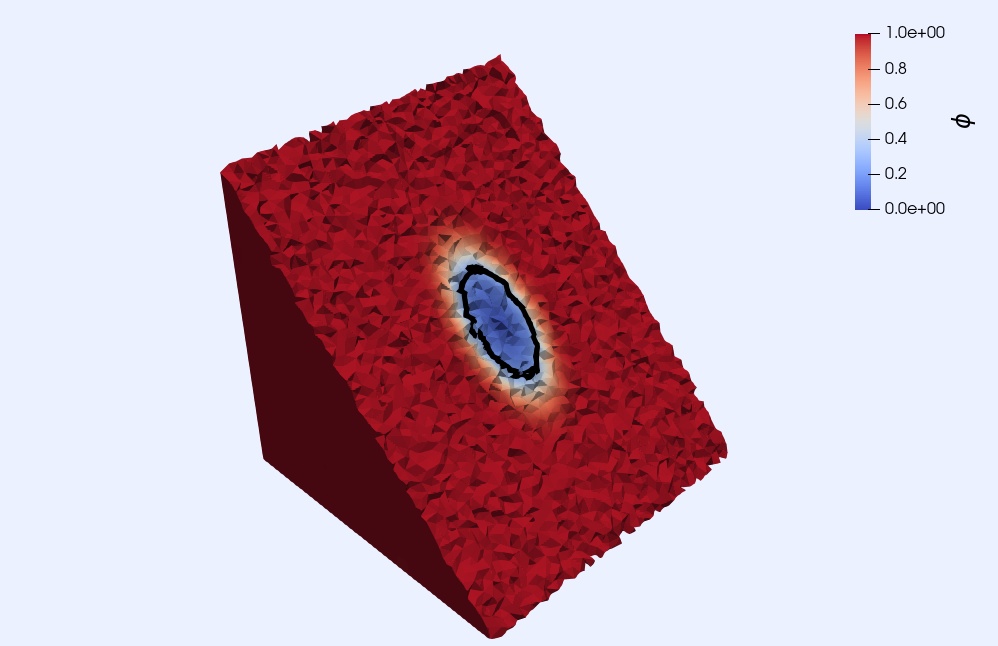}
  \includegraphics[trim = 70mm 0mm 90mm 18mm,  clip, width=0.32\linewidth]{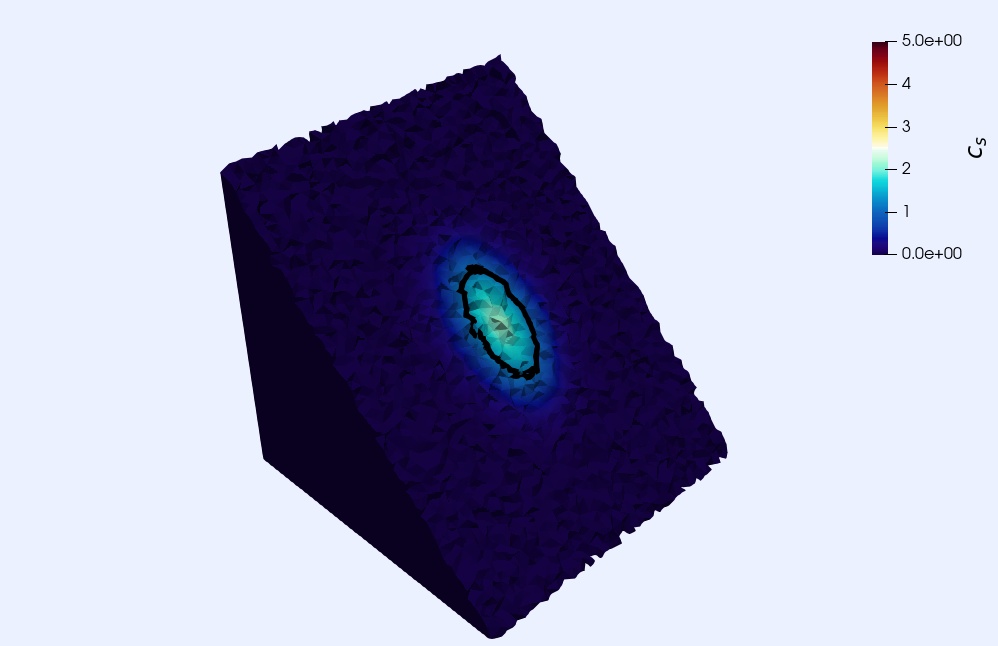}
  \includegraphics[trim = 70mm 0mm 90mm 18mm,  clip, width=0.32\linewidth]{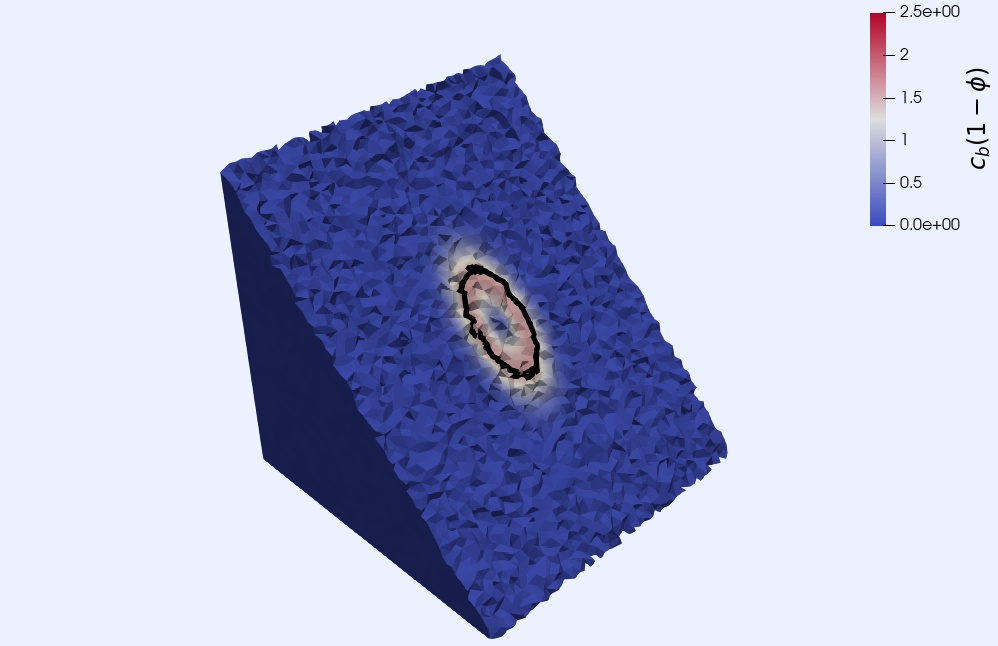}
    \includegraphics[trim = 70mm 0mm 90mm 18mm,  clip, width=0.32\linewidth]{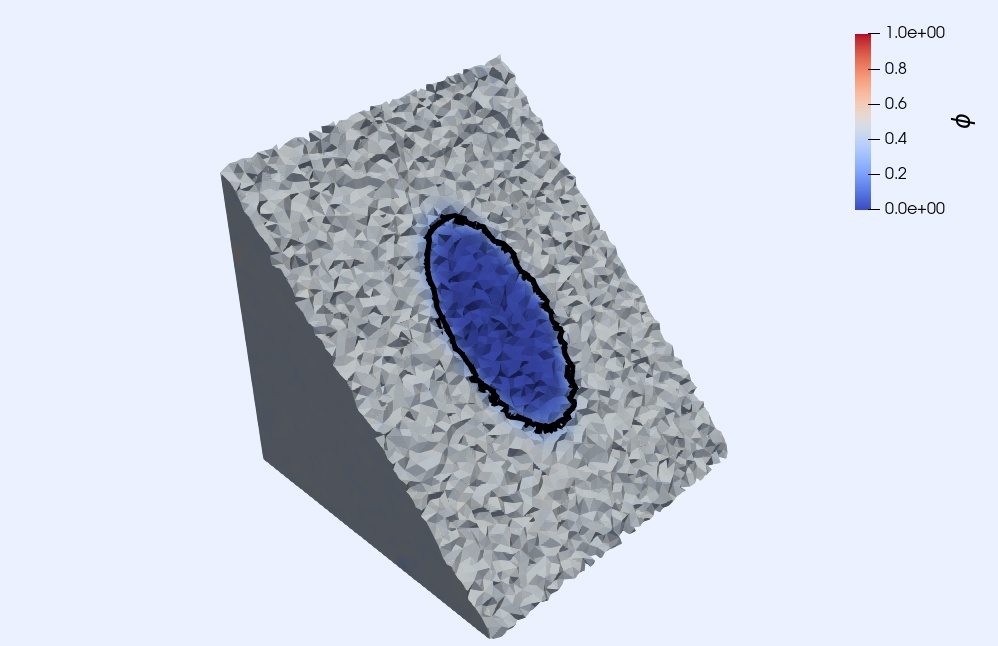}
  \includegraphics[trim = 70mm 0mm 90mm 18mm,  clip, width=0.32\linewidth]{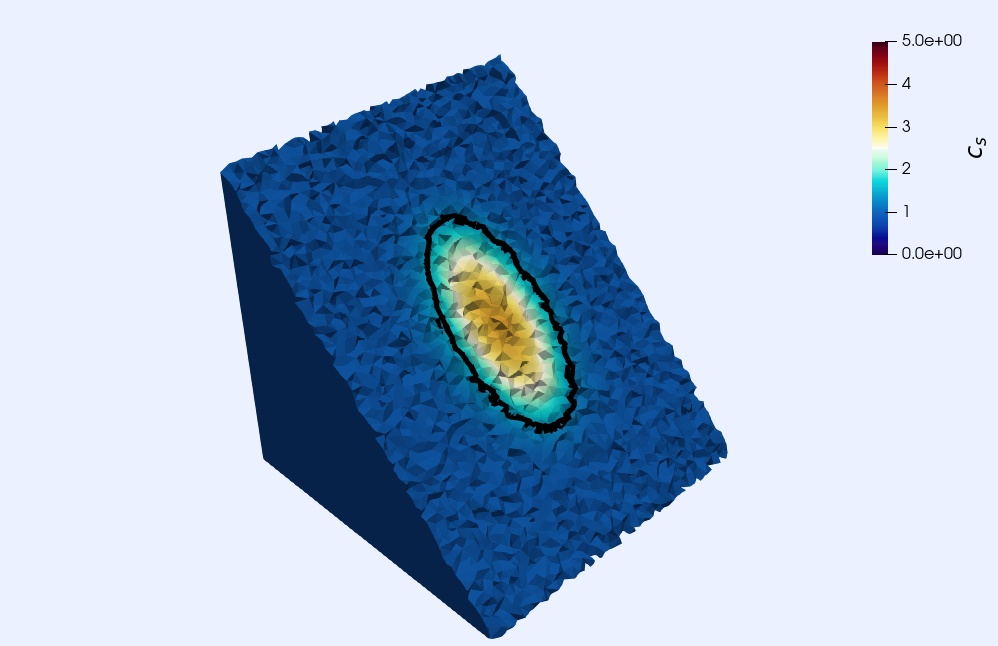}
  \includegraphics[trim = 70mm 0mm 90mm 18mm,  clip, width=0.32\linewidth]{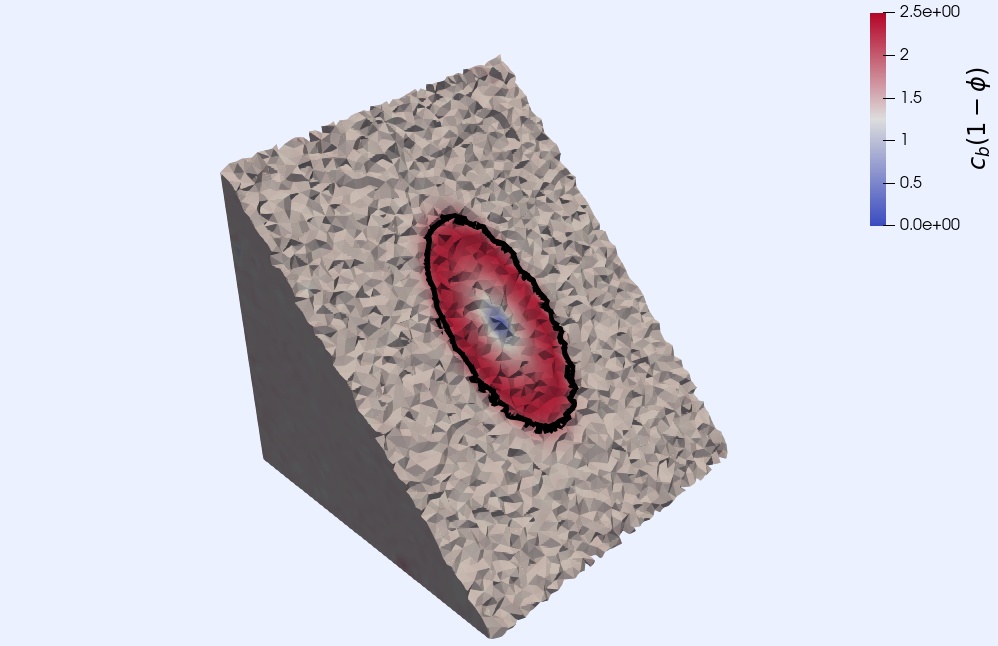}
    \includegraphics[trim = 70mm 0mm 90mm 18mm,  clip, width=0.32\linewidth]{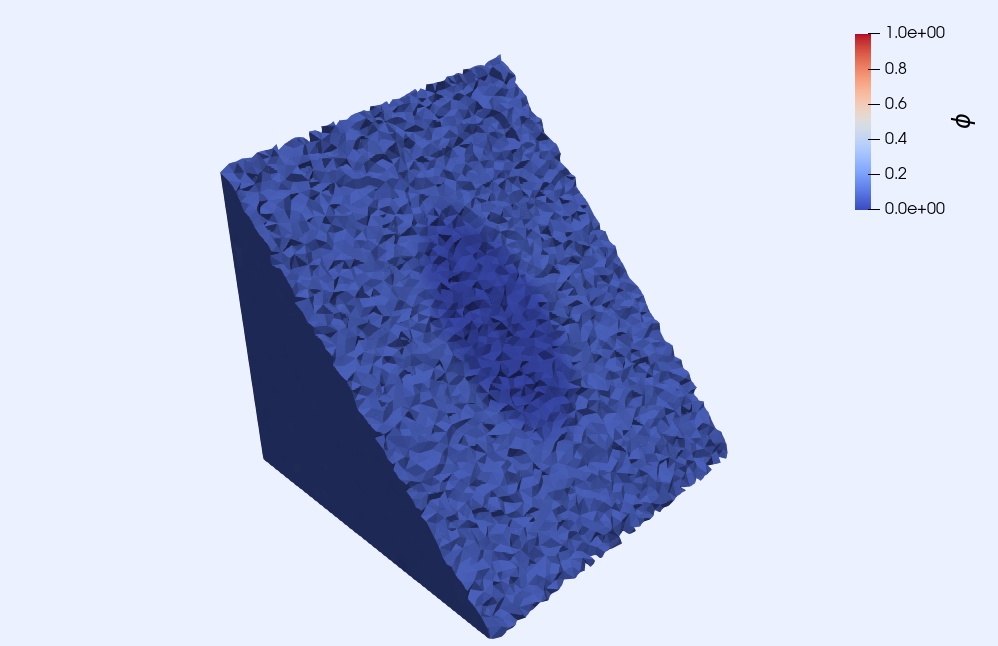}
  \includegraphics[trim = 70mm 0mm 90mm 18mm,  clip, width=0.32\linewidth]{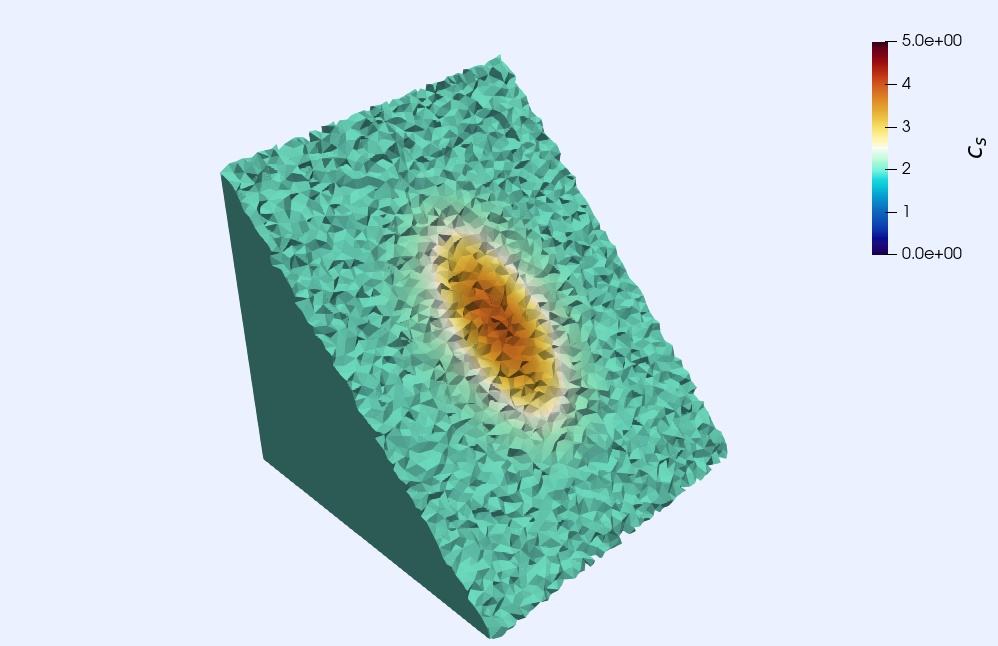}
  \includegraphics[trim = 70mm 0mm 90mm 18mm,  clip, width=0.32\linewidth]{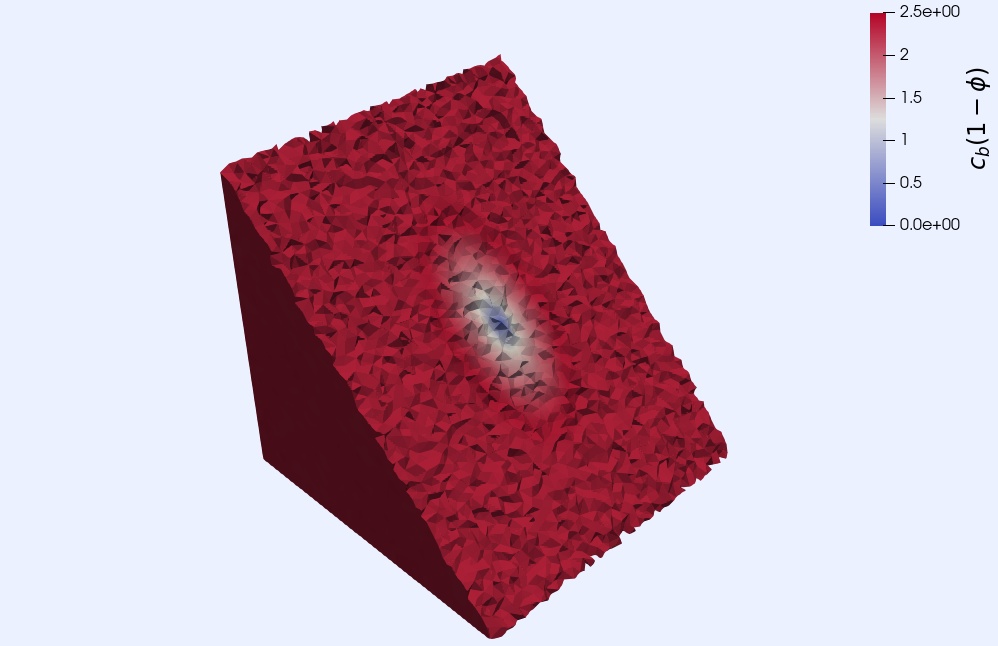}
\caption{$\mu_b=1,\mu_s=0$. Note for the results in this subfigure the soluble MMPs do not influence $\phi$.}
  \label{fig:3d_mub}
    \end{subfigure}
\hfill
  \begin{subfigure}[t]{0.49\textwidth}
    \includegraphics[angle=270,trim = 300mm 130mm 0mm 0mm,  clip, width=0.32\linewidth]{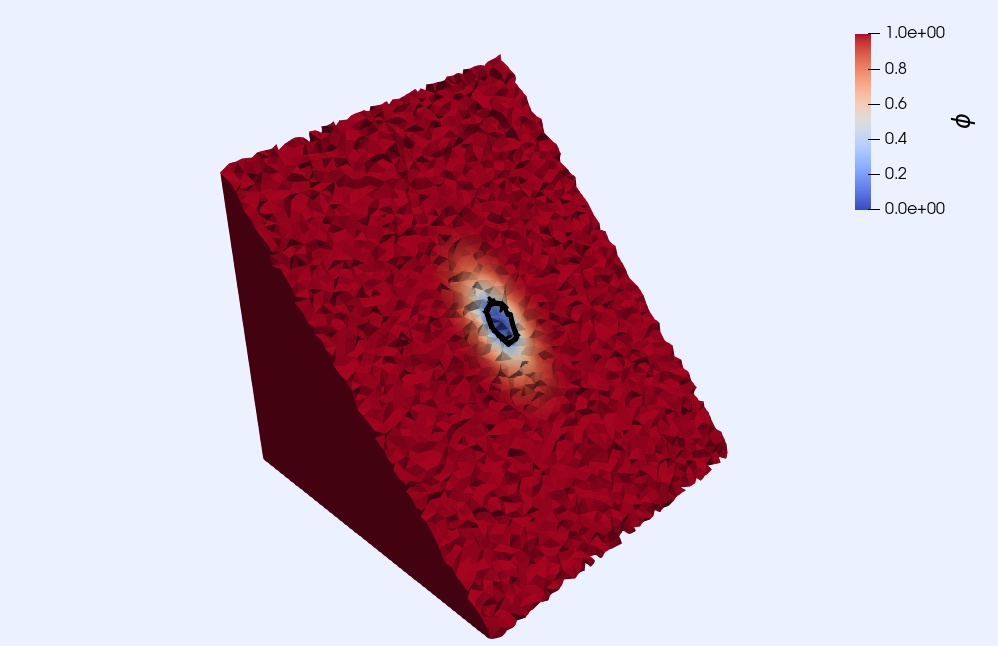}
  \includegraphics[angle=270,trim = 300mm 130mm 0mm 0mm,  clip, width=0.32\linewidth]{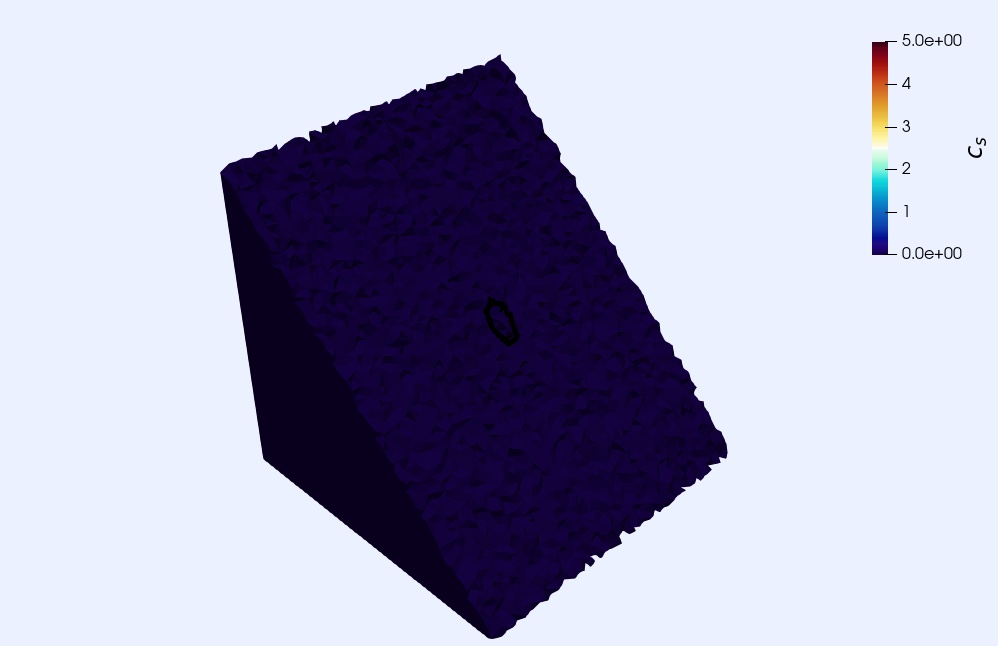}
  \includegraphics[angle=270,trim = 300mm 130mm 0mm 0mm,  clip, width=0.32\linewidth]{./3d_Images/mub/cb1minphi0000.jpeg}
  \includegraphics[trim = 70mm 0mm 90mm 18mm,  clip, width=0.32\linewidth]{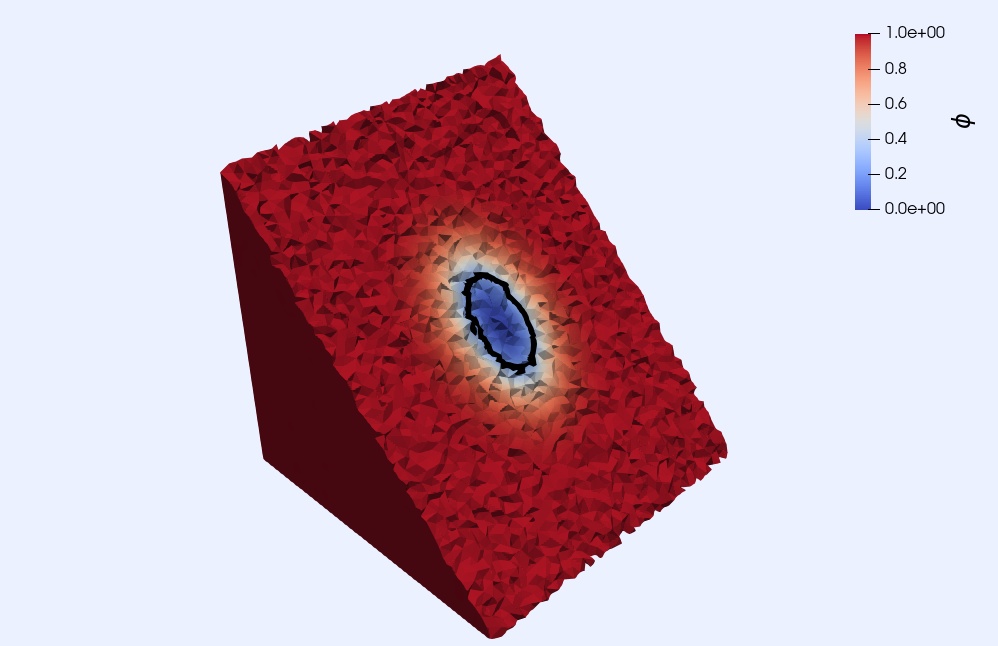}
  \includegraphics[trim = 70mm 0mm 90mm 18mm,  clip, width=0.32\linewidth]{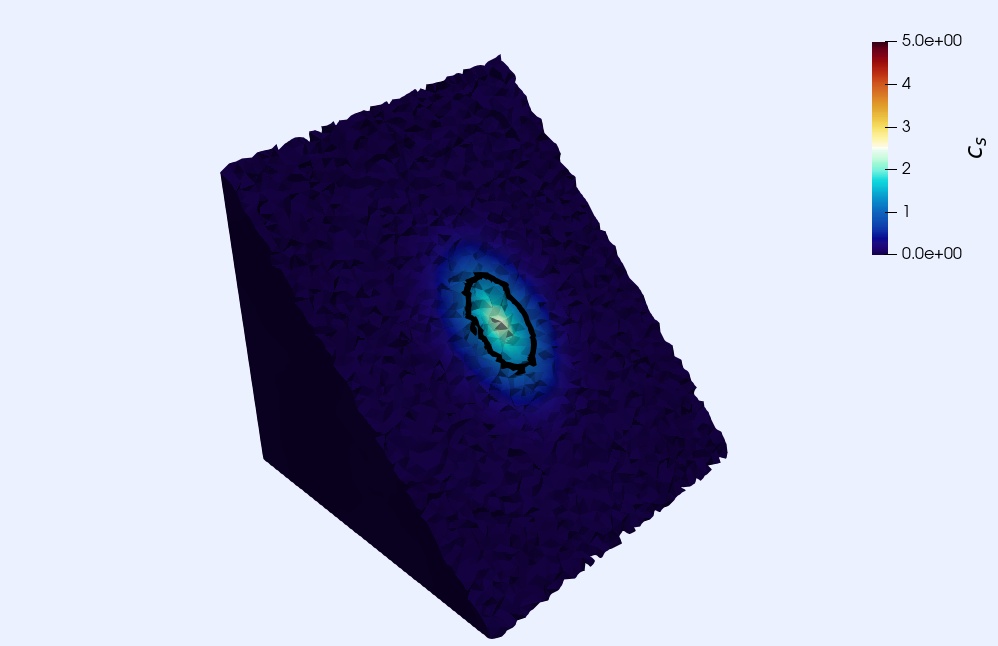}
  \includegraphics[trim = 160mm 0mm 185mm 18mm,  clip, width=0.32\linewidth]{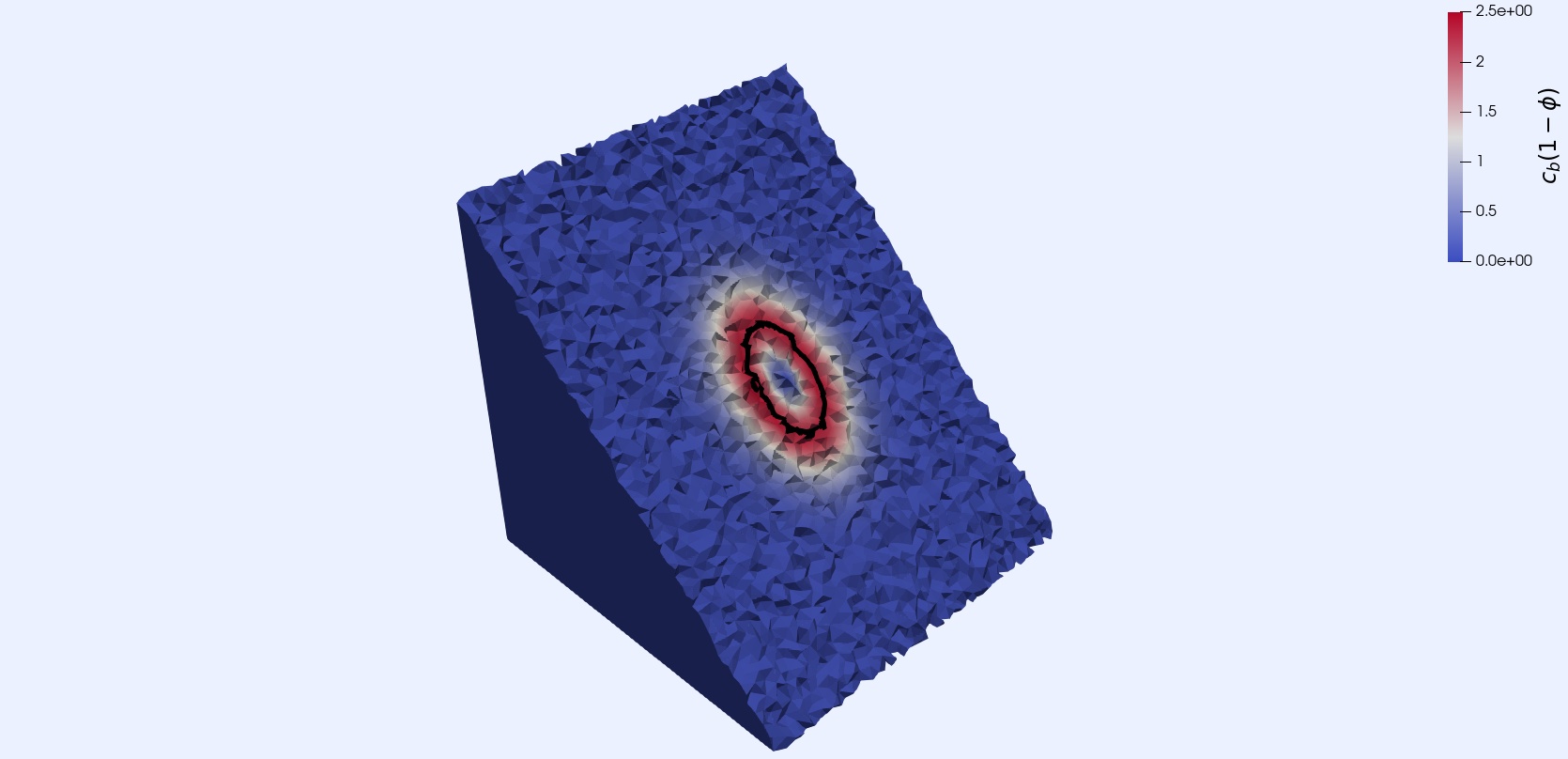}
    \includegraphics[trim = 70mm 0mm 90mm 18mm,  clip, width=0.32\linewidth]{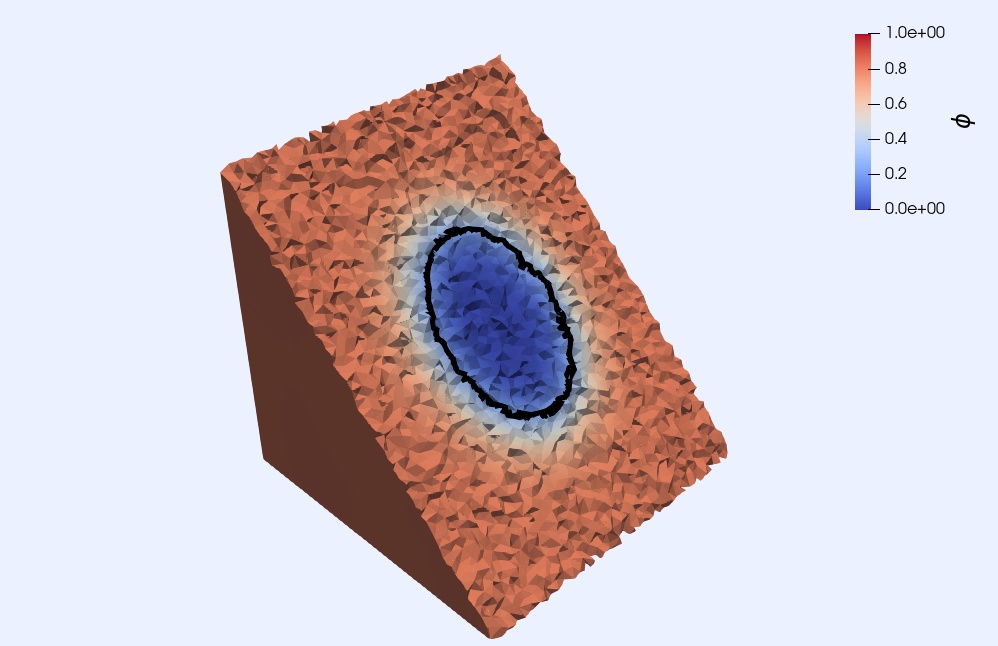}
  \includegraphics[trim = 70mm 0mm 90mm 18mm,  clip, width=0.32\linewidth]{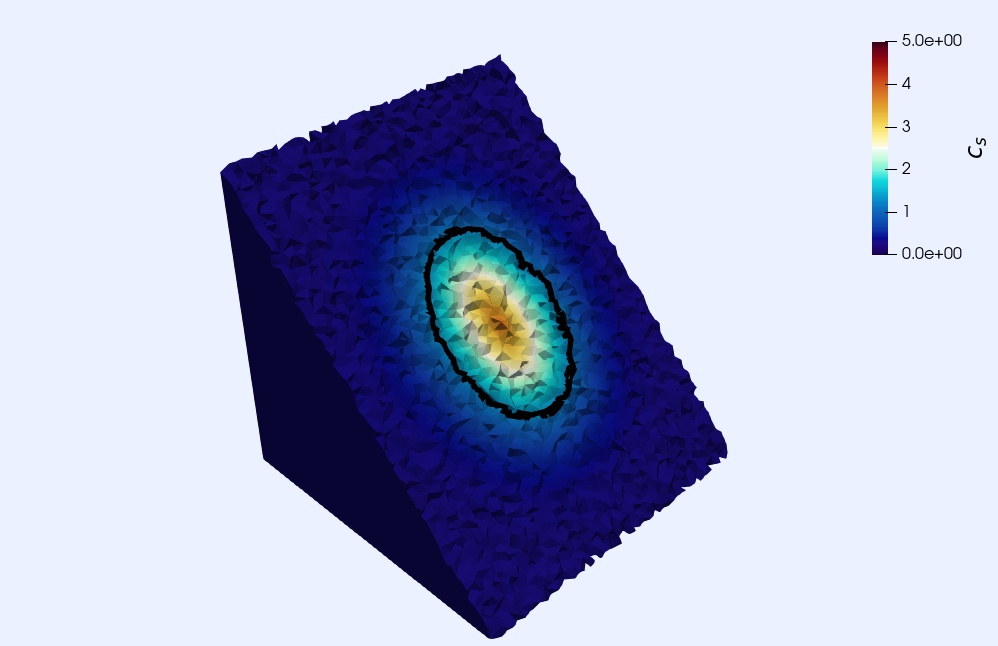}
  \includegraphics[trim = 160mm 0mm 185mm 18mm,  clip, width=0.32\linewidth]{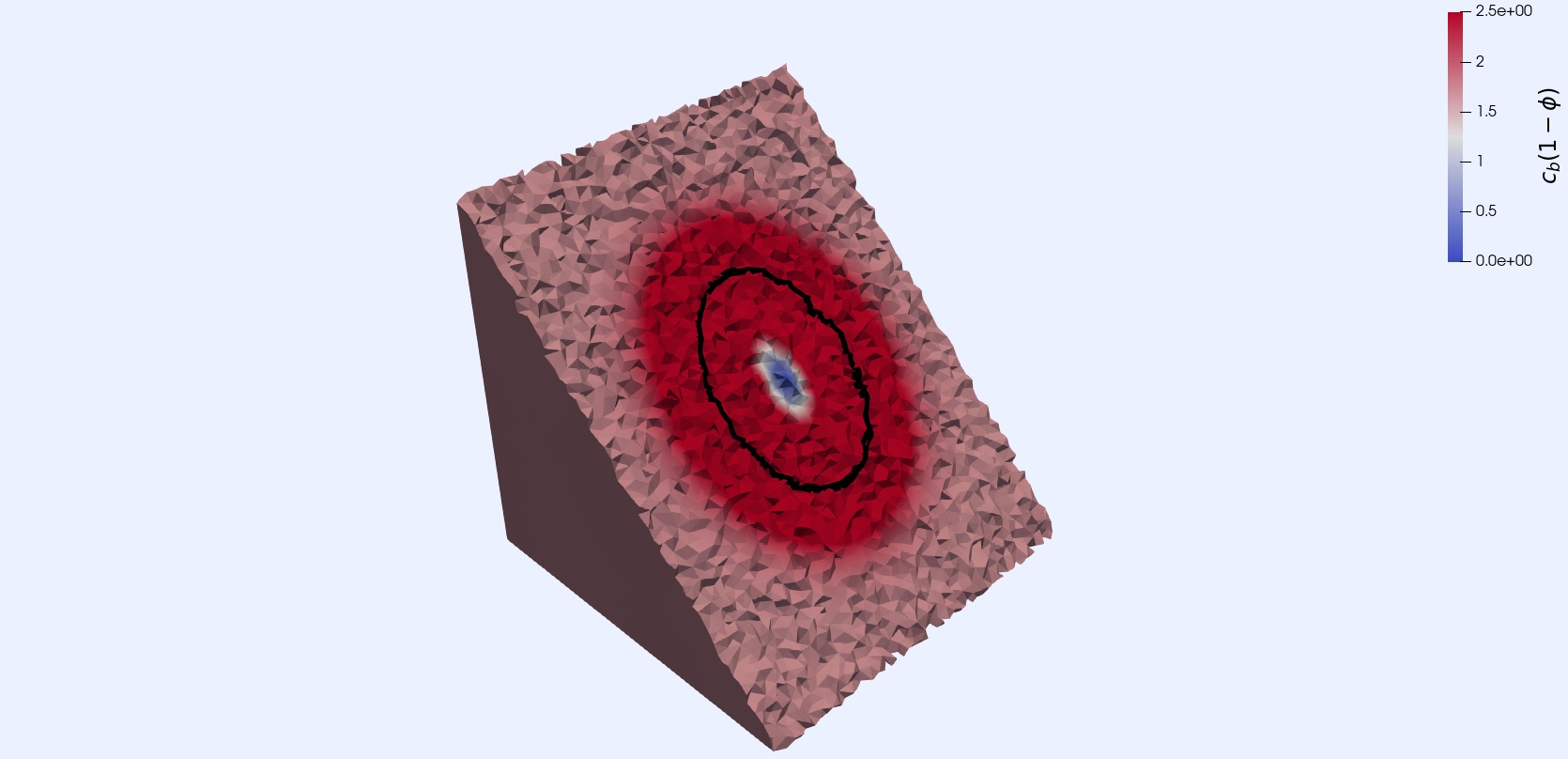}
    \includegraphics[trim = 70mm 0mm 90mm 18mm,  clip, width=0.32\linewidth]{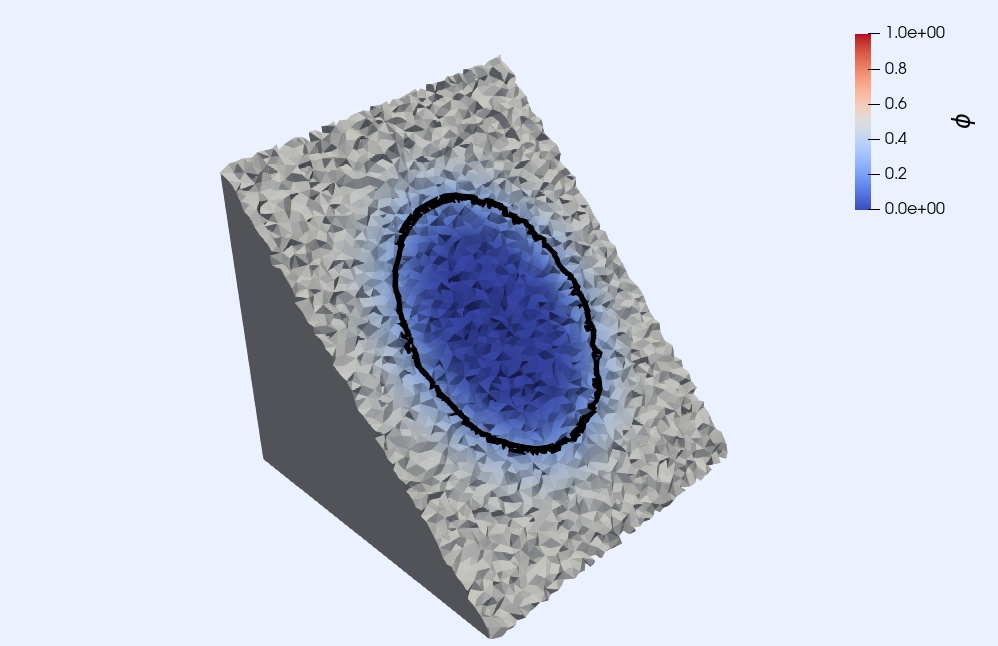}
  \includegraphics[trim = 70mm 0mm 90mm 18mm,  clip, width=0.32\linewidth]{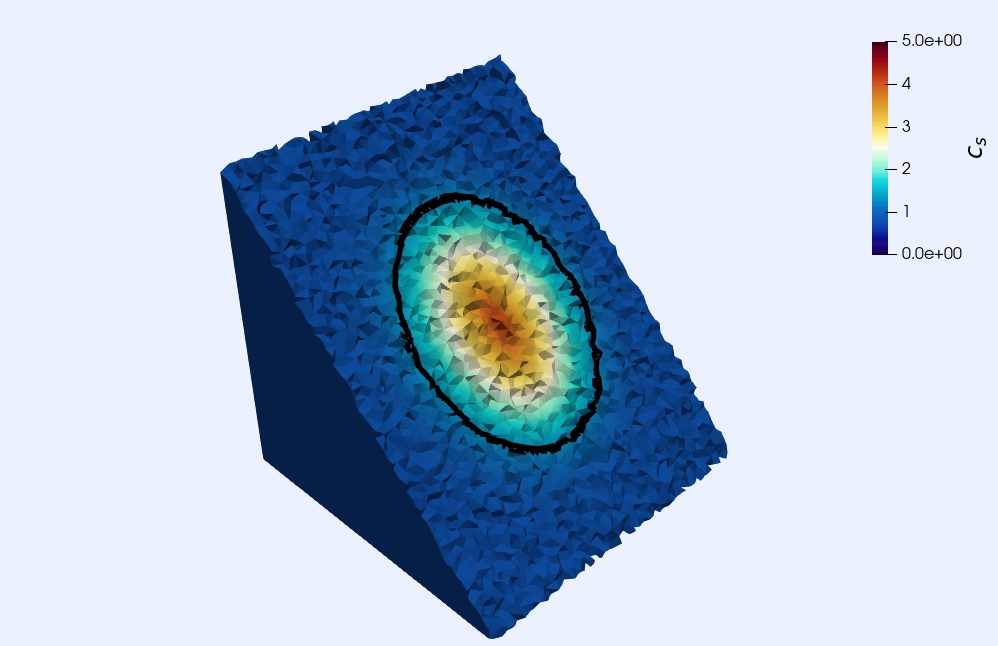}
  \includegraphics[trim = 160mm 0mm 185mm 18mm,  clip, width=0.32\linewidth]{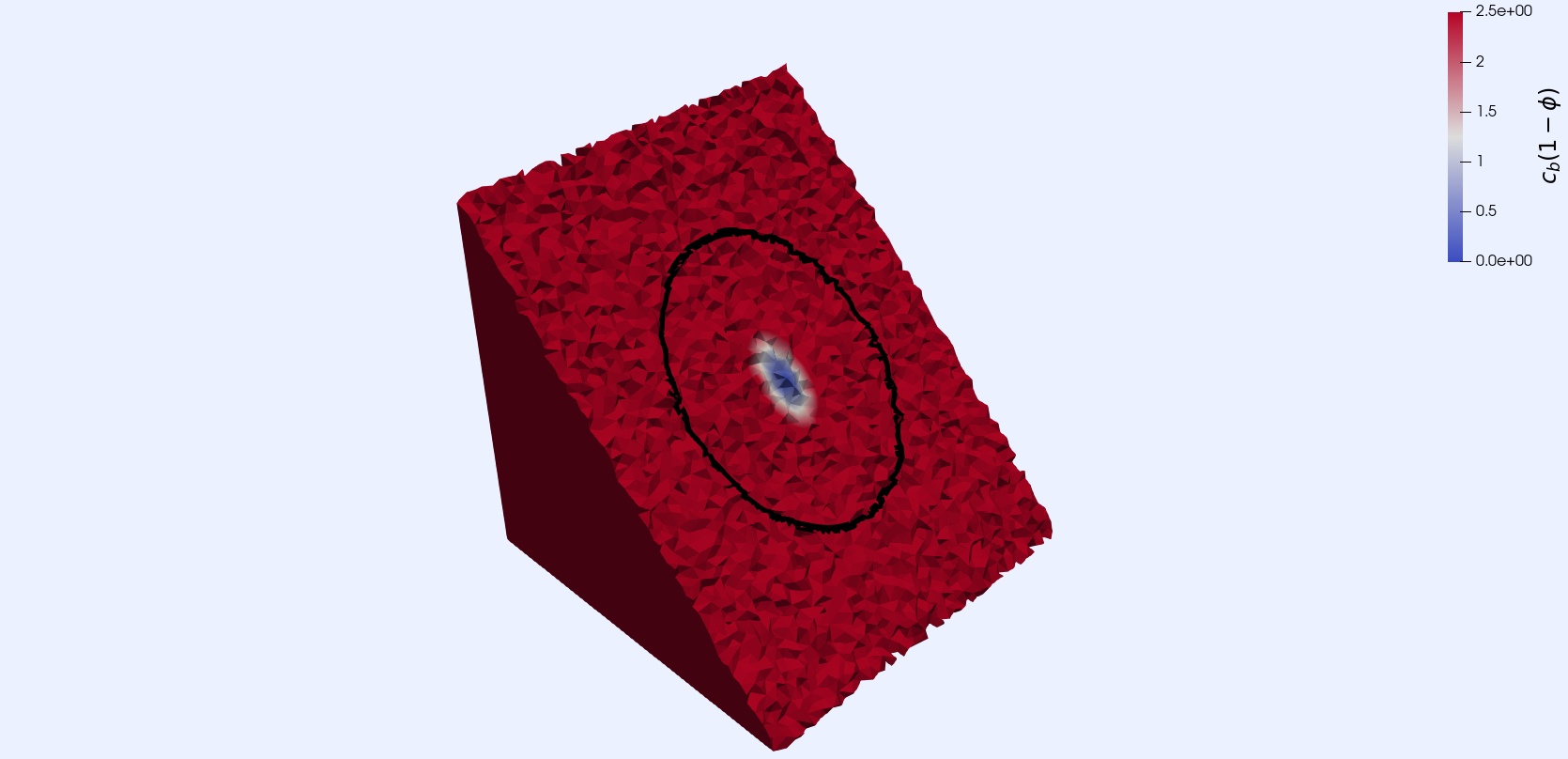}
\caption{$\mu_b=0,\mu_s=1$. Note for the results in this subfigure the bound MMPs do not influence $\phi$.}
  \label{fig:3d_mus}
    \end{subfigure}
    \caption{Simulations of the invasion model \eqref{eqn:model_macro}, \eqref{IC_BC} in $3{\rm d}$. In all  plots we have half of the domain transparent to aid visualisation.  In each subfigure, each row corresponds to times $t=2, 4$ and $5$. The black contour indicates the level set $\phi=0.25$ corresponding to a volume fraction of $75\%$ cancer cells.  The results appear analogous to the $2{\rm d}$ case despite the relatively higher diffusivity of soluble MMPs for $\phi>0$. The effect of membrane bound MMPs  increasing the speed of invasion is apparent whilst the degradation of the ECM by the soluble MMPs appears to generate more radially symmetric invasive profiles. Parameter values as in Table~\ref{table_1}.}
    \label{fig:3d}
  \end{figure}
  Figure \ref{fig:3d} shows results of  $3{\rm d}$ simulations. The results are analogous to the $2{\rm d}$ simulations, see  Figure~\ref{fig:2d}, and hence we infer that in this parameter regime the model  exhibits similar behaviour in $2{\rm d}$ and $3{\rm d}$ despite the qualitatively different effective diffusivities. In the following sections of the manuscript, we therefore report only on $2{\rm d}$ simulations and expect that analogous results would be obtained for the corresponding $3{\rm d}$ simulations.

 \subsection{Simulations of the model with varying matrix suitability~\eqref{het_ecm},~\eqref{IC_BC}}\label{subsec:het_ecm_numerics} 
We now report on simulations of the model~\eqref{het_ecm},~\eqref{IC_BC} with varying matrix suitability for invasion. For all  simulations in this section we take $\Omega=[-1,1]^2$, the time step $\tau=10^{-2}$ and use a mesh with $12909$ DOFs. We first consider initial conditions as in the $2{\rm d}$ simulations of section \ref{subsec:macro_model_numerics} along with initial conditions for $s$ of the form
\begin{equation}\label{ic_s_unsuitable}
s_0(x)=1-0.1\left(\cos(4\pi x_1)\cos(4\pi x_2)\right)^2,
\end{equation}
which corresponds to a matrix that is largely unsuitable for invasion at the initial time.  
\begin{figure}[htbp]
  \begin{subfigure}[t]{0.49\textwidth}
  \includegraphics[trim = 60mm 160mm 90mm 30mm,  clip, width=\linewidth]{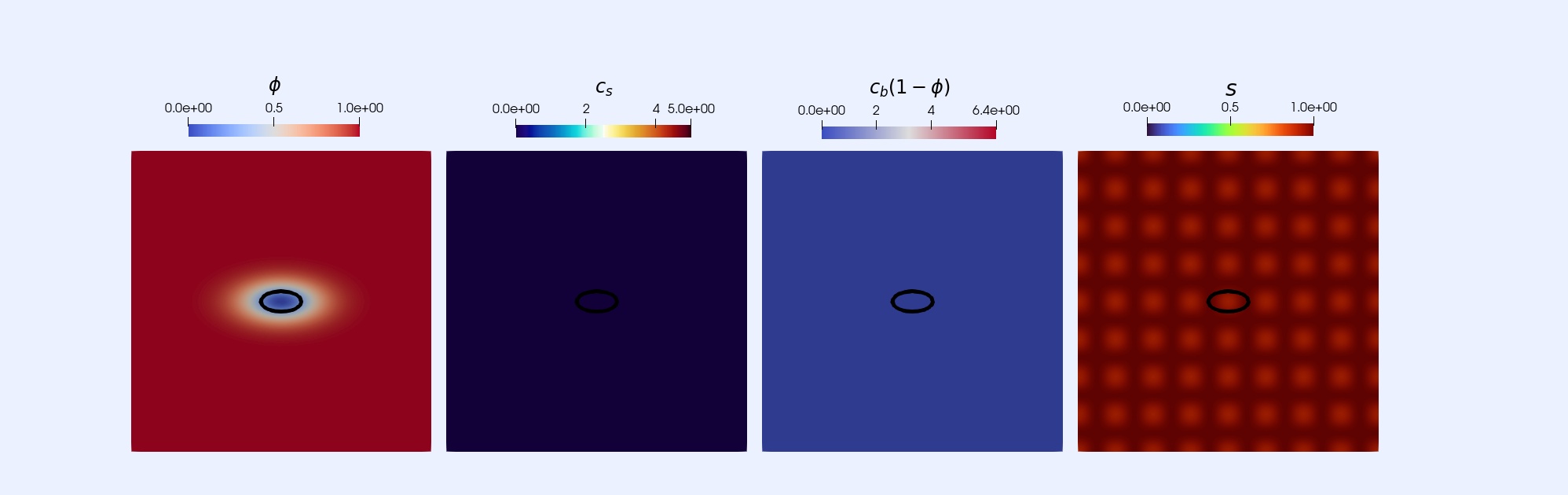}
  \includegraphics[trim = 60mm 0mm 90mm 65mm,  clip, width=\linewidth]{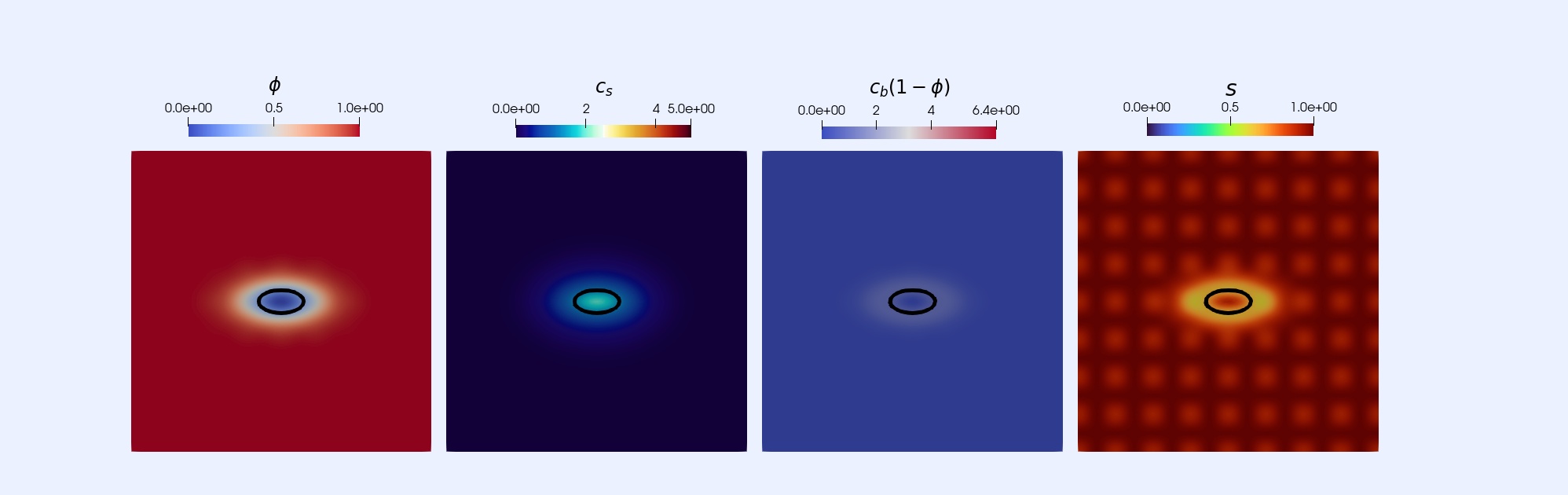}
  \includegraphics[trim = 60mm 0mm 90mm 65mm,  clip, width=\linewidth]{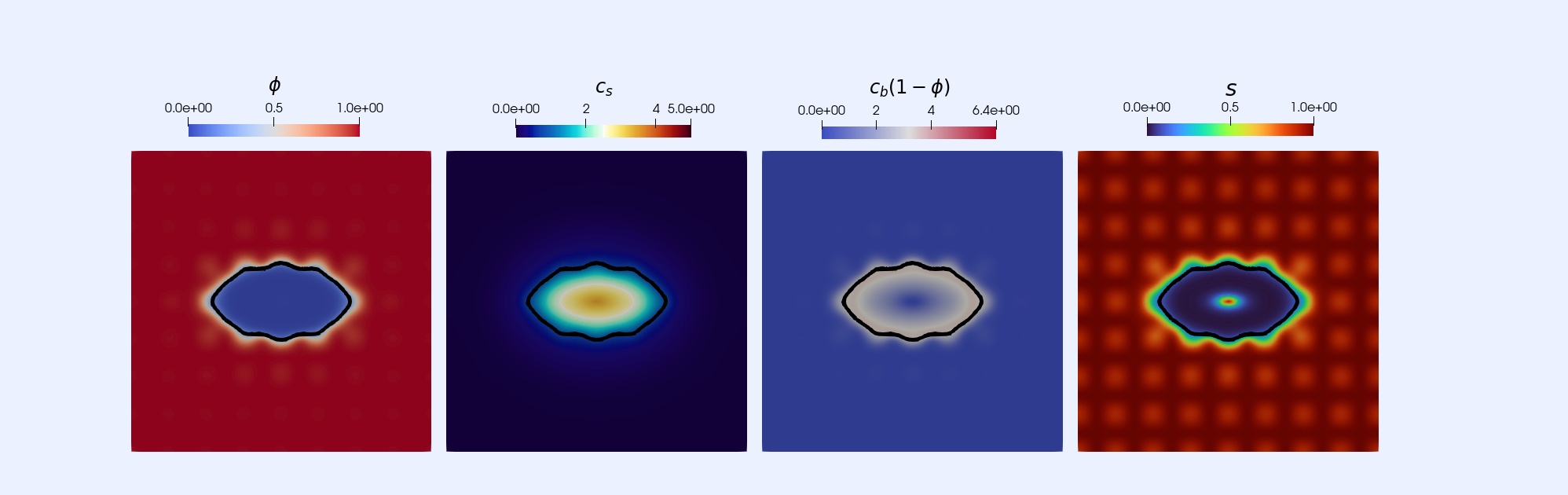}
  \includegraphics[trim = 60mm 0mm 90mm 65mm,  clip, width=\linewidth]{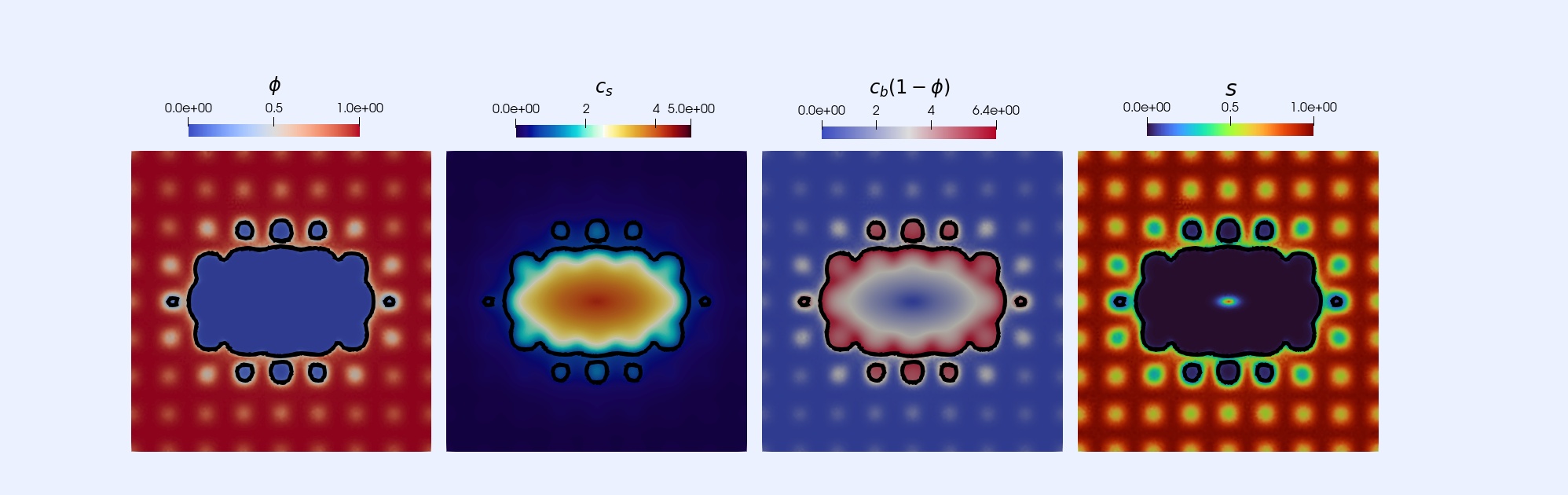}
\caption{$\mu_b=\mu_s=1$. Note for the results in this subfigure both soluble and bound MMPs degrade the matrix.}
  \label{fig:2d_suitability_both}
    \end{subfigure}
\hfill
  \begin{subfigure}[t]{0.49\textwidth}
  \includegraphics[trim = 60mm 160mm 90mm 30mm,  clip, width=\linewidth]{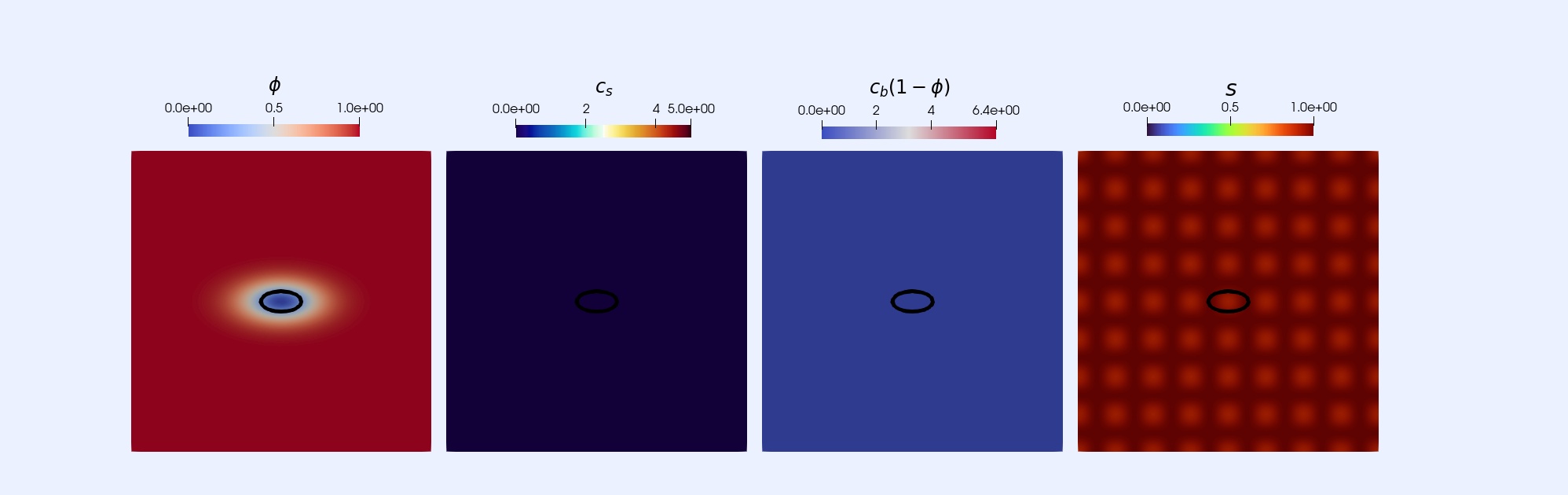}
  \includegraphics[trim = 60mm 0mm 90mm 65mm,  clip, width=\linewidth]{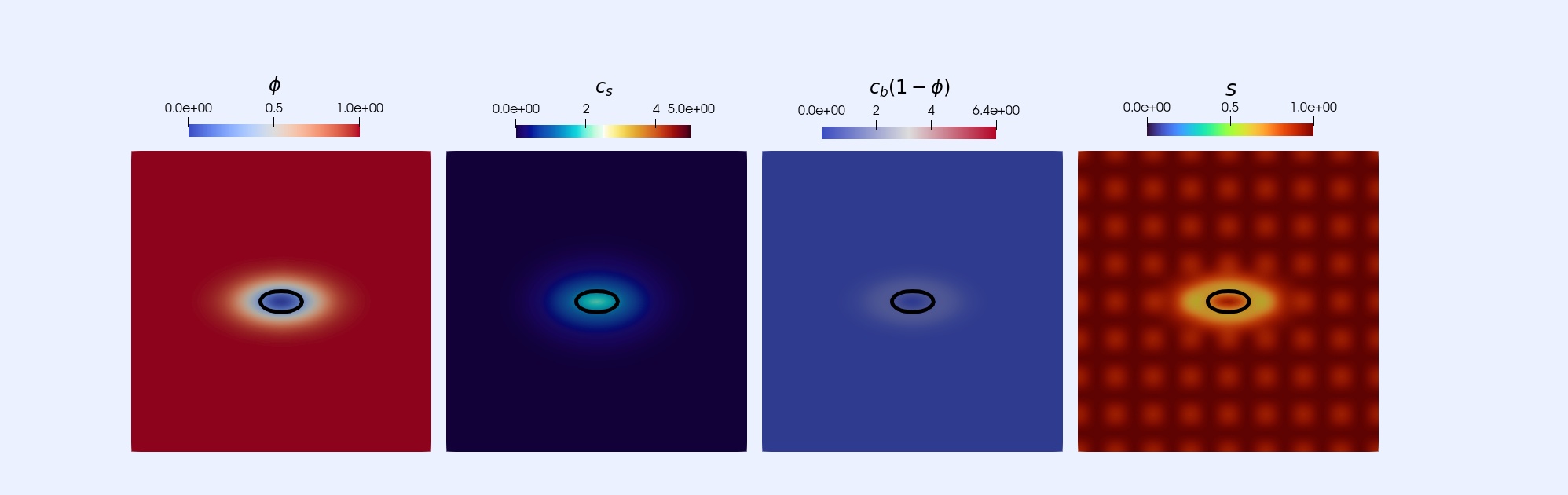}
  \includegraphics[trim = 60mm 0mm 90mm 65mm,  clip, width=\linewidth]{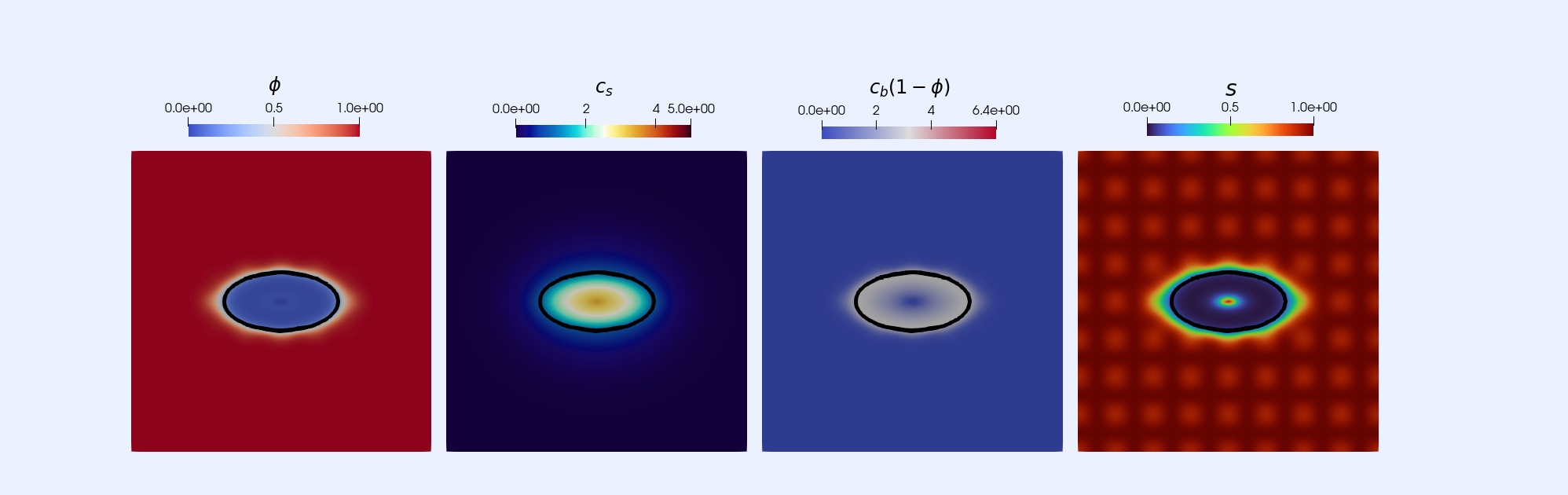}
  \includegraphics[trim = 60mm 0mm 90mm 65mm,  clip, width=\linewidth]{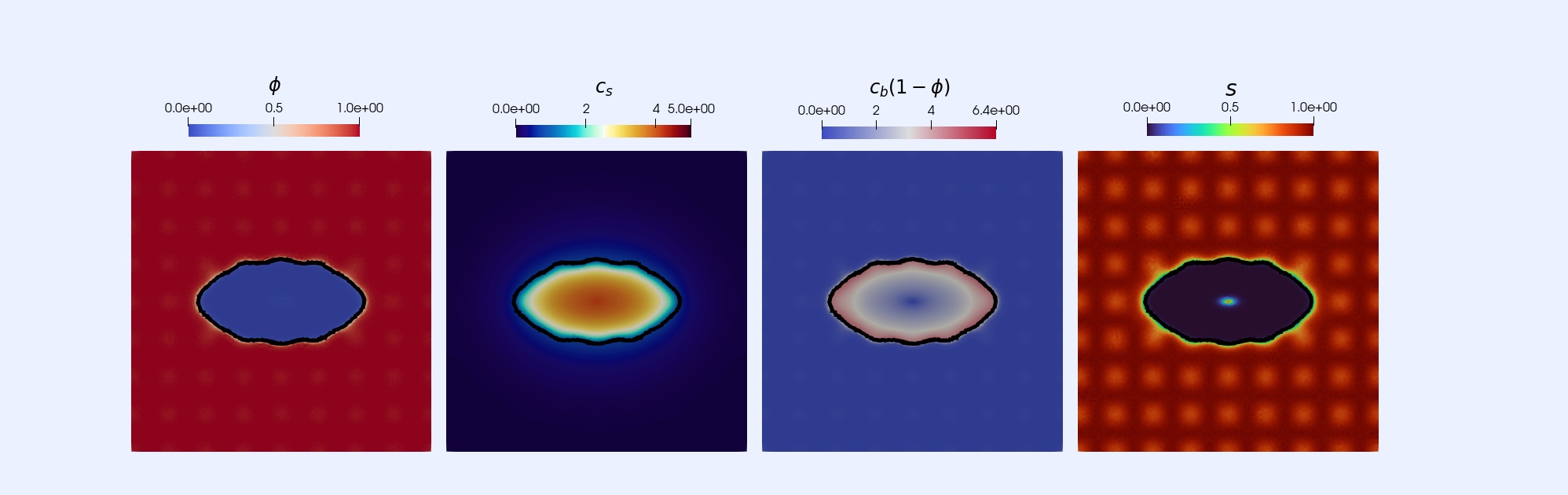}
\caption{$\mu_b=1,\mu_s=0$. Note for the results in this subfigure the soluble MMPs do not influence $\phi$.}
  \label{fig:2d_suitability_mub}
    \end{subfigure}
\hfill
  \begin{subfigure}[t]{0.49\textwidth}
  \includegraphics[trim = 60mm 160mm 90mm 30mm,  clip, width=\linewidth]{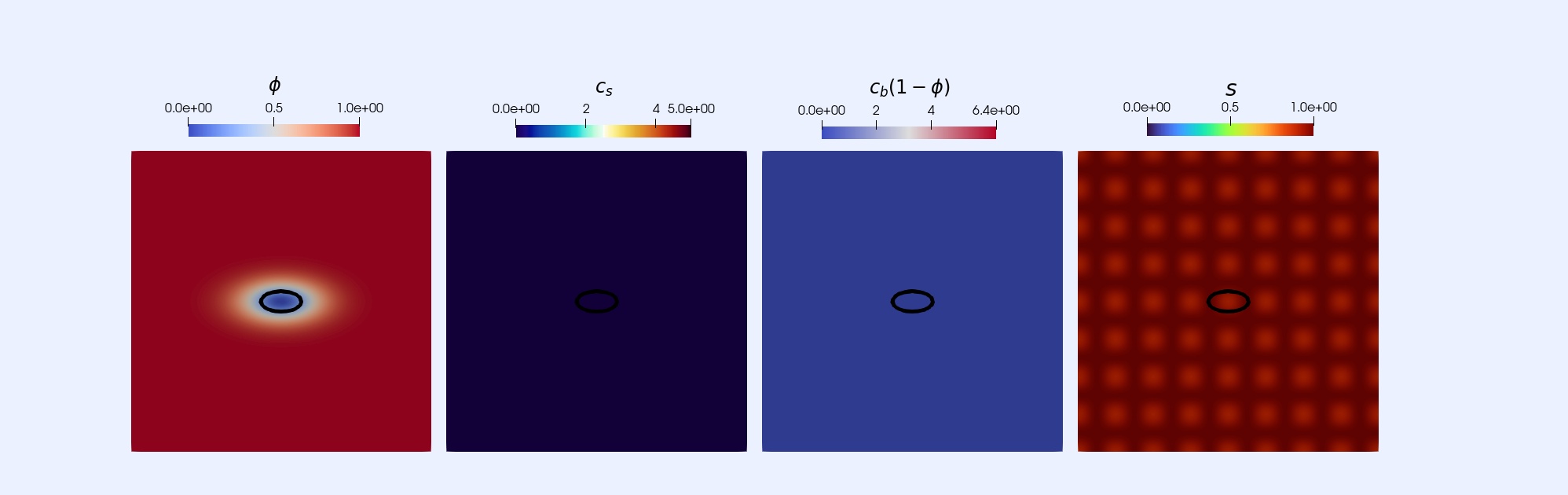}
  \includegraphics[trim = 60mm 0mm 90mm 65mm,  clip, width=\linewidth]{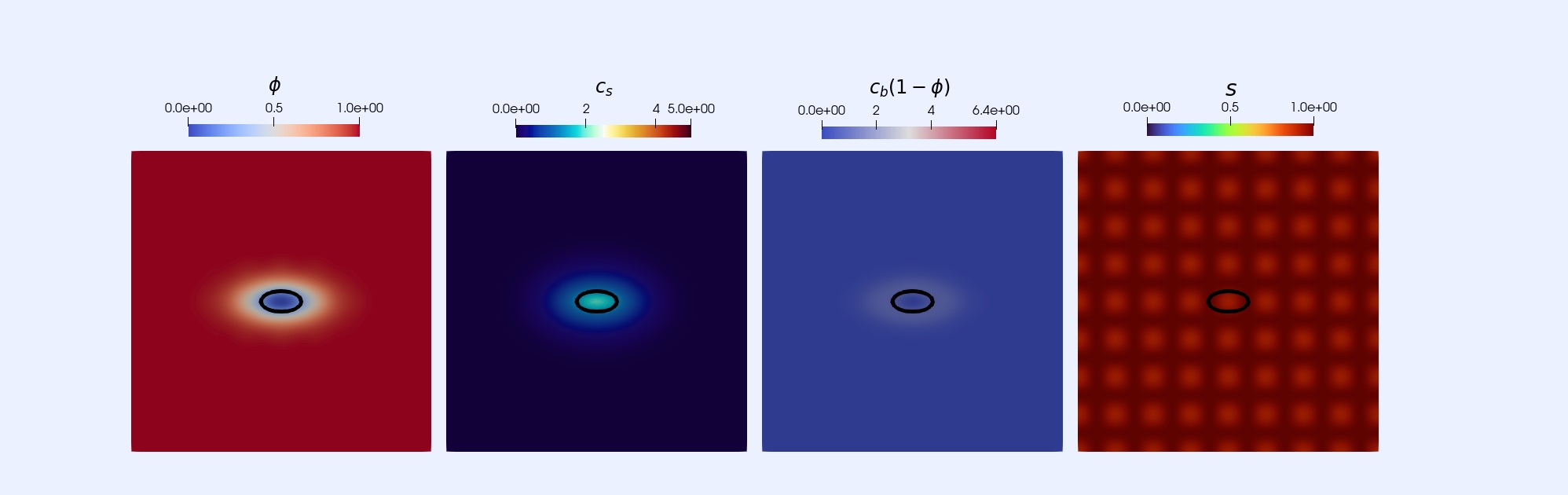}
  \includegraphics[trim = 60mm 0mm 90mm 65mm,  clip, width=\linewidth]{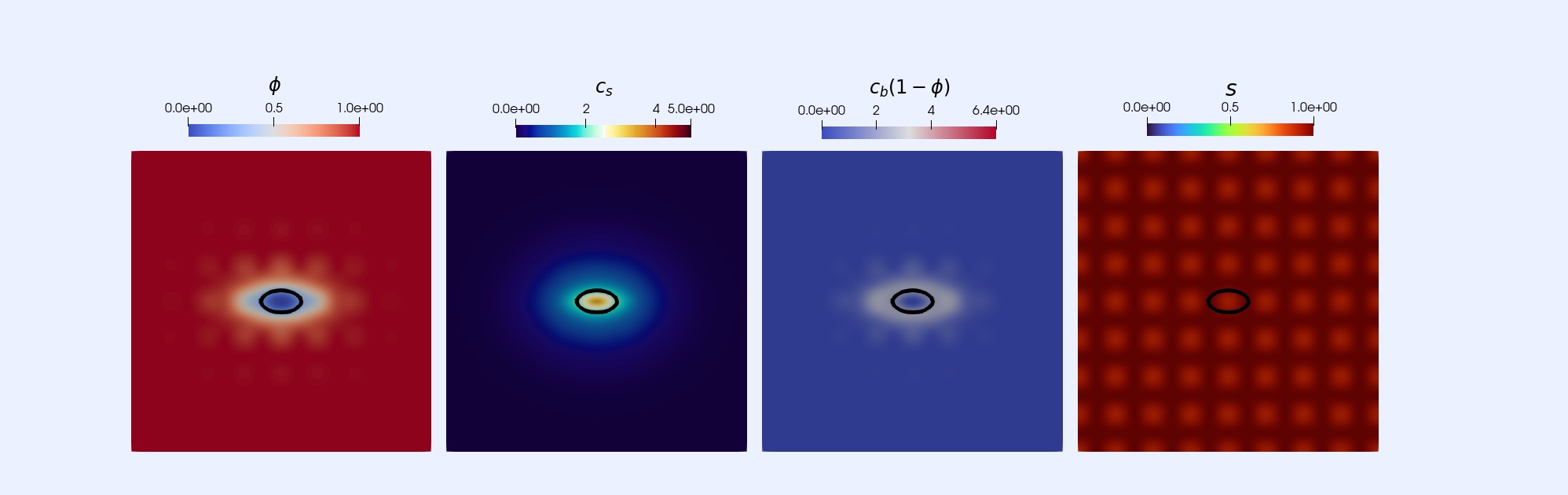}
  \includegraphics[trim = 60mm 0mm 90mm 65mm,  clip, width=\linewidth]{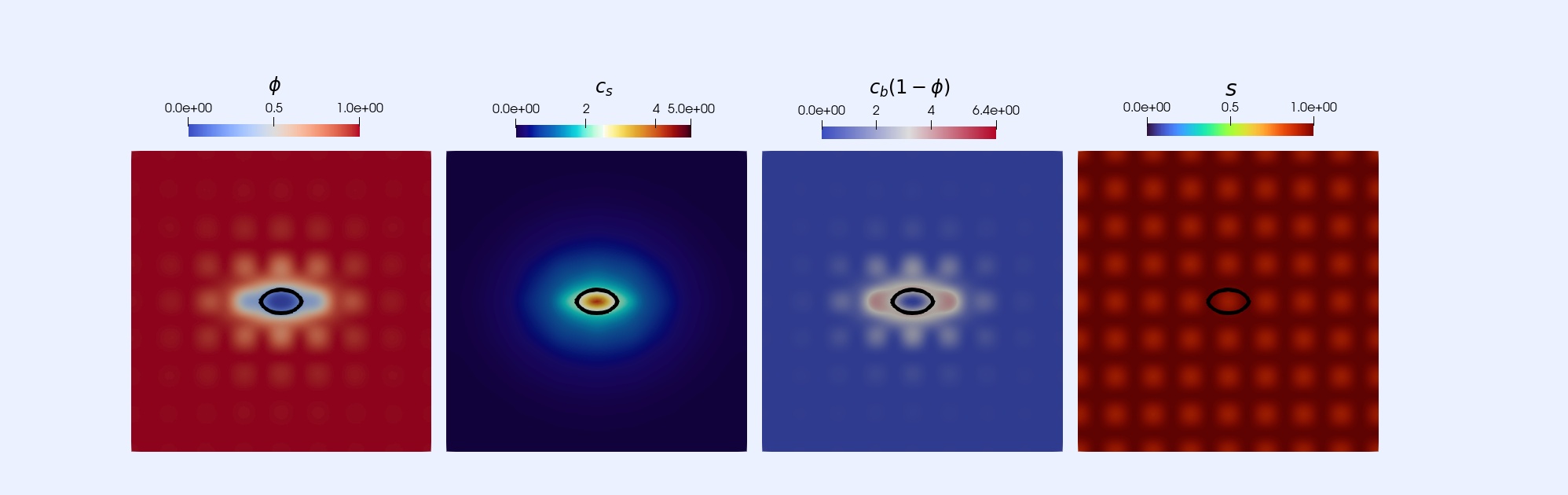}
\caption{$\mu_b=\delta_s=0,\mu_s=1$. Note for the results in this subfigure the bound MMPs do not influence $\phi$ and $s$ is constant as $\delta_s=0$.}
  \label{fig:2d_suitability_mus}
    \end{subfigure}
    \caption{Simulations of the heterogeneous ECM suitability invasion model~\eqref{het_ecm},~\eqref{IC_BC} in $2{\rm d}$ with initial conditions corresponding to low matrix suitability \eqref{ic_s_unsuitable}. In each subfigure, each row corresponds to times $t=1, 3$ and $5$. The black contour indicates the level set $\phi=0.25$ corresponding to a volume fraction of $75\%$ cancer cells. We note that, in both cases where the bound enzymes affect suitability,  spatial inhomogeneities arise in the invasive front. We also observe that the  invasiveness is greatly reduced in the case where the bound MMPs do not degrade the ECM or influence its suitability. Parameter values as in Table~\ref{table_1}.}
    \label{fig:2d_low_suitability}
  \end{figure}
Figure \ref{fig:2d_low_suitability} shows the corresponding simulation results. We note that in both cases where the bound MMPs affect the suitability (Figures \ref{fig:2d_suitability_both} and \ref{fig:2d_suitability_mub}) the ECM is made more suitable for invasion and then the ECM is degraded with  the matrix  suitability close to $0$ in the regions where $\phi<0.25$. On the other hand in the simulation results in Figure~\ref{fig:2d_suitability_mus} where the bound MMPs do not degrade the ECM or change the suitability (hence $s$ is constant in time) the invasive process is almost completely stalled with virtually no change in the level set  corresponding to $\phi=0.25$. As noted previously when discussing the simulations with a homogeneous ECM these results are consistent with the experimental observations of  \cite[Figure 1A]{sabeh2009protease}.

For the final set of simulation results  we consider a setup that is similar to that of \cite[Figure 7]{deakin2013mathematical} and take initial conditions of the form
\begin{equation}\label{IC_Deakin}
\begin{aligned}
\phi_0(x)=1-e^{-4(1+x_1)}\quad\text{and}\quad s_0(x)=\frac{1}{2}\left(1+\cos(4\pi x_1)\cos(4\pi x_2)\right),
\end{aligned}
\end{equation}
and set the initial conditions  for the MMPs to $0$. We note the initial conditions correspond to cancer cells invading from a region near the left boundary and a checkerboard initial  suitability with  suitable and unsuitable regions for invasion both present at the initial time.  \begin{figure}[htbp]
  \begin{subfigure}[t]{0.49\textwidth}
  \includegraphics[trim = 60mm 160mm 90mm 30mm,  clip, width=\linewidth]{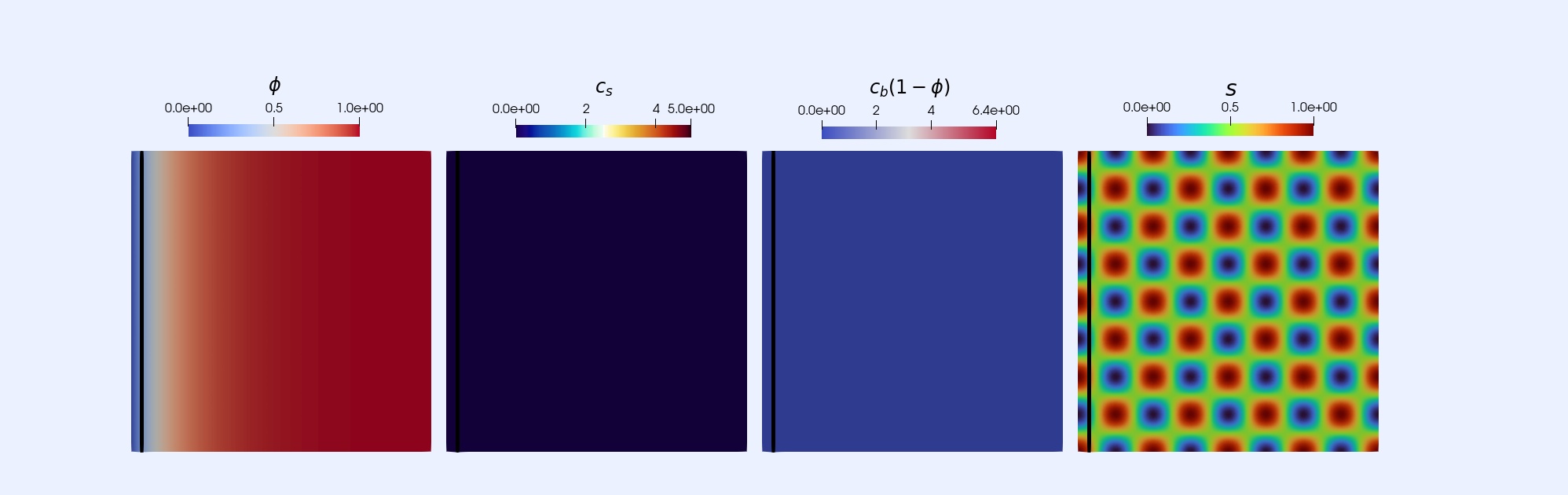}
  \includegraphics[trim = 60mm 0mm 90mm 65mm,  clip, width=\linewidth]{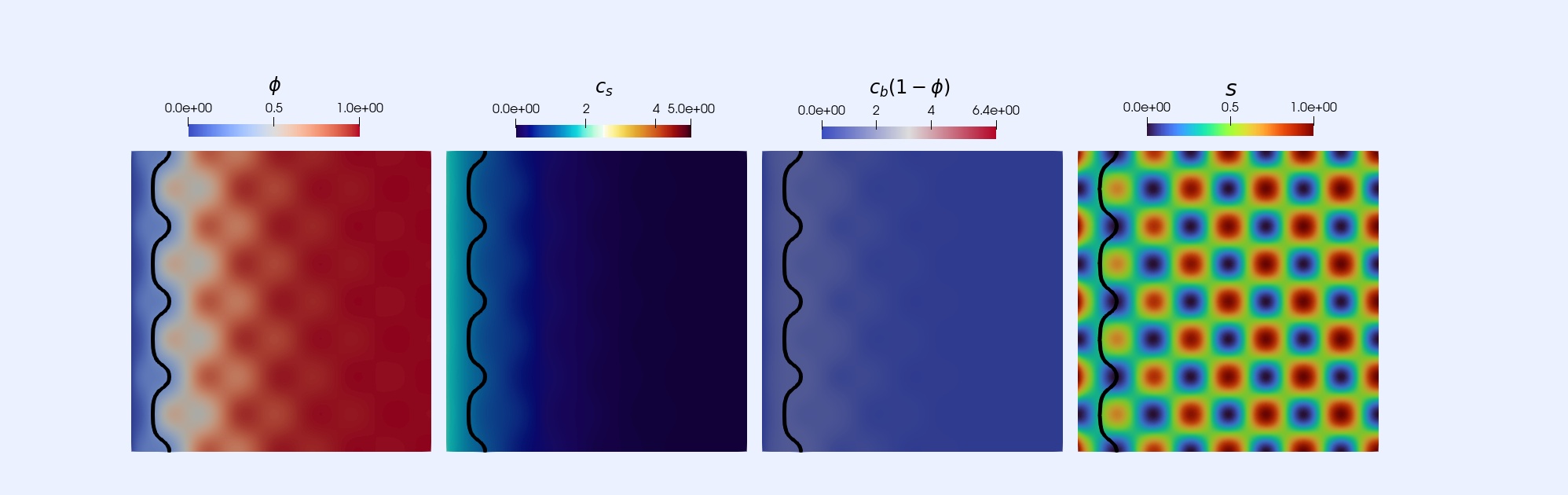}
  \includegraphics[trim = 60mm 0mm 90mm 65mm,  clip, width=\linewidth]{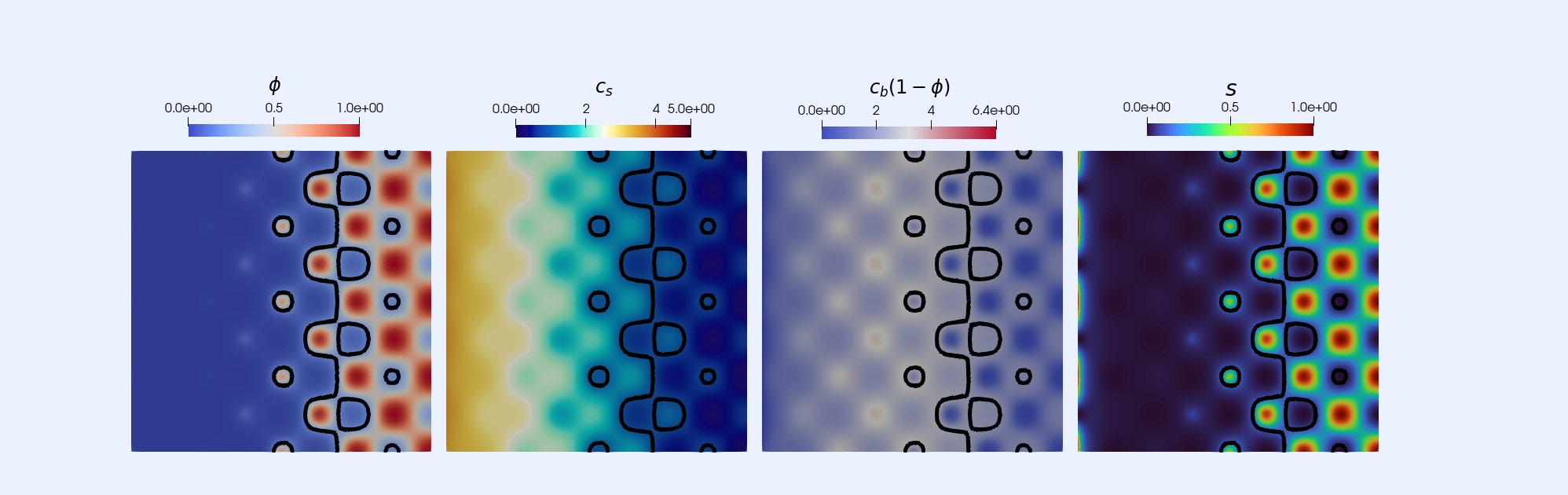}
  \includegraphics[trim = 60mm 0mm 90mm 65mm,  clip, width=\linewidth]{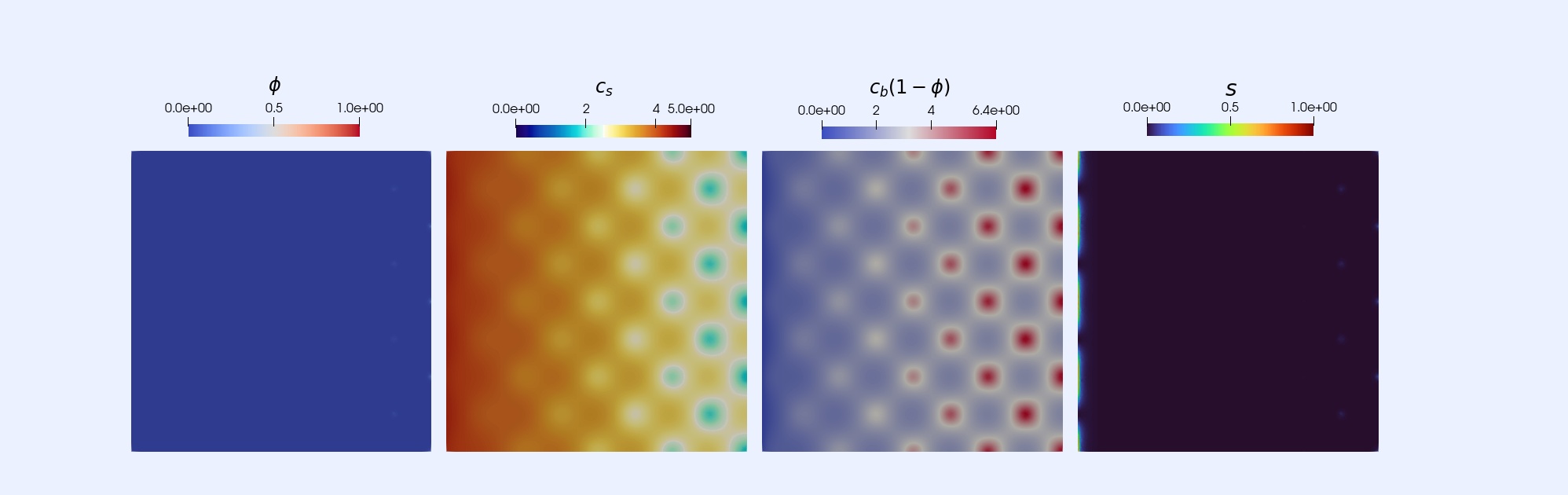}
\caption{$\mu_b=\mu_s=1$. Note for the results in this subfigure both soluble and bound MMPs degrade the matrix.}
  \label{fig:2d_Deakin_both}
    \end{subfigure}
    \hfill
  \begin{subfigure}[t]{0,49\textwidth}
  \includegraphics[trim = 60mm 160mm 90mm 30mm,  clip, width=\linewidth]{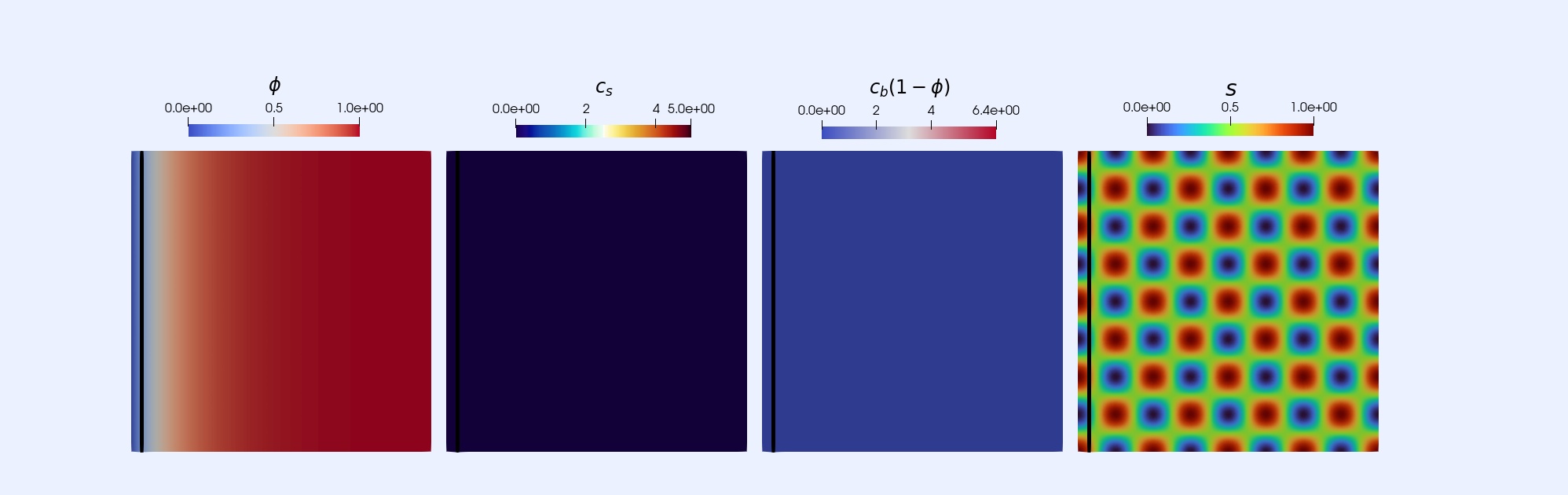}
  \includegraphics[trim = 60mm 0mm 90mm 65mm,  clip, width=\linewidth]{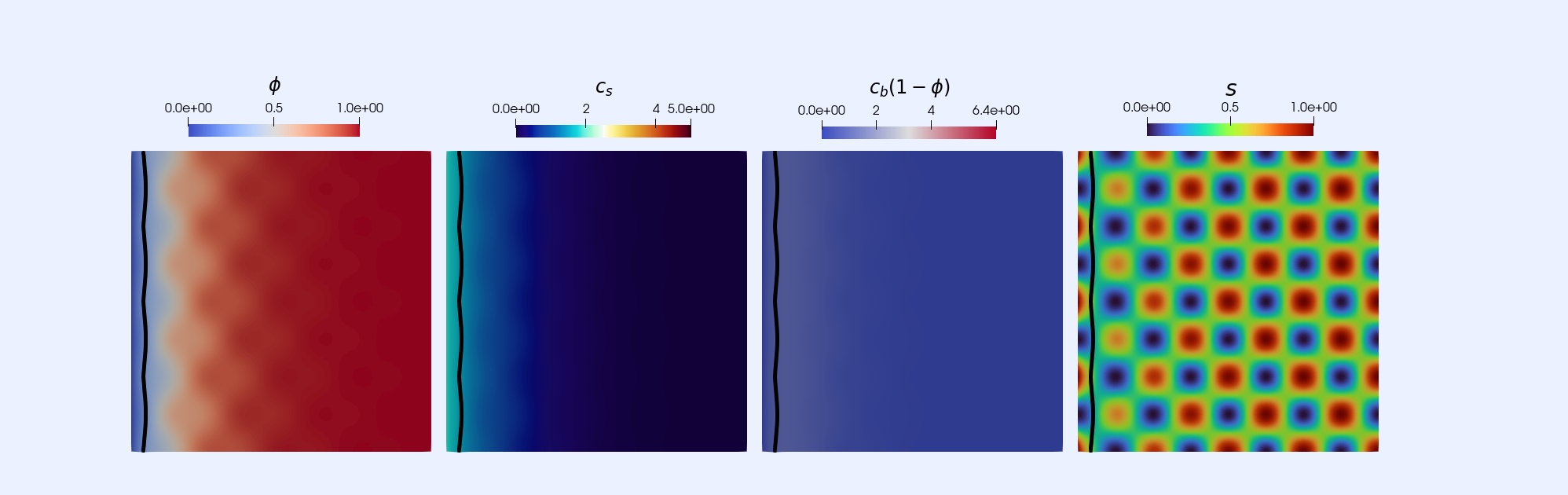}
  \includegraphics[trim = 60mm 0mm 90mm 65mm,  clip, width=\linewidth]{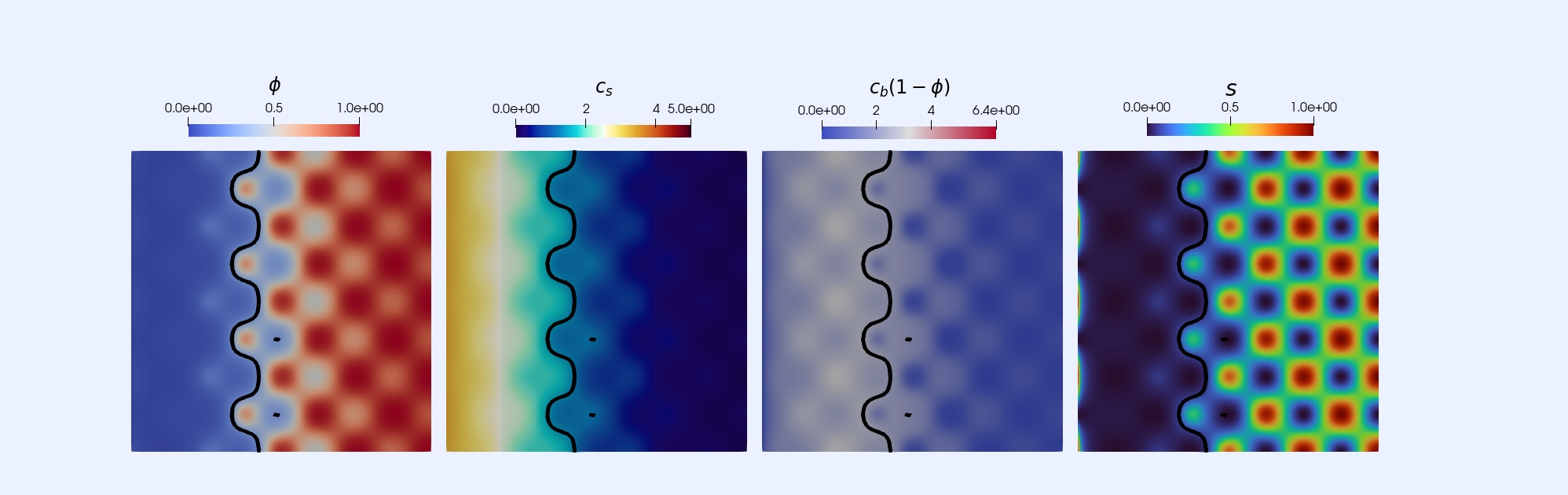}
  \includegraphics[trim = 60mm 0mm 90mm 65mm,  clip, width=\linewidth]{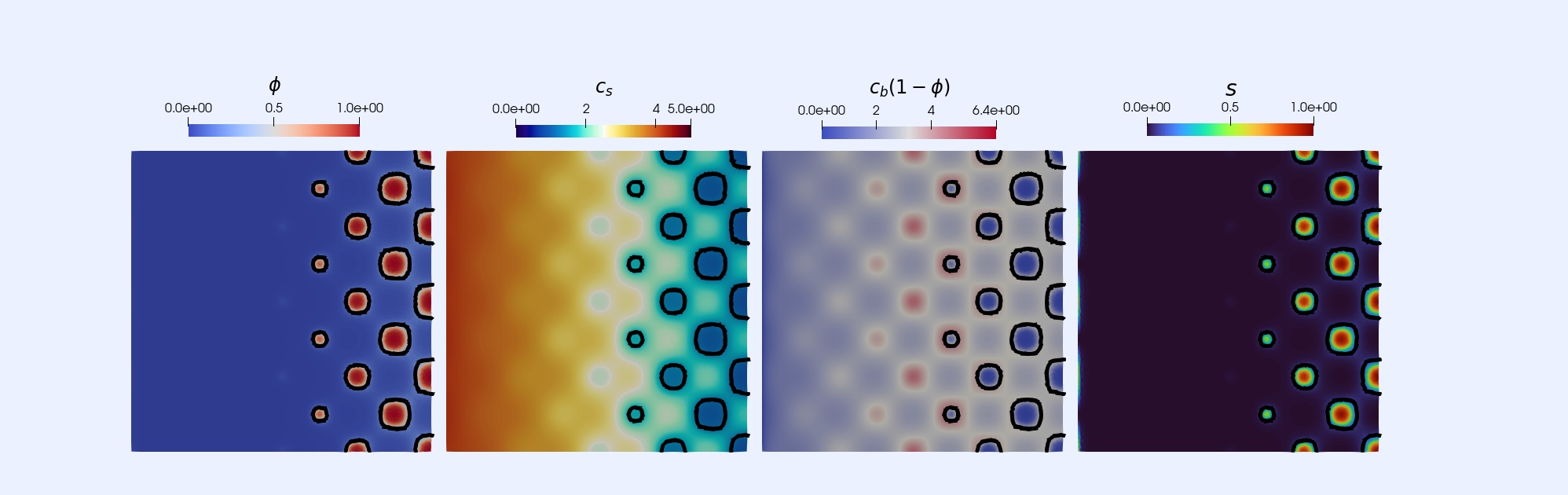}
\caption{$\mu_b=1,\mu_s=0$. Note for the results in this subfigure the soluble MMPs do not influence $\phi$.}
  \label{fig:2d_Deakin_mub}
    \end{subfigure}
    \hfill
  \begin{subfigure}[t]{0.49\textwidth}
  \includegraphics[trim = 60mm 160mm 90mm 30mm,  clip, width=\linewidth]{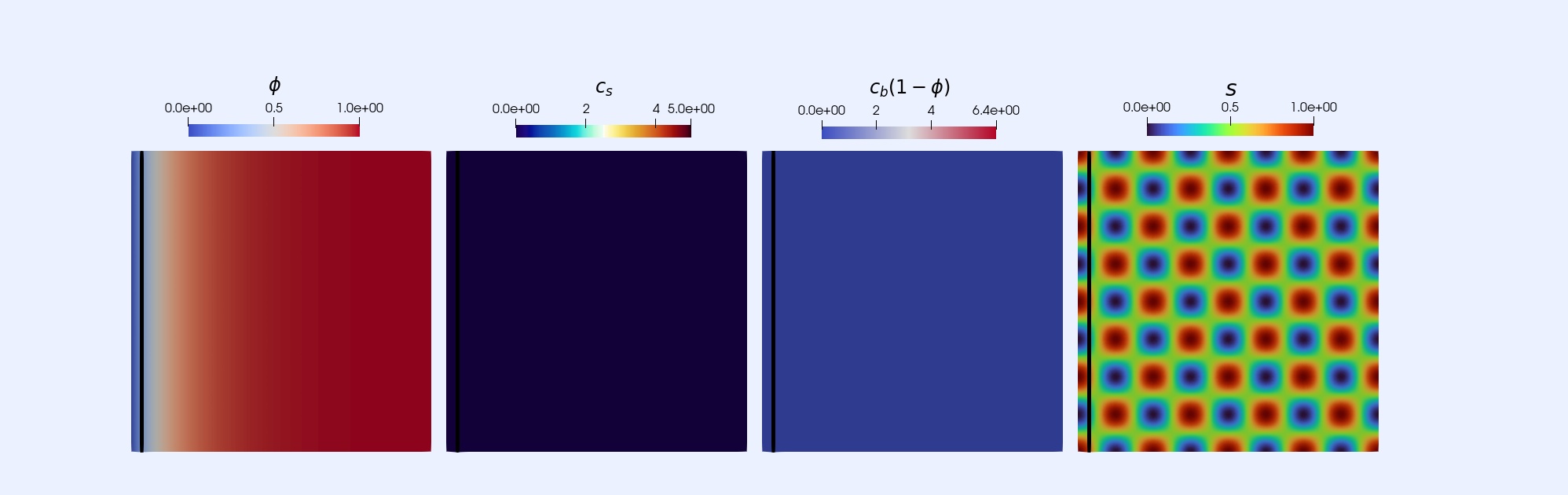}
  \includegraphics[trim = 60mm 0mm 90mm 65mm,  clip, width=\linewidth]{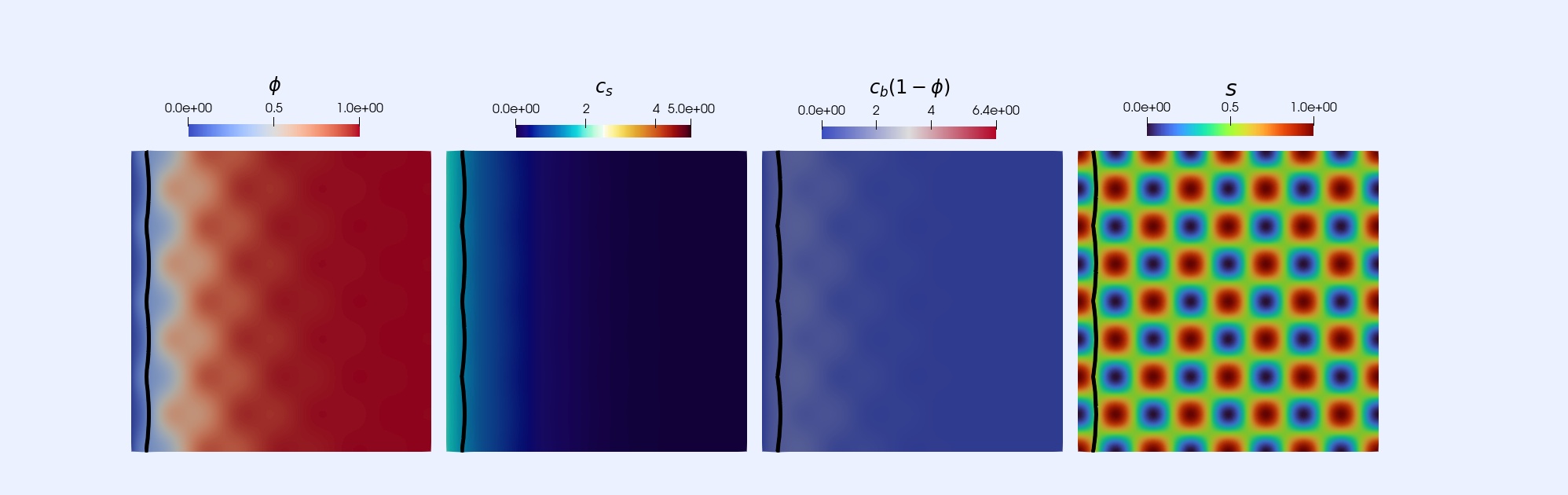}
  \includegraphics[trim = 60mm 0mm 90mm 65mm,  clip, width=\linewidth]{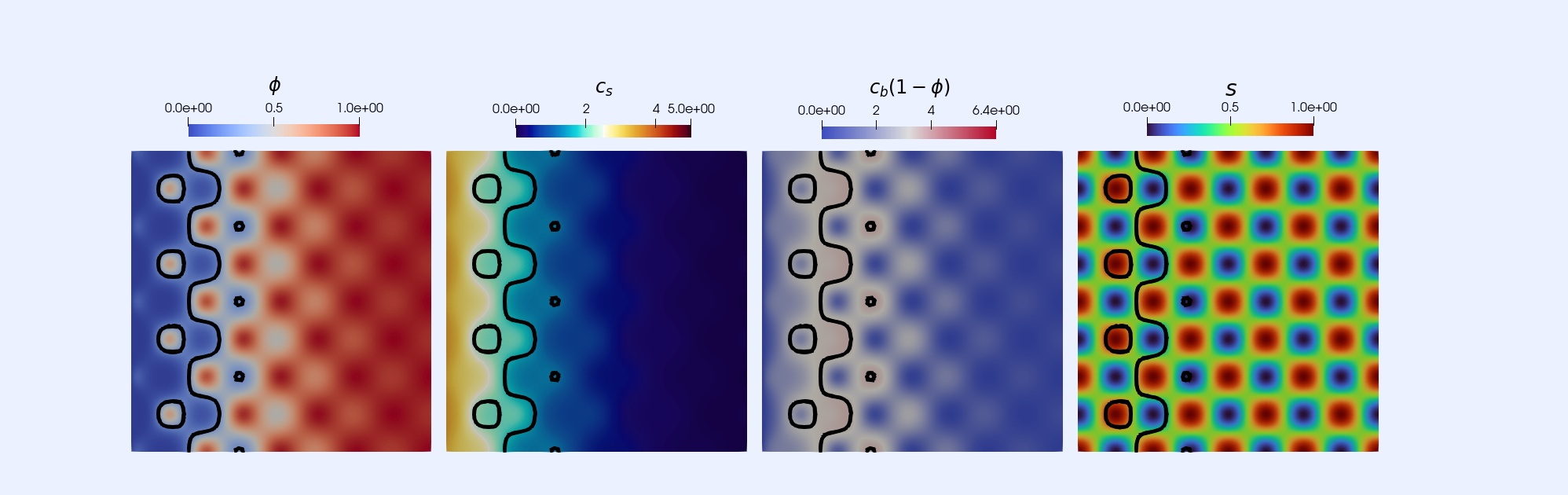}
  \includegraphics[trim = 60mm 0mm 90mm 65mm,  clip, width=\linewidth]{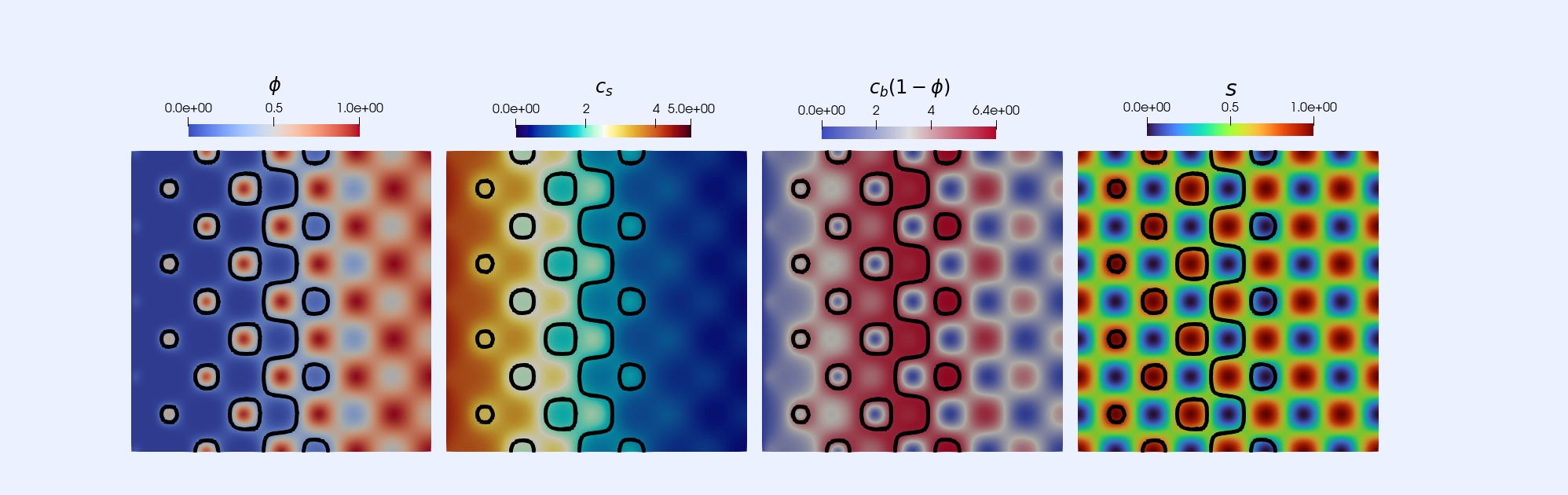}
\caption{$\mu_b=\delta_s=0,\mu_s=1$. Note for the results in this subfigure the bound MMPs do not influence $\phi$ and $s$ is constant as $\delta_s=0$.}
  \label{fig:2d_Deakin_mus}
    \end{subfigure}
    \caption{Simulations of the heterogeneous ECM suitability invasion model \eqref{het_ecm}, \eqref{IC_BC} in $2{\rm d}$ to be compared with \cite[Figure 7]{deakin2013mathematical}. In each subfigure, each row corresponds to times $t= 1, 3$ and $5$. The black contour indicates the level set $\phi=0.25$ corresponding to a volume fraction of $75\%$ cancer cells. The role of the bound MMPs in making the matrix more suitable for invasion is clear as is the reduced invasiveness in the case where the bound MMPs do not degrade the ECM or influence its suitability. Parameter values as in Table~\ref{table_1}.}
    \label{fig:2d_Deakin}
  \end{figure}
The results of the simulation are reported on in Figure~\ref{fig:2d_Deakin} and are broadly similar to those shown in~\cite[Figure 7]{deakin2013mathematical} and exhibit similar qualitative features to those reported on in Figure~\ref{fig:2d_low_suitability}. We observe a heterogeneous invasion front in all cases that largely matches the spatial variations in the suitability at the initial time.  In the simulations where the bound MMPs are active (Figures \ref{fig:2d_Deakin_both} and \ref{fig:2d_Deakin_mub}) we  again observe that the suitability is close to zero in regions where $\phi<0.25$. We also note that despite the presence of regions which are suitable for invasion at the initial time, the invasive speed, in the case where the bound MMPs are not active, is reduced  (Figure \ref{fig:2d_Deakin_mus}), albeit not as much as was the case in the simulations of Figure~\ref{fig:2d_low_suitability}.
 
\section{Discussion and possible extensions}\label{sec:conc}

In this work, we have presented a multiscale framework for the modelling of cancer cell invasion into the ECM through the degradation of the ECM by bound and soluble matrix-degrading enzymes. Our framework is novel in that we give a rational way of encoding microscopic details in a macroscopic (tissue-level) model, which accounts for the computation of an effective diffusivity that depends on the volume fraction of the ECM  as well as a means of incorporating details of cell-scale signalling processes in a tissue-level framework. We report on numerical simulations for the macroscopic model, which illustrate the role bound and soluble enzymes play in the invasion process and highlight interesting features such as the propensity of degradation of the ECM by membrane-bound enzymes to generate more spatially heterogeneous invasion profiles. The latter could be of importance in understanding the onset of complex tumour morphologies, which are associated with metastasis. 

One main area for future work is the rigorous derivation of a macroscopic invasion model from microscopic models such as those proposed in this work, 
which is of significant interest beyond this specific application. We believe this should be possible in the setting where cell movement drives invasion but proliferation is neglected, using the ideas from the multiscale analysis of free boundary problems, e.g.~\cite{kim2010, pozar2015, rodrigues1982, visintin2007}. As mentioned briefly previously, incorporating proliferation even at the microscopic level appears to require the development of a novel theory of continuum PDE models beyond topological changes, and this needs to be clarified before the derivation of macroscopic models involving proliferation can be considered.    

Several other extensions of this study are important to develop a detailed description and a better understanding of cancer invasion. From the modelling perspective, the two-scale framework warrants further study, especially if a detailed cell-signalling model such as those proposed in~\cite{ptashnyk2020multiscale} with non-trivial dynamics is considered. The structure of the ECM can be represented in greater detail, such as accounting for the orientation of the collagen fibres, which influence cell migration as well as a number of other factors, see \cite{crossley2024modeling} for a review. Depending on the timescales, ECM remodelling could also be included. The mixture-model-based formalism we employ could also be extended to allow for necrotic regions within the tumour. From the perspective of analysis, an obvious extension would be to extend the proof of the well-posedness to the degenerate case where we do not assume a strictly positive lower bound for the initial condition for the volume fraction of ECM. However, this appears challenging and requires the development of novel analytical techniques.  
\edit{From the applications perspective, the results above complement the work of \cite{deakin2013mathematical} by allowing for greater connection to microscopic cell scale models of matrix degradation. The work of \cite{deakin2013mathematical} was seeking to model the experiments of \cite{sabeh2009protease}.  Further numerical experiments in which we  simulate more closely the scenarios considered in the experiments of \cite{sabeh2009protease} and other related studies will allow us to assess the ability of the mathematical model to reproduce biological observations with the eventual goal of augmenting and guiding biological experiments through computational simulations.} 

\section*{Data Availability Statement}
Data sharing is not applicable to this article as no new data were created or analysed in this study.

\section*{Acknowledgments}
The authors wish to thank Mark Chaplain for many fruitful discussions both during the preparation of this manuscript and concerning the modelling of biological phenomena in general.

CV wishes to thank the ICMS, Edinburgh for support towards participation in the programme '3MC+PIMS+ICMS Winter school - Multiscale Modeling: Infectious Diseases, Cancer and Treatments' in which aspects of this work were finalised.

This work (CV) was partially supported by an individual grant from the Dr Perry James (Jim) Browne Research Centre on Mathematics and its Applications (University of Sussex).

Both authors wish to thank the anonymous reviewers for their insightful comments which improved this manuscript.
\bibliographystyle{alpha}      

\bibliography{./CPV_MMP.bib}
\edit{
\appendix
\section{Additional computations of effective diffusion tensors}\label{app:eff_diff}
We report on computations of the effective diffusion tensor considering different microstructures.
\subsection{Alternative cell shape}\label{app:square_eff_diff}
In order to consider different cell shapes and allow smaller volume fractions of ECM to be modelled without issues arising due to meshing complex geometries we  consider the following geometric setup for the microstructure. The domain corresponds to a region of the form $[-2,2]^d\setminus[-Y,Y]^d$ where $d=2,3$ and $Y=0.35n, n=1,2,3,4,5$, i.e., the domain is a hypercube of side length $4$ with  a single hypercube  of side length $0.7n$ removed (i.e., the cell is a hypercube)  both centred at  the origin.
 The corresponding volume fractions of ECM are $0.97, 0.88, 0.72, 0.51$ and $0.23$ in $2{\rm d}$ and $0.99, 0.96, 0.86, 0.66$ and $0.33$ in $3{\rm d}$ for $n=1,2,3,4,5$ and $6$ respectively.  We approximate the solution to \eqref{eqn:cell_problems} using a finite element method implemented in the FEniCS software~\cite{logg2012automated}. We take $D_s(0)=0$ as a further fitting point as we assume the soluble MMPs do not diffuse through the cancer cells.  Fitting a cubic polynomial  in the $2{\rm d}$ and $3{\rm d}$ case (to the average of the computed diffusivities, noting that the off-diagonal components are zero)  yields the following expressions
 \begin{align}
 {D_s(\phi)}=\big(0.67\phi^3-0.27\phi^2+0.60\phi\big){D_s(1)},
 \end{align}
 in $2{\rm d}$ and
  \begin{align}
 {D_s(\phi)}=\big(1.32\phi^3-2.94\phi^2+2.62\phi\big){D_s(1)},
 \end{align}
 in $3{\rm d}$. Figure \ref{fig:square_d_eff} shows plots of the effective diffusivity as the volume fraction of ECM changes. 
\begin{figure}[htbp]
  \includegraphics[trim = 0mm 0mm 0mm 0mm,  clip, width=0.45\linewidth]{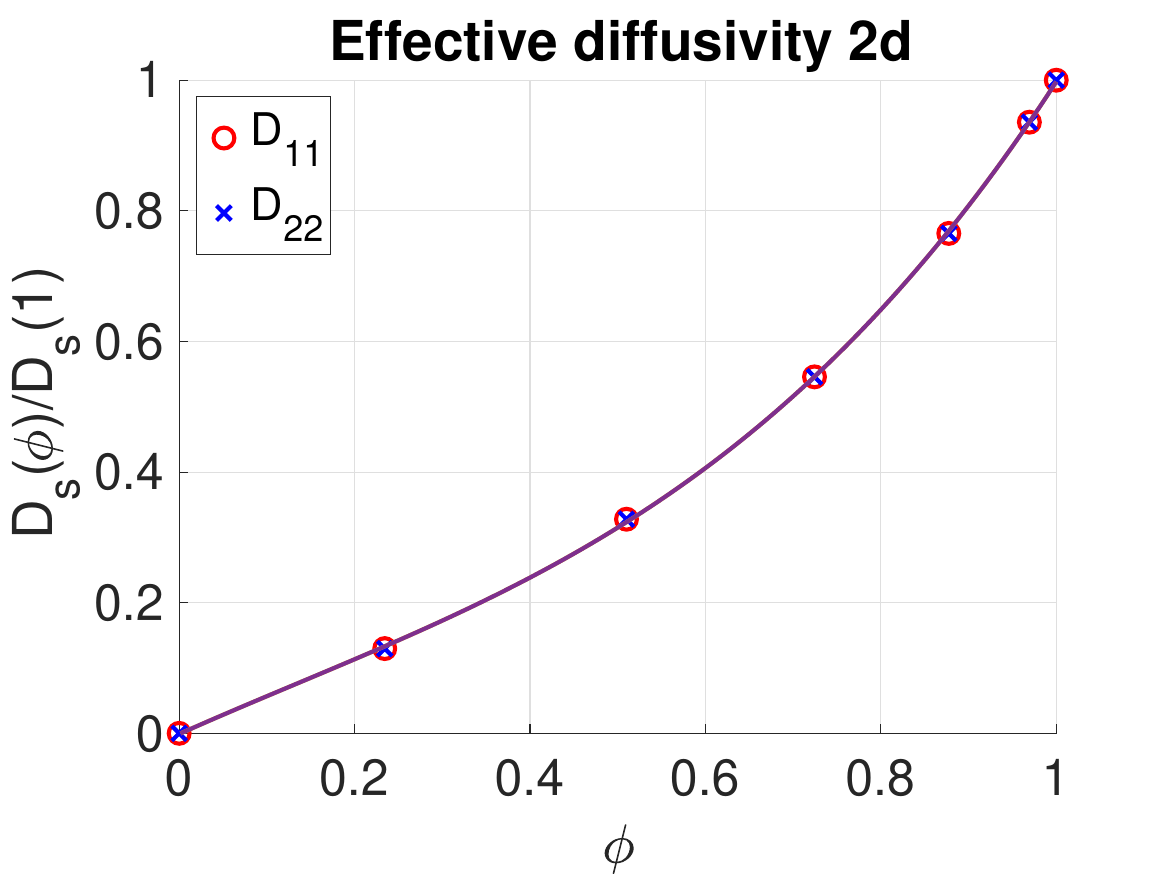}
  \includegraphics[trim = 0mm 0mm 0mm 0mm,  clip, width=0.45\linewidth]{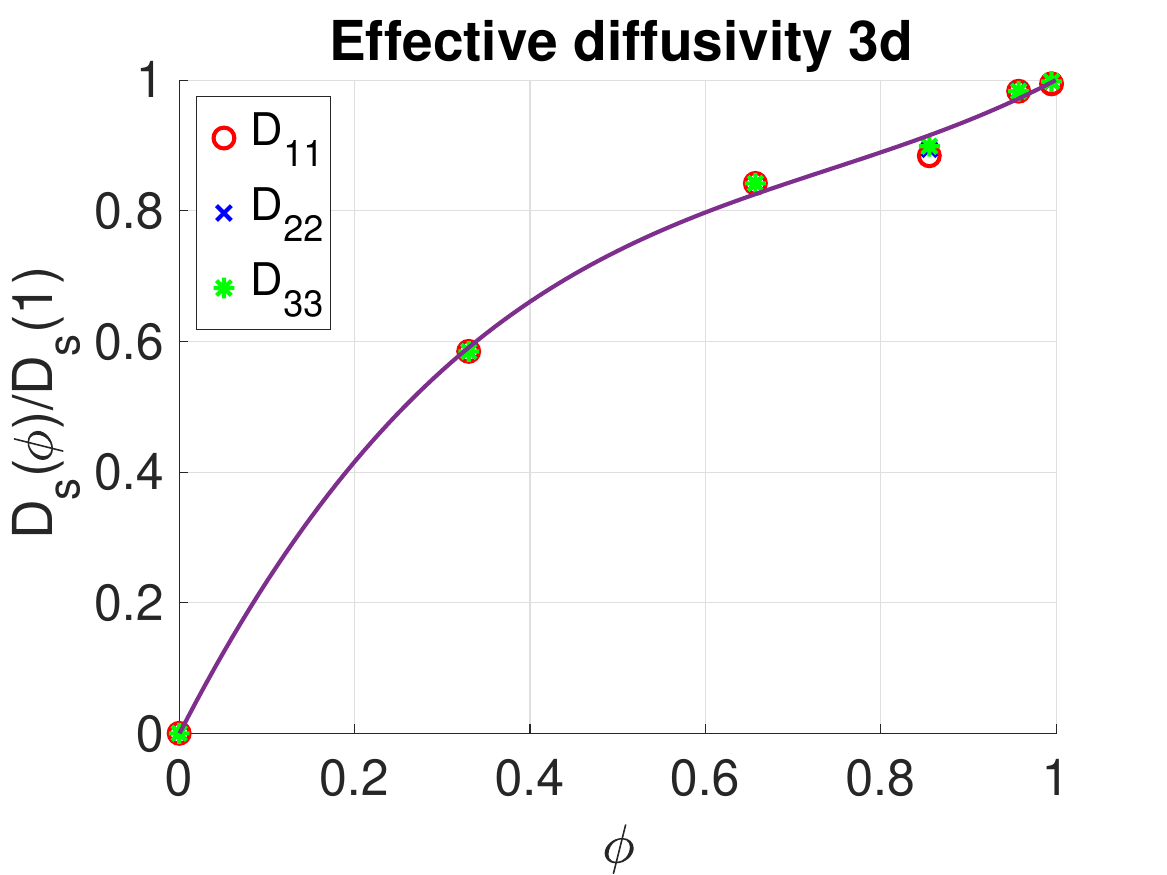}
\caption{Fitted effective diffusivity in $2{\rm d}$ and $3{\rm d}$ for a microstructure consisting of a single square (cube in $3{\rm d}$) shaped  cell of varying length. The red circles and blue crosses correspond to computed effective diffusivities obtained by solving~\eqref{eqn:cell_problems} on domains described in section~\ref{app:square_eff_diff}. The curves are  least squares cubic polynomial fits using the average of the estimated diffusivities.}
  \label{fig:square_d_eff}
  \end{figure}
  We note that in the relevant range of volume fractions, $\phi\in[0,1]$, the fitted polynomial for the effective diffusivity is very close to the previously obtained polynomials, see section~\ref{sec:homg_diff_sims}, as shown in Figure~\ref{fig:d_eff_comp}.
\begin{figure}[htbp]
  \includegraphics[trim = 0mm 0mm 0mm 0mm,  clip, width=0.45\linewidth]{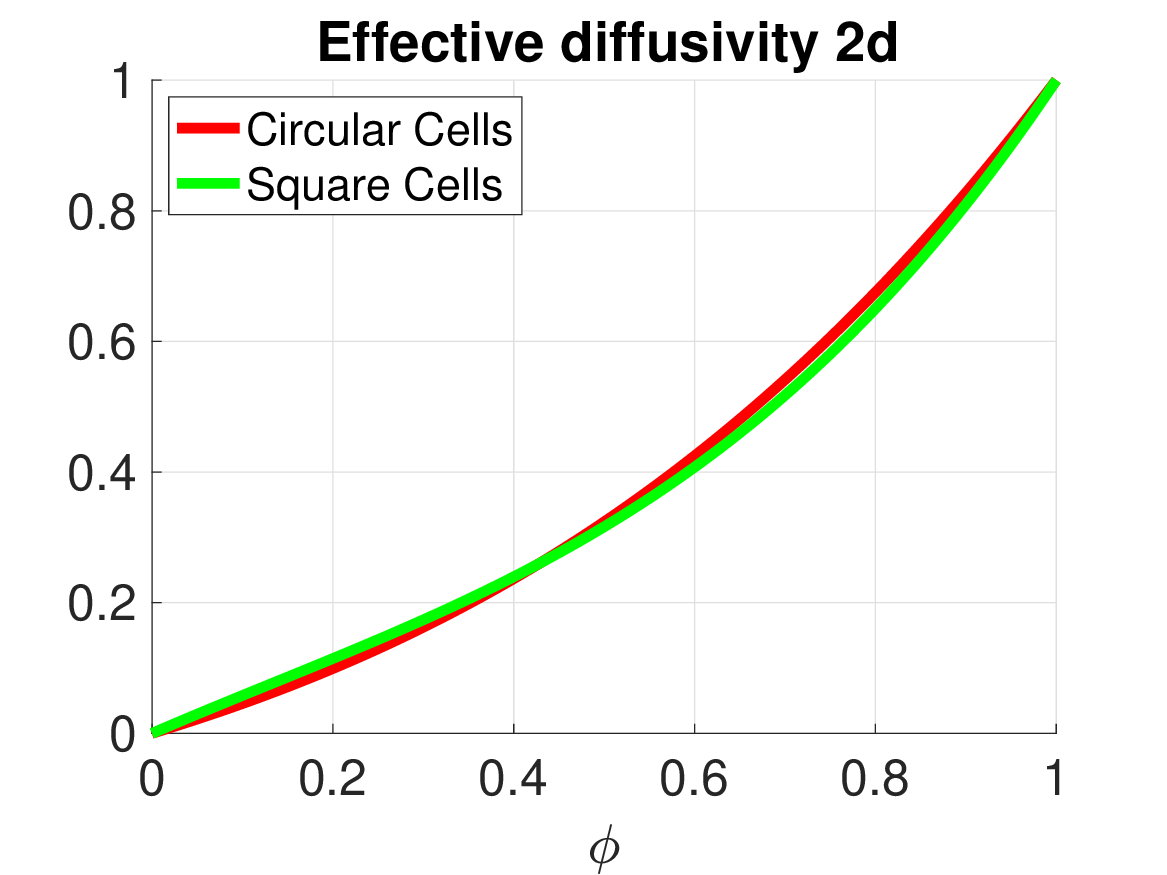}
\hskip 0.1em
  \includegraphics[trim = 0mm 0mm 0mm 0mm,  clip, width=0.45\linewidth]{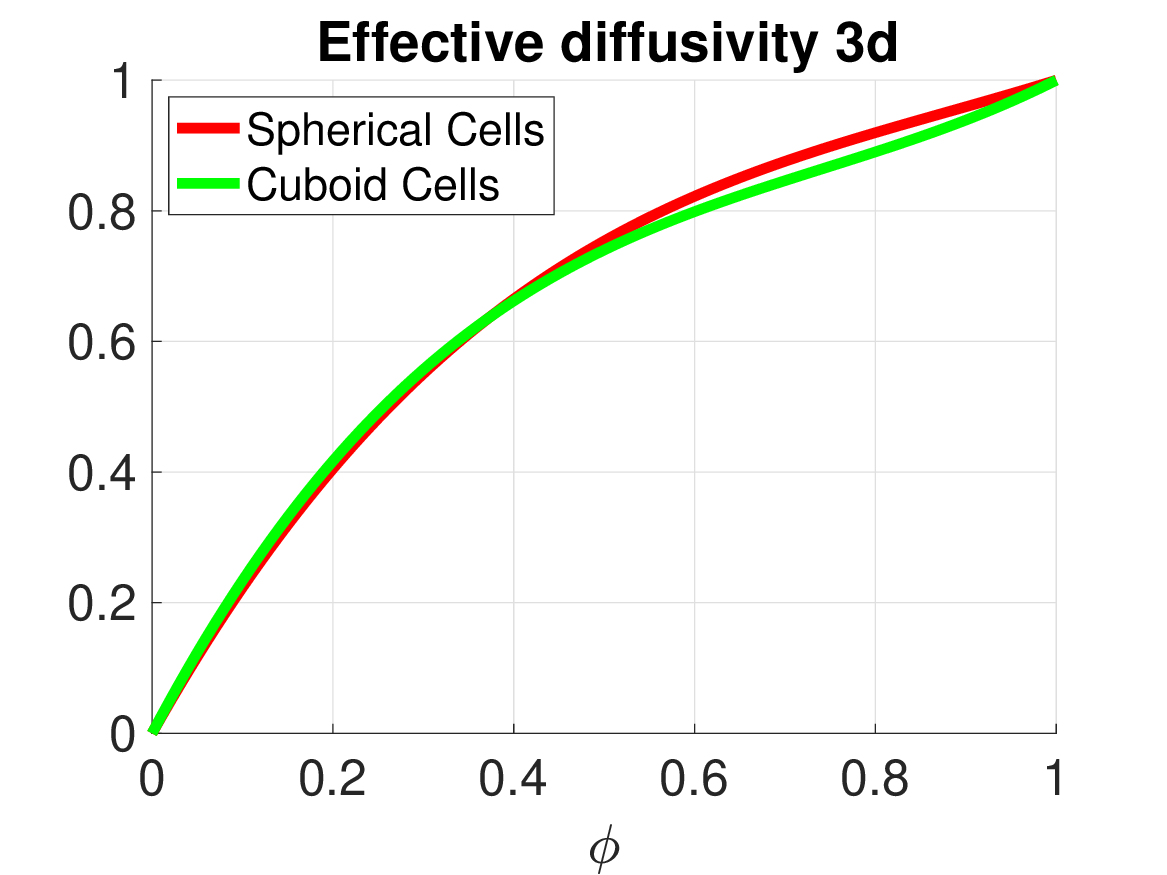}
\caption{Comparison of the polynomials describing the effective diffusivity for the microstructures considered in sections~\ref{sec:homg_diff_sims} and \ref{app:square_eff_diff}.}
  \label{fig:d_eff_comp}
  \end{figure}
  \subsection{Microstructure that results in a full effective diffusion tensor}\label{app:full_eff_diff}
In this section we consider a microstructure that results in a full effective diffusion tensor, specifically  the domains we use correspond to squares of the form $[0,4]^2$ with $n^2$, for $n=1,2,3,4$,  ellipse(s) of major and minor axes $0.6$ and $0.25$ removed. In order to obtain a full effective diffusion tensor the orientation of the ellipses is such that the major axis of each ellipse  is parallel to the  line $y=x$. The centres of the ellipses are taken to form a square lattice, with the resulting geometry for a single ellipse as shown in Figure~\ref{fig:hom_ellipse_2d}. The corresponding volume fractions of ECM are $0.970,0.880,0.729$ and $0.518$  for $n=1,2,3$ and $4$ respectively. 
  \begin{figure}[htbp]
  \includegraphics[trim = 70mm 0mm 70mm 0mm,  clip, width=0.25\linewidth]{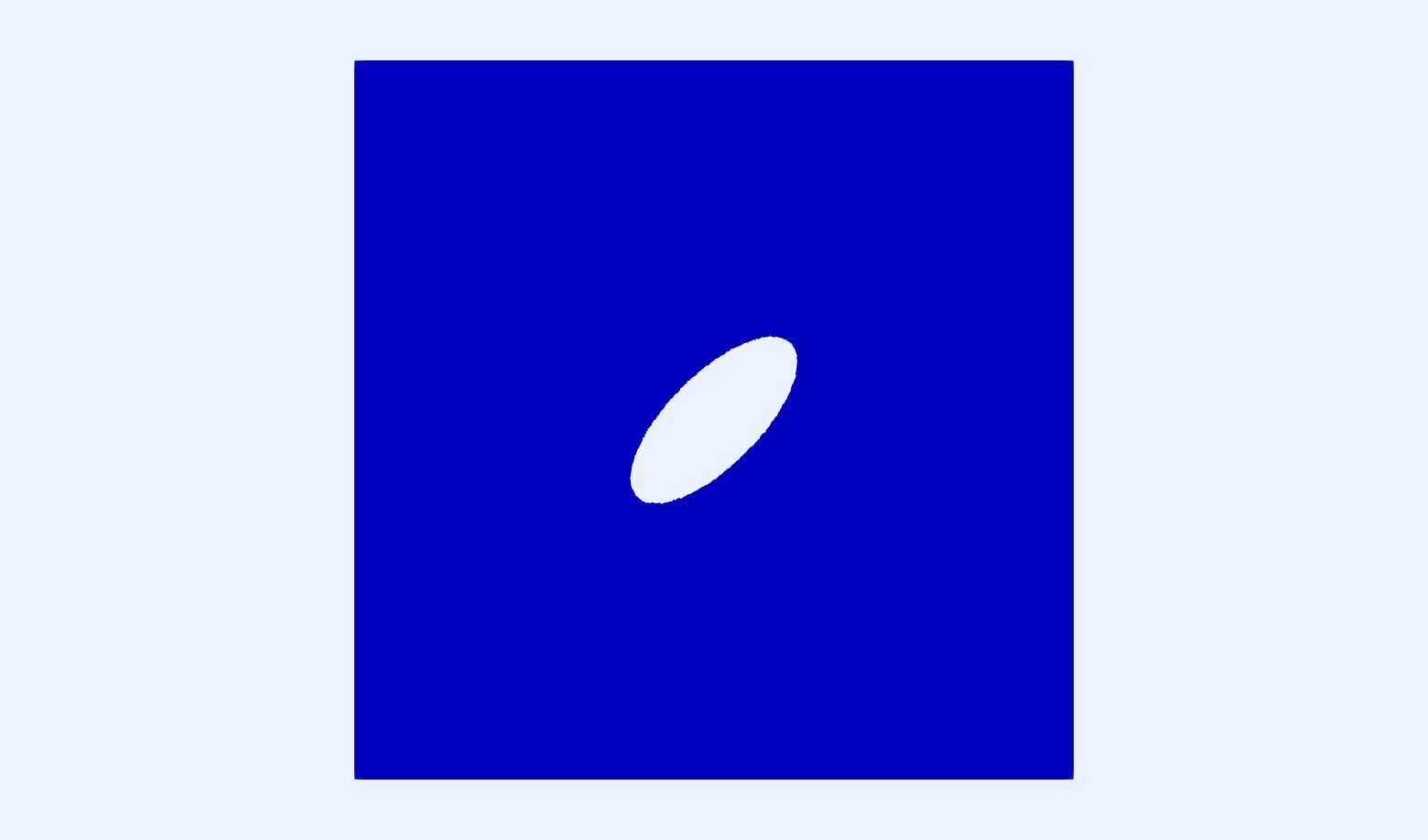}
\caption{An example of a microscopic geometry that leads to a full diffusion tensor. }
  \label{fig:hom_ellipse_2d}
  \end{figure}
 For the choice of microscopic geometry considered here  we note that,  the effective diffusion tensor is not expected to be  diagonal.  This is reflected in the numerical results, for which an unstructured triangulation is used, in that the off diagonal components of the computed diffusion tensors for all geometries are of the same order of magnitude as the diagonal components and the diagonal components of the computed diffusion tensors on each specific geometry exhibit only small differences ($<1\%$). To obtain the effective diffusion tensor as a function of
  $\phi$ we note that a constraint we must enforce is that the diffusion tensor remains positive definite for all $\phi$ strictly positive. In higher spatial dimensions naive polynomial fits of the individual coefficients or polynomial interpolation do not guarantee this. We instead follow \cite{arsigny2006log} and use Log-Euclidean interpolation to interpolate the tensor and also take ${\vec D_s}_{i,j}(0)=0$, $i,j=1,2$ as a further fitting point since there is no diffusion of MMPs within the cells. As  
 $\vec D_s(0)$ is not strictly positive definite  we use linear interpolation in the interval with $\phi=0$ as the left hand end point. For completeness we state the Log-Euclidean interpolation formula below and refer to \cite{arsigny2006log} for theoretical considerations. Given symmetric positive definite diffusion tensors $\vec D(\phi_1)$ and $\vec D(\phi_2)$ assuming without loss of generality that $0<\phi_1<\phi_2\leq 1$ we set
 \[
 \vec L_1=\operatorname{logm}(\vec D(\phi_1))\quad\text{and}\quad  \vec L_2=\operatorname{logm}(\vec D(\phi_2))
 \]
 and for $\phi\in[\phi_1,\phi_2]$ we set
 \[
 \vec L(\phi)=\operatorname{expm}\left(\frac{(\phi_2-\phi)}{(\phi_2-\phi_1)}\vec L_1+\frac{(\phi-\phi_1)}{(\phi_2-\phi_1)}\vec L_2\right),
 \] 
 where $\operatorname{logm}$ and $\operatorname{expm}$ denote the matrix logarithm and matrix exponential respectively.
 Figure~\ref{fig:d_eff_full} shows plots of the effective diffusivity as the volume fraction of ECM changes together with the  interpolated values (Log-Euclidean for $\phi >0.518$ and linear for $\phi<0.518$).
\begin{figure}[htbp]
  \includegraphics[trim = 0mm 0mm 0mm 0mm,  clip, width=0.45\linewidth]{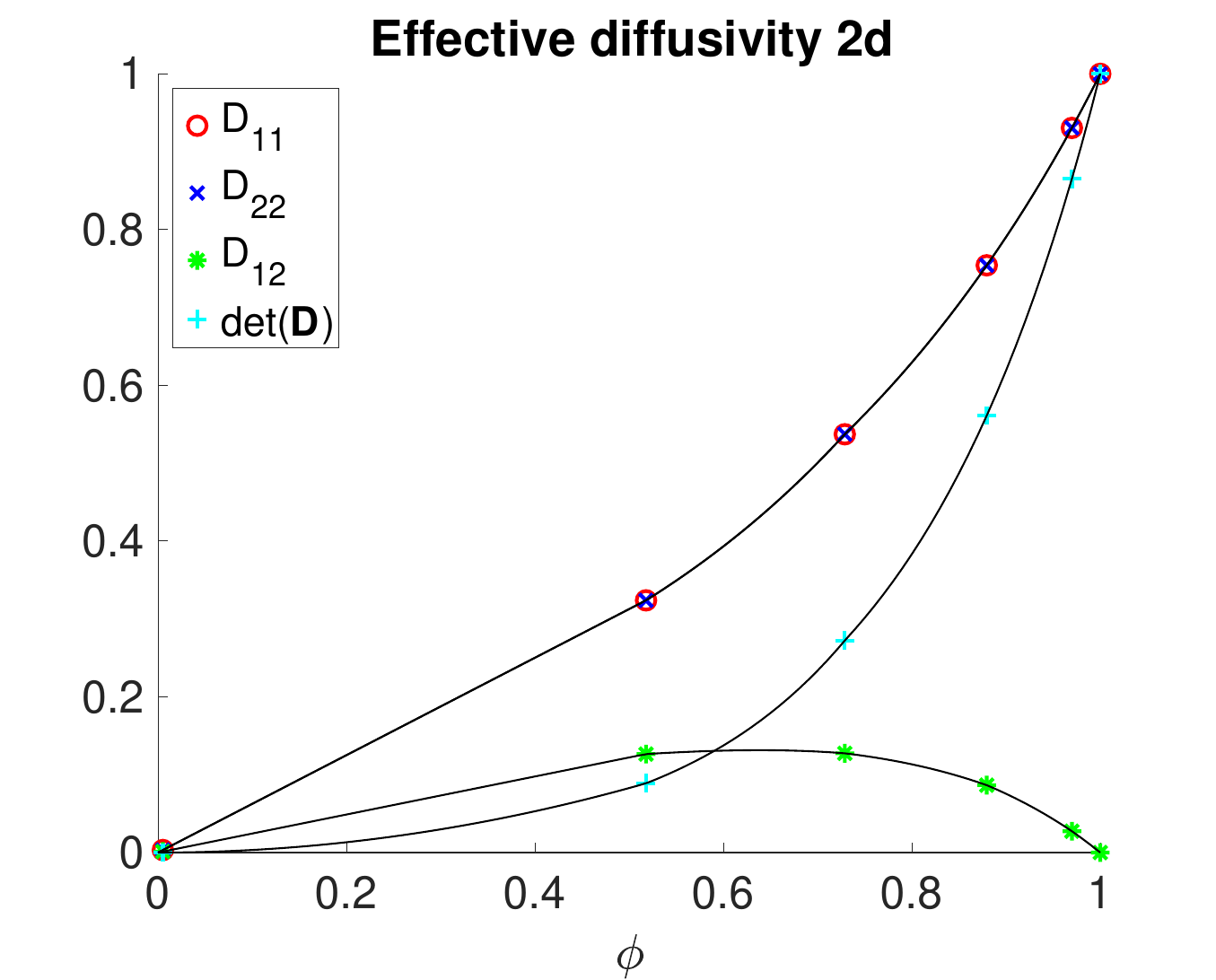}
\caption{Fitted effective diffusivity in $2{\rm d}$ for the example in section~\ref{app:full_eff_diff}. The red circles, blue crosses and green asterisks correspond to computed effective diffusivities obtained by solving~\eqref{eqn:cell_problems} on domains described in section~\ref{app:full_eff_diff}. The  curves for the diffusivities are  interpolated values (Log-Euclidean for $\phi >0.518$, linear for $\phi<0.518$). The determinant is computed directly using the interpolated diffusivity values.}
  \label{fig:d_eff_full}
  \end{figure}
We note that the diagonal components of the effective diffusivity appear to be a convex functions of the volume fraction whilst the off diagonal components appear to be concave (in $2$d).  
 }

\end{document}